\documentclass[12pt,letterpaper,leqno]{amsart}
\usepackage{microtype}
\usepackage{amssymb}
\usepackage{amsfonts}
\usepackage{geometry}
\usepackage{hyperref}

\renewcommand{\Re}{\operatorname{Re}}

\newtheorem{theorem}{Theorem}
\theoremstyle{plain}

\newtheorem{corollary}[theorem]{Corollary}

\newtheorem{lemma}[theorem]{Lemma}

\newtheorem{proposition}[theorem]{Proposition}

\numberwithin{equation}{section}
\numberwithin{theorem}{section}  
\IfFileExists{tcilatex.tex}{\input{tcilatex}}{}
\usepackage{tikz,longtable,booktabs,array,needspace,placeins}
\usetikzlibrary{arrows.meta,calc}
\numberwithin{figure}{section}
\newcommand{\R}{\mathbb R}
\newcommand{\Sph}{\mathbb S}
\newcommand{\dd}{\,\mathrm d}

\newcolumntype{L}[1]{>{\raggedright\arraybackslash}p{#1}}

\begin{document}

\title[Failure of fixed-profile soliton resolution]{A Counterexample to Fixed-Profile Soliton Resolution for the Three-Dimensional Quintic Wave Equation}
\author[X. Chen]{Xuwen Chen}
\address{Department of Mathematics, University of Rochester, Rochester, NY 14627, USA}
\email{xuwenmath@gmail.com}

\author[J. Holmer]{Justin Holmer}
\address{Department of Mathematics, Brown University, Providence, RI 02912, USA}
\email{justin\_holmer@brown.edu}

\author[S. Shen]{Shunlin Shen}
\address{School of Mathematical Sciences, University of Science and Technology of China, Hefei, 230026, China}
\email{slshen@ustc.edu.cn}

\date{}
\begin{abstract}
We construct a bounded forward-global solution of the three-dimensional focusing energy-critical wave equation whose single sign-changing core has two subsequential limits in distinct similarity orbits within one conformal class. This disproves the fixed-profile soliton resolution conjecture and suggests a reformulation in which each soliton has a fixed velocity and a profile that may vary within a fixed conformal class. Under full conformal nondegeneracy, we establish without symmetry assumptions that a one-core decomposition along a single sequence of times extends to all late times with the same class, velocity and free wave. In the counterexample, quadratic impulses from radiation packets with summable amplitudes produce fixed changes of shape over increasingly long intervals, and an energy-space norm identity rules out every fixed list of profiles. To prove the one-core resolution result, we identify the canonical radiation, use a linearized exterior energy estimate to detect growing modes, and derive decay of the error transverse to the orbit from the modulation equations.
\end{abstract}
\maketitle
\tableofcontents

\section{Introduction}\label{main:start}
We consider the three-dimensional focusing energy-critical wave equation
\begin{equation}\label{main:nlw}
 \partial_t^2u-\Delta u=u^5,\qquad (t,x)\in[0,\infty)\times\R^3 .
\end{equation}
We construct a bounded forward-global solution with a single stationary
core.\footnote{We use \emph{core} for a component of the profile
decomposition that we follow between observation times. In the radial
setting, one usually calls such a component a bubble.}
For $\vec u=(u,\partial_tu)$, we use the energy norm and conserved energy
\begin{align}
 \mathcal H&=\dot H^1(\R^3)\times L^2(\R^3),\qquad
 \|(f,g)\|_{\mathcal H}^2=\|\nabla f\|_2^2+\|g\|_2^2, \label{main:energy-space}\\
 E(f,g)&=\frac{1}{2}\|(f,g)\|_{\mathcal H}^2-\frac{1}{6}\|f\|_6^6 .
       \label{main:energy}
\end{align}
The conformal group $G=\operatorname{Conf}(\Sph^3)$ acts on stationary
states by stereographic projection. We denote by $\mathrm{Sim}$ the group
of spatial similarities, including reflections. For the solution we
construct, there is a free wave $v_+$ such that
\begin{equation}\label{main:intro-objective}
 \begin{gathered}
 \vec u(t)=(Q_{p(t)},0)+\vec v_+(t)+o_{\mathcal H}(1),\qquad Q_{p(t)}\in GQ,\\
 Q_{p(T_{2k})}\longrightarrow Q,\qquad
 Q_{p(T_{2k+1})}\longrightarrow Q_{(de,0)}\in GQ\setminus(\mathrm{Sim}\cdot Q)
 \quad\text{in }\dot H^1.
 \end{gathered}
\end{equation}

\subsection{Fixed profiles and the counterexample}\label{main:fixed-profiles}
In the conjecture of Duyckaerts, Kenig and Merle
\cite[Conjecture~5.5]{DKM}, one fixes the profiles, velocities and free
wave before choosing the times. Their formulation uses translations
and dilations. We also allow orthogonal transformations. For a bounded
forward-global solution, the quantifiers have the order
\begin{equation}\label{main:quantifiers}
 \exists\ J,\ (Q_j,\ell_j)_{j=1}^J,\ v_L,\qquad
 \forall\ t_n\to\infty,\qquad
 \exists\ (a_{j,n},\lambda_{j,n},O_{j,n})_{j,n}.
\end{equation}
Here $J\ge0$ is finite, $v_L$ is a free wave, and each traveling wave
has a nonzero rest profile $Q_j$ and velocity $|\ell_j|<1$. We take
$O_{j,n}\in O(3)$.
With $S_{j,n}=S_{a_{j,n},\lambda_{j,n},O_{j,n}}$, the decomposition
and separation conditions are
\begin{align}
 \vec u(t_n)
 &=\vec v_L(t_n)+\sum_{j=1}^J
       S_{j,n}\vec Q_{j,\ell_j}+o_{\mathcal H}(1),
 \label{main:FP}\\
 \frac{\lambda_{i,n}}{\lambda_{j,n}}
 +\frac{\lambda_{j,n}}{\lambda_{i,n}}
 +\frac{|a_{i,n}-a_{j,n}|^2}{\lambda_{i,n}\lambda_{j,n}}
 &\longrightarrow\infty,\qquad i\ne j.
 \label{main:profile-separation}
\end{align}

\begin{theorem}[Failure of fixed-profile soliton resolution]\label{main:counterexample}
There exist a sign-changing stationary state $Q$, a smooth chart
$\{Q_p: p\in\R^3,\ |p|<\delta\}\subset GQ$ with $\delta>0$ and $Q_0=Q$,
and a unit vector $e\in\R^2$ with the following property.
For each sufficiently small $d>0$ there are a real forward-global solution
$u$ of \eqref{main:nlw}, even separately in $x_1,x_2$, a $C^1$ path
$p$ and a free wave $v_+$ such that
\begin{align}
 \|\vec u(0)-(Q,0)\|_{\mathcal H}
       +\|\vec v_+(0)\|_{\mathcal H}
 &\lesssim d^3,\qquad
 \sup_{t\ge0}\|\vec u(t)\|_{\mathcal H}<\infty,
 \label{main:counterexample-bounds}\\
 \|p\|_\infty&\lesssim d,\qquad
 \vec u(t)=(Q_{p(t)},0)+\vec v_+(t)+o_{\mathcal H}(1).
 \label{main:scattering}
\end{align}
For $T_j=4^jd^{-5}$, the stationary component satisfies
\begin{equation}\label{main:two-limits}
 \begin{gathered}
 Q_{p(T_{2k})}\longrightarrow Q,\qquad
 Q_{p(T_{2k+1})}\longrightarrow Q_{(de,0)}
       \quad\text{in }\dot H^1,\\
 \operatorname{dist}_{\dot H^1}
       (Q_{(de,0)},\mathrm{Sim}\cdot Q)\gtrsim d.
 \end{gathered}
\end{equation}
No fixed list in \eqref{main:quantifiers} satisfies
\eqref{main:FP}--\eqref{main:profile-separation} even on $t_n=T_n$.
In particular, $u$ admits no full-time resolution by a fixed finite list.
\end{theorem}

Suppose that a fixed list satisfies \eqref{main:FP} and
\eqref{main:profile-separation} on $T_n$. Applying the norm identity
\eqref{fp:exhaustion} to the even and odd subsequences in turn, we obtain
\[
 \begin{aligned}
 t=T_{2k}   &\quad\Longrightarrow\quad v_L=v_+,\quad J=1,\\
 t=T_{2k+1} &\quad\Longrightarrow\quad Q_{(de,0)}\in\mathrm{Sim}\cdot Q .
 \end{aligned}
\]
The second conclusion contradicts \eqref{main:two-limits}. On the other
hand, $Q_{p(t)}\in GQ$ for all $t$, which suggests that we allow each
core to move within a fixed conformal class. If we choose the profile
separately on the even and odd subsequences, the two endpoint profiles
suffice. At all late times, even a finite choice of lists fails, as we
show in Corollary~\ref{fp:catalogue} using the connected set of limits.

As in five dimensions \cite{CHS5D}, a packet changes the core velocity
at quadratic order. We space the packets so that this small change
produces a fixed displacement. In the five-dimensional construction,
each tangent to the stationary family has an $O(|x|^{-3})$ tail and
belongs to $L^2$. In three
dimensions the double-even sector contains two $L^2$ modes $T,K$ and
a resonance with a $|x|^{-1}$ tail. Motion in the resonance direction
therefore contributes a term outside $L^2$ to the time derivative of the
core. The dipole of $T$ also contributes a rank-one term to the retarded
response. We compute the residues in \eqref{fr:eq:residue}, express
them in chart coordinates in \eqref{non:eq:residue-coordinates}, and
use \eqref{non:eq:threshold-cancellation} to recover an $L^2$ velocity.

\subsection{Resolution in a fixed conformal class}\label{main:classes}
For a stationary state
\begin{equation}\label{main:stationary}
 P\in\dot H^1(\R^3)\setminus\{0\},\qquad -\Delta P=P^5 ,
\end{equation}
we denote by $\mathcal K_P=O(4)P$ the compact orbit under isometries
of the round sphere. Every element of $GP$ is a translate and dilate
of some $q\in\mathcal K_P$. We use the conventions
\begin{align}
 \bigl[S_{a,\lambda,O}(f,g)\bigr](a+\lambda Oy)
   &=\bigl(\lambda^{-\frac{1}{2}}f(y),\lambda^{-\frac{3}{2}}g(y)\bigr),
       \label{main:similarity}\\
 P_\ell(y)&=P\bigl((I-\ell\otimes\ell)^{-\frac{1}{2}}y\bigr),\qquad
 \vec P_\ell=(P_\ell,-\ell\cdot\nabla P_\ell).
       \label{main:lorentz}
\end{align}
Then $P_\ell(x-t\ell)$ solves \eqref{main:nlw}. With
$S_{a,\lambda}=S_{a,\lambda,I}$, we set
\begin{equation}\label{main:class-orbit}
 \mathcal M_\ell(P)
    =\{S_{a,\lambda}\vec q_\ell:q\in\mathcal K_P,\ a\in\R^3,\ \lambda>0\}.
\end{equation}

\noindent\textit{Conformal-class conjecture.}
For each bounded forward-global solution, we ask whether there are
finitely many fixed classes $GQ_j$, fixed velocities $|\ell_j|<1$ and
one free wave $v_L$ such that
\[
 \vec u(t)=\vec v_L(t)+\sum_{j=1}^J\mathbf p_j(t)+o_{\mathcal H}(1),
 \qquad \mathbf p_j(t)\in\mathcal M_{\ell_j}(Q_j)\qquad(t\to\infty),
\]
where the separation condition \eqref{main:profile-separation} holds
for the centers and scales as $t\to\infty$.
For a single core, the full conformal nondegeneracy condition of
Duyckaerts, Kenig and Merle \cite{DKMcompact} allows us to pass from
one sequence to all late times.

\begin{theorem}[Full-time resolution in one conformal class]\label{cap:main}
Let $u$ be a forward-global strong solution of \eqref{main:nlw} bounded
in $\mathcal H$. Fix $|\ell|<1$ and a nonzero stationary state $P$ satisfying
\[
 \ker_{\dot H^1(\R^3)}(-\Delta-5P^4)=T_P(GP).
\]
Suppose there are a free energy wave $v$ and a sequence $t_n\to\infty$
such that
\[
 \operatorname{dist}_{\mathcal H}
       \bigl(\vec u(t_n)-\vec v(t_n),\mathcal M_\ell(P)\bigr)\longrightarrow0.
\]
Then $v$ is the canonical free radiation of $u$. On a late half-line
there are $C^1$ parameters $a(t)\in\R^3$, $\lambda(t)>0$ and
$q(t)\in\mathcal K_P$ such that, with this same free wave,
\[
 \vec u(t)=\vec v(t)+S_{a(t),\lambda(t)}\vec q_\ell(t)
                   +o_{\mathcal H}(1)\qquad(t\to\infty).
\]
Here $\vec q_\ell(t)$ is the phase in \eqref{main:lorentz} evaluated at
$q(t)$. The parameters satisfy
\[
 \lambda'(t)\longrightarrow0,\qquad a'(t)\longrightarrow\ell,\qquad
 \lambda(t)\|q'(t)\|_{\dot H^1}\longrightarrow0 .
\]
\end{theorem}

For Theorem~\ref{main:counterexample}, we may take $P=Q$, $\ell=0$,
$v=v_+$ and $t_n=T_{2n}$. By stereographic covariance and
Theorem~\ref{cert:thm:gap},
\[
 \ker_{\dot H^1}(-\Delta-5Q^4)=T_Q(GQ),\qquad \dim T_Q(GQ)=8 .
\]
Thus the counterexample satisfies the full kernel hypothesis of
Theorem~\ref{cap:main}. Its odd subsequence rules out a fixed representative
modulo similarities even under this hypothesis.

\smallskip
\noindent\textit{The three-dimensional radial case.}
Duyckaerts, Kenig and Merle \cite[Theorem~1]{DKMradial} proved that a
forward-global radial finite-energy solution resolves into a free wave
and signed ground states at separated scales. Shen \cite[Theorem~1.2]{Shen}
proved that at most one ground state occurs. Thus the solution either
scatters or admits a radial free wave $v$, a fixed sign
$\iota\in\{-1,1\}$ and a scale $\lambda(t)>0$ such that
\[
 \vec u(t)=\vec v(t)+(\iota W_{\lambda(t)},0)+o_{\mathcal H}(1),
 \qquad \frac{\lambda(t)}{t}\longrightarrow0,
\]
where
\[
 W(x)=\left(1+\frac{|x|^2}{3}\right)^{-1/2},\qquad
 W_\lambda(x)=\lambda^{-1/2}W(x/\lambda).
\]
After a fixed dilation, $W$ corresponds to a constant on $\Sph^3$.
The four first spherical harmonics span the kernel and generate
translations and dilation. Hence $P=\iota W$ satisfies the hypothesis
of Theorem~\ref{cap:main}. On the radial nonscattering branch we take
$\ell=0$, $a(t)=0$ and $q(t)=\iota W$, so only $\lambda(t)$ varies.

\subsection{Relation to earlier work}\label{main:literature}
In five dimensions, Jendrej and Martel \cite[Theorem~1]{JM5D}
constructed, for every integer $N\ge2$, a global nonradial solution
with $N$ ground-state bubbles concentrating at $N$ distinct prescribed
points. Each bubble is radial about its own center. Shen's conjectured
dimension-dependent bound \cite[Remark~1.5]{Shen} concerns solutions
radial about a common center. In our five-dimensional counterexample
\cite[Theorem~1.1 and Remark~1.2]{CHS5D}, one nonradial core has a continuum of limiting
similarity classes within a single conformal orbit. This motion of one
core disproves the fixed-profile assertion of Duyckaerts, Kenig and Merle
\cite[Conjecture~5.5]{DKM} in dimension five.

\smallskip
\noindent\textit{Remark (Conformal zero modes and radial bubbles).}
For the state used in five dimensions, all conformal zero modes are
$O(r^{-3})$ at infinity and belong to $L^2$ \cite[Section~3]{CHS5D}. In three dimensions varying
the coefficient of $r^{-1}$ gives a resonance in $\dot H^1\setminus L^2$.
At fixed coefficient, the $L^2$ mode $K$ is still transverse to the
similarity orbit. Thus the endpoints in \eqref{main:two-limits} have
the same stationary energy and the same leading tail.

For two concentric positive ground states, let $\delta$ denote the ratio
of the inner scale to the outer scale. Let $I_N(\delta)$ denote the leading
term of the interaction energy in dimension $N$. In three dimensions,
we use the stationary equation, $W^4\in L^1$ and
$W_\delta(x)\lesssim\delta^{1/2}/|x|$ to obtain
\[
 \begin{aligned}
 \int\nabla W\cdot\nabla W_\delta\dd x
 &=\int W^5W_\delta\dd x=\int W W_\delta^5\dd x,\\
 \int\bigl(W^2W_\delta^4+W^4W_\delta^2\bigr)\dd x&=O(\delta).
 \end{aligned}
\]
Cauchy--Schwarz gives the same bound for $\int W^3W_\delta^3$.
Expanding \eqref{main:energy} and changing variables, we obtain
\[
 E(W+W_\delta,0)-2E(W,0)
 =-\delta^{1/2}\int W(y)^5W(\delta y)\dd y+O(\delta).
\]
By dominated convergence, $c_3=W(0)\int W^5\dd x>0$.
In five dimensions, the two-bubble energy expansion of Jendrej and
Lawrie \cite{JL} gives $c_5>0$. The leading terms satisfy
\[
 \begin{aligned}
 I_3(\delta)&=-c_3\delta^{1/2},&
 |I_3'(\delta)|&\longrightarrow\infty,\\
 I_5(\delta)&=-c_5\delta^{3/2},&
 |I_5'(\delta)|&\longrightarrow0
 \end{aligned}
 \qquad(\delta\downarrow0).
\]
In three dimensions, Shen rules out two persistent radial scales by
showing that they would force radiation to concentrate at every
sufficiently late time \cite[Theorem~1.2 and Sections~1.3 and~6]{Shen}.
In five dimensions, an inner oscillation of relative amplitude
$\delta^{3/4}$ has a quadratic mean of order $\delta^{3/2}$, the order
of the outer approximation at the origin in inner variables. Matching
the mean correction there to the full outer value at each new scale
suggests how a radial solution might sustain any prescribed finite
number of concentric bubbles.
\smallskip

Duyckaerts, Jia, Kenig and Merle \cite{DJKM} proved soliton resolution
along a sequence of times for solutions bounded in the energy space
in dimensions three through five, without radial symmetry. If their decomposition has one profile whose full
energy kernel consists of conformal tangents, Theorem~\ref{cap:main}
extends it to all late times with the same radiation.

Duyckaerts, Kenig and Merle \cite{DKMcompact} developed modulation on
the stationary conformal orbit and proved rigidity for exact solutions
with the compactness property under full conformal nondegeneracy.
Their exterior radiation theorem \cite{DKMradiation} identifies the
canonical free wave $v$ of a bounded forward-global solution. Once we
subtract $v$, the equation for the remainder becomes
\[
 U=u-v,\qquad \square U=(U+v)^5 .
\]
To control $U$ between observation times, we use the exterior energy
of the growing modes in Lemma~\ref{cap:linear}. For a single fully
nondegenerate profile, \eqref{cap:late} then keeps the transverse error
small after one sufficiently late entry, including when the profile
changes sign and $v\ne0$.
Collot, Duyckaerts, Kenig and Merle \cite{CDKM} study compactly supported
nonradial data with ground-state traveling profiles of distinct velocities.
They prove a resolution/collision dichotomy and obtain full-time resolution
when there are at most two such profiles.

Ding \cite{Ding} proved the existence of sign-changing stationary
solutions. Del Pino, Musso, Pacard and Pistoia \cite{DMPP} constructed
large-energy nodal families, and Musso and Wei \cite{MW} studied their
conformal nondegeneracy. Our stationary state belongs to the single-node
isoparametric family of Fern\'andez and Petean \cite{FP}. Neves and
Pistoia \cite[Theorems~A and~B]{NP} proved nondegeneracy in the
$O(k)\times O(k)$-invariant sector for $k\ge2$ and full nondegeneracy
for sufficiently large $k$, where the dimension is $2k-1$. In
Appendix~\ref{cert:start} we treat the full kernel for $k=2$. For related reductions by harmonic morphisms,
see Clapp and Pistoia \cite{CP} and Fern\'andez, Palmas and Petean \cite{FPP}.

For the three-dimensional quintic equation, Krieger and Schlag \cite{KS}
and Beceanu \cite{Beceanu} treated stable dynamics near the ground state
in the presence of a zero-energy resonance. Krieger, Schlag and Tataru
\cite{KST} constructed type-II solutions with controlled concentration.
Donninger and Krieger \cite{DK} constructed global solutions with continuing
scale evolution, and Pillai \cite{Pillai} prescribed admissible scales,
including oscillatory ones. In these examples the core remains a rescaled
ground state. Jendrej \cite{Jendrej} prescribed the regular part at
blow-up for the five-dimensional wave equation. Jendrej, Lawrie and
Rodriguez \cite{JLR} prescribed the radiation for equivariant wave maps.

For the harmonic map heat flow between two-spheres, Topping
\cite{Topping97} proved uniqueness of the infinite-time bubble locations
when the limiting map and all bubbles have the same orientation.
He also constructed finite-time singularities with nonunique limiting
bubbles for other targets \cite{Topping04}. Engelstein and Mendelson
\cite{EM} constructed wave maps with nonunique limiting bubbles.
Gustafson, Nakanishi and Tsai \cite{GNT} constructed degree-two
harmonic map heat flows with oscillating scale.
Jendrej, Lawrie and Schlag \cite{JLS} distinguish convergence to a family
of bubble configurations from selection of fixed profiles.

Kenig and Merle \cite{KM} proved scattering below the ground-state energy
when the initial $\dot H^1$ norm is below that of the ground state. For
radial soliton resolution in higher dimensions, see Duyckaerts, Kenig,
Martel and Merle \cite{DKMM}, Duyckaerts, Kenig and Merle \cite{DKModd},
and Jendrej and Lawrie \cite{JL}.

\subsection{Construction and proof of the counterexample}\label{main:outline}
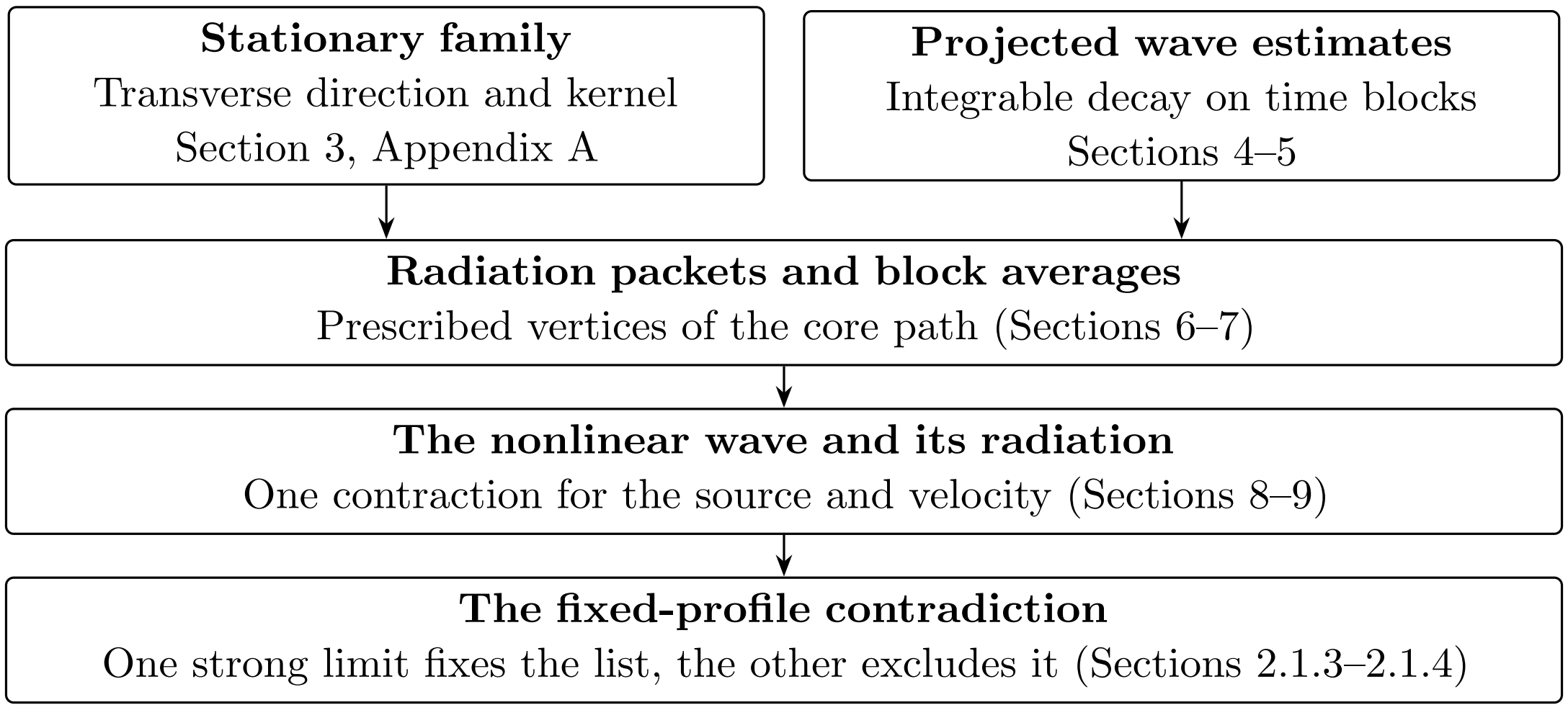
\begin{figure}[!htbp]
\centering
\begin{tikzpicture}[>=Stealth,font=\small,semithick,
 box/.style={draw,rounded corners=2pt,align=center,inner sep=4pt},
 half/.style={box,text width=6.65cm,minimum height=1.52cm},
 full/.style={box,text width=14cm,minimum height=1.10cm}]
 \node[half] (geo) at (-3.65,0) {\textbf{Stationary family}\\
 Transverse direction and kernel\\Section~\ref{geo:start}, Appendix~\ref{cert:start}};
 \node[half] (lin) at (3.65,0) {\textbf{Projected wave estimates}\\
 Integrable decay on time blocks\\Sections~\ref{lin:start}--\ref{lin:old-start}};
 \node[full] (pack) at (0,-1.9) {\textbf{Radiation packets and block averages}\\
 Prescribed vertices of the core path (Sections~\ref{cat:start}--\ref{non:start})};
 \node[full] (wave) at (0,-3.45) {\textbf{The nonlinear wave and its radiation}\\
 One contraction for the source and velocity (Sections~\ref{end:start}--\ref{end:global-section})};
 \node[full] (contra) at (0,-5.0) {\textbf{The fixed-profile contradiction}\\
 One strong limit fixes the list, the other excludes it (Sections~\ref{fp:start}--\ref{fp:energy-shape})};
 \draw[->] (geo.south)--(pack.north -| geo.south);
 \draw[->] (lin.south)--(pack.north -| lin.south);
 \draw[->] (pack.south)--(wave.north);
 \draw[->] (wave.south)--(contra.north);
\end{tikzpicture}
\caption{Construction of the solution and proof of Theorem~\ref{main:counterexample}.}
\label{main:proof-diagram}
\end{figure}

In five dimensions, we used incoming packets with summable amplitudes
to move a stationary core by a fixed distance on each successive
interval \cite{CHS5D}. In three dimensions we must also control the
$|x|^{-1}$ resonance and the dipole of an $L^2$ zero mode. We give the
four steps of the proof in Figure~\ref{main:proof-diagram}.

\smallskip
\noindent\textbf{Step 1. The stationary chart and the projected evolution.}
Let $w$ be the balanced single-node profile of Appendix~\ref{cert:start}.
On $\Sph^3\subset\mathbb C^2$, its Hopf lift $q$ satisfies
\begin{equation}\label{main:hopf-lift}
 \begin{gathered}
 \pi_H(z_1,z_2)=(2z_1\overline{z_2},\,|z_1|^2-|z_2|^2),\qquad
 q(\cos\chi\,e^{i\alpha},\sin\chi\,e^{i\beta})
       =\left(\frac{3}{4}\right)^{\frac{1}{4}}w(2\chi),\\
 q^{-1}(0)=\left\{(z_1,z_2)\in\Sph^3:
                   |z_1|=|z_2|=\frac{1}{\sqrt2}\right\}.
 \end{gathered}
\end{equation}
We rotate this nodal torus and apply the normalization
\eqref{geo:eq:normalize}. At fixed monopole and zero dipole,
Proposition~\ref{geo:prop:shape} identifies a direction transverse to
$\mathrm{Sim}\cdot Q$. For $L_*=-\Delta-5Q^4$, the chart even in each
of $x_1,x_2$ satisfies
\begin{equation}\label{main:kernel-types}
 \ker_{\dot H^1,\mathrm{even}}L_*=\operatorname{span}\{E_1,E_2,R_*\},
 \qquad \ker_{L^2,\mathrm{even}}L_*=\operatorname{span}\{E_1,E_2\}.
\end{equation}
The simple-pole residue contains both the monopole term in $R_*$
and the dipole term in the $L^2$ kernel. The projection $\Pi$ removes
the double pole and the entire simple-pole residue. In
Section~\ref{lin:start} we prove Theorem~\ref{fr:thm:full} and the
projected decay estimate
\[
 \|S_*(t)\Pi F\|_X\lesssim\langle t\rangle^{-\frac{3}{2}}\|F\|_Y .
\]

\smallskip
\noindent\textbf{Step 2. Packet impulses and the reference vertices.}
For $V=5Q^4$, we differentiate the conformally invariant identity
\eqref{ray:eq:riesz} along the chart. Since the line transform is
injective and $XV>0$, the balance \eqref{ray:eq:balance} implies that
the ray impulses generate $\R^3$ as a positive cone. We realize them
by packets of amplitude $\eta_j$ on intervals of length comparable to
$T_j$, where
\begin{equation}\label{main:scales-intro}
 \eta_j=d^32^{-j},\qquad T_j=4^jd^{-5},\qquad
 \sum_{j\ge0}\eta_j=2d^3,\qquad T_j\eta_j^2=d .
\end{equation}
Thus a velocity change of order $\eta_j^2$ produces a displacement of
order $d$ on the $j$th interval. Proposition~\ref{construction:controls}
bounds the integrated impulse error. We choose the packet coefficients by
prescribing the average velocity on each interval, so that the finite
coordinates satisfy
\[
 y(e_j)=\bar y(e_j),\qquad e_j=T_{j+1}.
\]

\smallskip
\noindent\textbf{Step 3. Threshold cancellation and the exact wave.}
We introduce the retarded correction to the moving core. By
\eqref{non:eq:threshold-cancellation}, its frozen threshold terms match the
two residues from Step~1. We then solve for the projected source $f$
and the velocity correction $b$ in one contraction. At the fixed point
the residual satisfies
\[
 \Pi(\square u-u^5)=0,\qquad \langle Z_p,\square u-u^5\rangle=0.
\]
Since the tangent-duality matrix is invertible, $\square u=u^5$.
The integrated momentum identity \eqref{end:eq:momentum} controls the
velocity on each block. The scattering estimate \eqref{end:eq:physicalscattering}
identifies a single wave $v_+$. Alternating the reference vertices
between $0$ and $de$, we obtain the two strong limits in
\eqref{main:two-limits}.

\smallskip
\noindent\textbf{Step 4. The two subsequences and the fixed list.}
We use free-wave and profile orthogonality
\eqref{fp:dislocation}--\eqref{fp:free-separation}, in the setting of
the profile decomposition of Bahouri and G\'erard~\cite{BG}. On the
even subsequence, at most one similarity in a separated list can stay
bounded, and the nonzero limit requires one such index $j_*$. Writing
$\mathbf f_j=\vec Q_{j,\ell_j}$, we obtain from \eqref{fp:exhaustion}
\[
 0=\|\vec v_L(0)-\vec v_+(0)\|_{\mathcal H}^2+
       \sum_{j\ne j_*}\|\mathbf f_j\|_{\mathcal H}^2 .
\]
Thus $v_L=v_+$ and $J=1$. On the odd subsequence the same profile would
represent $Q_{(de,0)}$, in contradiction with \eqref{main:two-limits}.
In Section~\ref{fp:energy-shape} we also use continuity and precompactness
of the path to show that the distance of its limiting profiles from
$\mathrm{Sim}\cdot Q$ takes every value between the endpoint distances.

\FloatBarrier

For this solution, the even subsequence determines the class and the
free wave, and the odd subsequence requires a change of shape within
that class. For a general one-core solution, full conformal
nondegeneracy identifies the neutral directions with the tangent space
to the orbit. We use the exterior energy to detect a growing transverse
mode, then use the modulation equations to control the other components.

\subsection{Proof of the single-core theorem}\label{main:cap-outline}
\begin{figure}[!htbp]
\centering
\begin{tikzpicture}[>=Stealth,font=\small,semithick,
 box/.style={draw,rounded corners=2pt,align=center,inner sep=4pt},
 half/.style={box,text width=6.65cm,minimum height=1.52cm},
 full/.style={box,text width=14cm,minimum height=1.10cm}]
 \node[half] (rad) at (-3.65,0) {\textbf{One core along a sequence}\\
 Canonical radiation and limiting action\\Section~\ref{cap:sec-radiation}};
 \node[half] (ext) at (3.65,0) {\textbf{Full conformal nondegeneracy}\\
 Linearized control on the exterior cone\\Section~\ref{cap:sec-exterior}};
 \node[full] (mod) at (0,-1.9) {\textbf{Orthogonal modulation}\\
 Control of the error and damping of the mode coefficients (Section~\ref{cap:sec-modulation})};
 \node[full] (con) at (0,-3.45) {\textbf{Convergence on a half-line}\\
 Late entry, continuation and decay with the same radiation (Section~\ref{cap:sec-convergence})};
 \draw[->] (rad.south)--(mod.north -| rad.south);
 \draw[->] (ext.south)--(mod.north -| ext.south);
 \draw[->] (mod.south)--(con.north);
\end{tikzpicture}
\caption{Estimates leading from one sequence to all late times in Theorem~\ref{cap:main}.}
\label{cap:proof-diagram}
\end{figure}
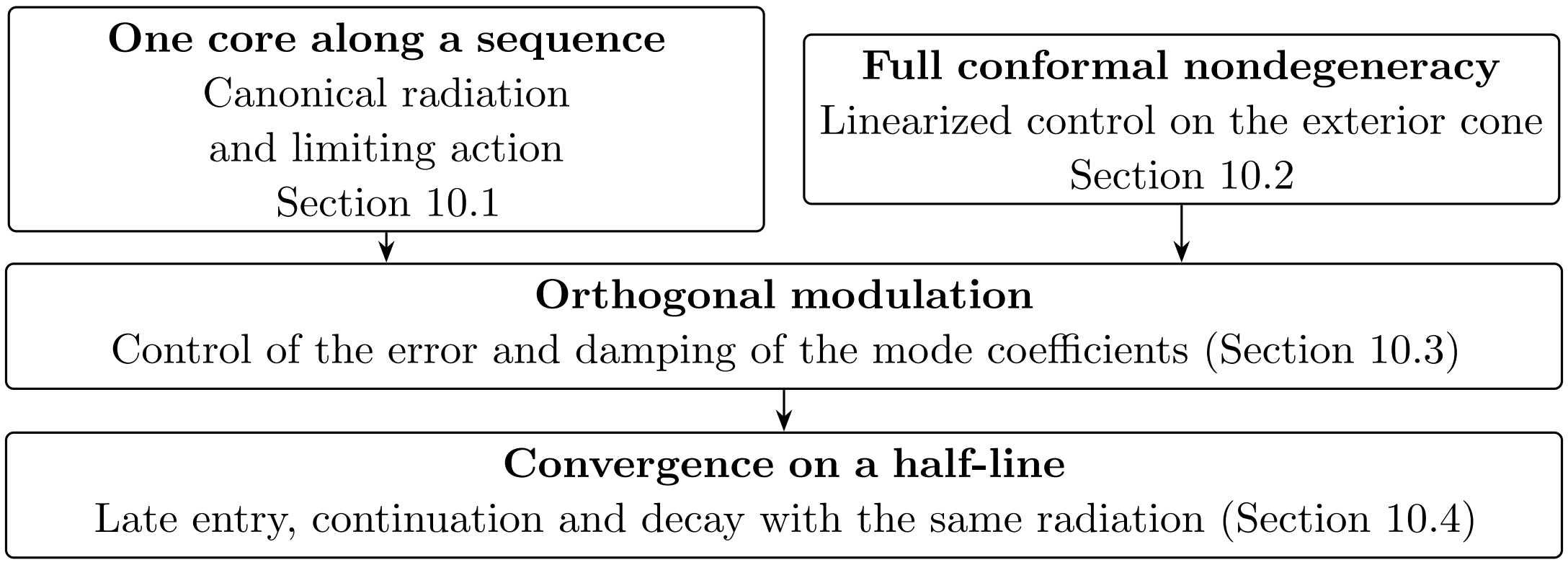

We fix the free wave $v$ and the velocity $\ell$ in the hypothesis of
Theorem~\ref{cap:main}. Throughout the proof we write
\[
 U=u-v,\qquad d(t)=\operatorname{dist}_{\mathcal H}
       \bigl(\vec U(t),\mathcal M_\ell(P)\bigr).
\]
The hypothesis gives $d(t_n)\to0$. To estimate $d(t)$ between these
times, we use the exterior energy of $U$ and the linearized equation
at a traveling wave, as indicated in Figure~\ref{cap:proof-diagram}.

\smallskip
\noindent\textbf{Step 1. The free radiation and the limiting action.}
By \eqref{fp:free-separation}, free waves are asymptotically orthogonal
to all translates and dilates of $\{\vec q_\ell:q\in\mathcal K_P\}$.
Testing the decomposition at $t_n$ against free waves identifies $v=v_{\rm can}$.
We apply the exterior radiation theorem \eqref{cap:radiation} to $U$.
By conservation of energy and momentum and decoupling from $v$, the action
$\mathcal J_\ell$ of Section~\ref{cap:sec-radiation} satisfies
\[
 \mathcal J_\ell(\vec U(t))\longrightarrow J_*=
                    \sqrt{1-|\ell|^2}\,E(P,0).
\]
We collect the error in the action and the free-wave tails in the
decreasing quantity $\eta(t)\to0$ of \eqref{cap:eta}.

\smallskip
\noindent\textbf{Step 2. The linearized estimate on the exterior cone.}
By full conformal nondegeneracy, the kernel of the Hessian $\mathfrak b_q$
is tangent to the traveling-wave orbit. Let $z$ solve the linearized
equation with forcing $F$ and data $H=(h,k)$, with $h$ orthogonal to the
spatial tangents. For $\chi=\mathbf1_{\{|y|>s\}}$, Lemma~\ref{cap:linear}
asserts that
\[
 \|H\|_{\mathcal H}^2\lesssim|\beta_q(H)|^2+[\mathfrak b_q(H,H)]_+
 +\|\chi F\|_{L_s^1L_y^2}^2+
            \limsup_{s\to\infty}\|\chi\nabla_{s,y}z(s)\|_2^2 .
\]
Here $\beta_q$ tests against the boosted modes of the negative spectrum.
The exterior energy detects the growing modes in $\ker\beta_q$.

\FloatBarrier

\smallskip
\noindent\textbf{Step 3. Modulation and control of the error.}
While $d(t)$ is small, we choose $(a,\lambda,q)$ by spatial
orthogonality and set
\[
 \varepsilon=\|H\|_{\mathcal H}\simeq d(t),\qquad
 b=|\beta_q(H)|,\qquad \dd\tau=\frac{\dd t}{\lambda(t)} .
\]
After we rescale at time $t$, Step~1 gives vanishing exterior energy
at infinity in Step~2. We expand the action at $\vec q_\ell$ and
differentiate the mode coefficients. By \eqref{cap:control} and
\eqref{cap:modulation}, we obtain
\[
 \varepsilon(t)\lesssim b(t)+\eta(t),\qquad
 D_\tau^+b+c_0b\lesssim\varepsilon^2+\eta,\qquad c_0>0.
\]

\smallskip
\noindent\textbf{Step 4. Convergence after one late time.}

We choose $T=t_n$ sufficiently large to absorb the quadratic error in
Step~3 while the core lies in the modulation neighborhood. By the integrated estimate
\eqref{cap:late}, we keep the core in this neighborhood for all $t\ge T$.
The bound on $\lambda'$ gives $\tau(t)\to\infty$, so the transverse
error and the modulation rates tend to zero for the same free wave $v$.
\smallskip
\noindent\textit{Organization and notation.}
In Section~\ref{construction:negative} we derive the obstruction from
\eqref{fp:exhaustion}. Sections~\ref{geo:start}--\ref{lin:old-start}
give the stationary and linear estimates. We construct the packets and
the modulated core in Sections~\ref{cat:start}--\ref{non:start}, then
solve the coupled equations and identify the free wave in
Sections~\ref{end:start}--\ref{end:global-section}. The proof in
Section~\ref{proof:capture} uses \eqref{fp:free-separation} and the
exterior and modulation estimates of Section~\ref{cap:start}.
Appendix~\ref{cert:start} gives the kernel calculation for both theorems,
and Appendix~\ref{notation:start} collects the notation.

All spatial norms are over $\R^3$ unless we specify another domain.
We write $A\lesssim B$ when $A\le CB$. We fix the stationary chart and
the packet family before choosing $d$, so $C$ is independent of $d$
and of the block index. A subscript records any additional dependence.
We sum over repeated chart indices and write
$\langle x\rangle=(1+|x|^2)^{1/2}$. We denote the free wave group on
$\mathcal H$ by $U_0(t)$ and use
\[
 \square=\partial_t^2-\Delta,\qquad D=\sqrt{-\Delta},\qquad
 S_{\rm free}(t)=\mathbf 1_{\{t>0\}}\frac{\sin(tD)}{D}.
\]

\section{Proof of the main theorems}\label{construction:start}
For Theorem~\ref{main:counterexample}, we prescribe alternating endpoint
shapes and test a fixed list against both limits. For
Theorem~\ref{cap:main}, we use the exterior radiation and the modulation
equations to continue a one-core decomposition from a sequence to a
half-line.

\subsection{Construction of the counterexample}\label{construction:negative}
We deduce Theorem~\ref{main:counterexample} from the packet and wave constructions
in Sections~\ref{cat:start}--\ref{end:global-section}.
We use $Q,Q_p,e$ from Proposition~\ref{geo:ass:stationary}.

\subsubsection{Packet controls and their integrated error}\label{construction:controls-section}
We approximate the changes in core velocity by packet impulses and
estimate the error after one time integration.

We choose the scales
\begin{equation}\label{main:scales}
 \eta_0=d^3,\qquad \eta_j=\eta_02^{-j},\qquad T_j=4^jd^{-5},
 \qquad \sum_{j\ge0}\eta_j=2\eta_0,\qquad T_j\eta_j^2=d.
\end{equation}
We take the intervals and their right endpoints to be
\[
 I_0=[0,4T_0],\qquad I_j=[T_j,4T_j]\quad(j\ge1),\qquad
 e_j=4T_j=T_{j+1}\quad(j\ge0).
\]
We write $p=(y,a)\in\R^2\times\R$, where $a$ is the monopole increment,
and denote the first two coordinates of $z\in\R^3$ by $z_y$. At the
base point, we set
\begin{equation}\label{construction:eq:frozen-frame}
 Z_*=(E_{*,1},E_{*,2},R_*),\qquad
 \mathsf M=\operatorname{diag}
       \bigl((\langle E_{*,i},E_{*,j}\rangle)_{ij},4\pi\bigr),
 \qquad m_0(v)=\langle Z_*,10Q^3v^2\rangle.
\end{equation}
\begin{proposition}[Packet impulses and integrated tails]\label{construction:controls}
For controls $c=(c_j)$ in a fixed bounded ball of
$\ell^\infty(\R^2)$, there is a smooth packet family such that
\begin{equation}\label{cat:eq:train}
 v_c(t)=\sum_{j\ge0}\eta_j\mathcal V_{c_j}(t-e_j),\qquad
 P(c)=\mathbf1_{[0,\infty)}\mathsf M^{-1}m_0(v_c).
\end{equation}
The real double-even frozen wave $v_c$ scatters freely, and
\begin{equation}\label{cat:eq:train-bounds}
 \begin{gathered}
 \|\vec v_c\|_{L^\infty\mathcal H}
 +\|v_c\|_{\mathcal T(X)}+\|\rho v_c\|_\infty
 +\|\partial_tv_c\|_{\mathcal T(\mathcal O)}\lesssim\eta_0,\\
 \|P(c)\|_{L^1(I_j)}\lesssim\eta_j^2,\qquad
 \|P(c)\|_{\mathcal T}\lesssim\eta_0^2.
 \end{gathered}
\end{equation}
The integrated error satisfies
\begin{equation}\label{cat:eq:integrated-error}
 e_c(t)=\int_t^\infty[P(c)(s)]_y\dd s
        -\sum_{e_i>t}\eta_i^2c_i,\qquad
 \|e_c\|_{L^1(0,\infty)}\lesssim\eta_0^2.
\end{equation}
All bounds hold for differences, with the additional factor
$\|c-\widetilde c\|_{\ell^\infty}$.
\end{proposition}

\subsubsection{Realization of the reference polygon}\label{construction:response-section}
We prescribe a continuous reference polygon $\bar y$ with
\begin{equation}\label{construction:eq:polygon}
 \bar y(0)=0,\qquad \bar y|_{I_j}\text{ affine},\qquad
 |\bar y(e_j)|\le d.
\end{equation}
In particular, $|\bar y'|\lesssim\eta_j^2$ on $I_j$.

\begin{proposition}[Realization of the reference polygon]\label{construction:response}
For sufficiently small $d$, every reference polygon satisfying
\eqref{construction:eq:polygon} has a real double-even forward-global
solution $u$, a $C^1$ path $p$, and a free wave $v_+$ satisfying
\eqref{main:counterexample-bounds}--\eqref{main:scattering}.
Moreover,
\begin{equation}\label{construction:eq:realization}
 \begin{gathered}
 \|p-(\bar y,0)\|_{L^\infty(0,\infty)}\lesssim d^2,\qquad
 \|p'\|_{L^2(0,\infty)}\lesssim d^3,\\
 p(e_j)-(\bar y(e_j),0)\longrightarrow0.
 \end{gathered}
\end{equation}
\end{proposition}
In Sections~\ref{non:start}--\ref{end:global-section}, we determine
the controls by block averages and solve for the source and velocity
by a single contraction.

\subsubsection{Identification of the radiation and profile from one strong limit}\label{fp:start}\label{fp:radiation-energy}
We first give the two orthogonality relations needed to compare a
separated list of profiles with a nonzero strong limit.

For each fixed $\mathbf f\in\mathcal H$, the similarity action satisfies
\begin{equation}\label{fp:dislocation}
 \begin{aligned}
 S_n\to S&\quad\Longrightarrow\quad S_n\mathbf f\to S\mathbf f,\\
 S_n\text{ escapes compact parameter sets}
   &\quad\Longrightarrow\quad S_n\mathbf f\rightharpoonup0.
 \end{aligned}
\end{equation}
For scales tending to zero or infinity, we approximate in Fourier
space by annuli. For bounded scales and escaping centers, we use compactly
supported data. In both cases the rotations range over a compact group.

For a free energy wave $v$ and fixed $\mathbf f\in\mathcal H$, we also have
\begin{equation}\label{fp:free-separation}
 \sup_{S\in\mathrm{Sim}}
 |\langle\vec v(t),S\mathbf f\rangle_{\mathcal H}|
 \longrightarrow0\qquad(t\to\infty).
\end{equation}
By the free Strichartz estimate \cite[Corollary~1.3]{KT}, we have
\begin{equation}\label{fp:eq:free-strichartz}
 \|u\|_{L^5(I;L^{10})}
 +\sup_{t\in I}\|\vec u(t)\|_{\mathcal H}
 \lesssim\|\vec u(t_0)\|_{\mathcal H}
                +\|\square u\|_{L^1(I;L^2)},\qquad t_0\in I.
\end{equation}
The constant is independent of the interval $I$.
For arbitrary $t_n\to\infty$ and $S_n=S_{a_n,\lambda_n,O_n}$,
let $w_n$ be the free wave with
$\vec w_n(0)=S_n^{-1}\vec v(t_n)$. By scaling, we have
\begin{equation}\label{fp:eq:strichartz-tail}
 \|w_n\|_{L^5([0,\infty);L^{10})}
 =\|v\|_{L^5([t_n,\infty);L^{10})}\longrightarrow0.
\end{equation}
Energy conservation and similarity invariance give
$\|\vec w_n(0)\|_{\mathcal H}=\|\vec v(0)\|_{\mathcal H}$.
Suppose that $\vec w_n(0)\rightharpoonup\mathbf h$ along a subsequence,
and let $w$ be the free wave with data $\mathbf h$. The free solution
map in \eqref{fp:eq:free-strichartz} is bounded and linear, so
\[
 w_n\rightharpoonup w
 \quad\text{in }L^5([0,\infty);L^{10}).
\]
By \eqref{fp:eq:strichartz-tail}, $w_n\to0$ in the same space. Hence
$w=0$ on the half-line and $\mathbf h=0$. This proves
\eqref{fp:free-separation} for every choice of $S_n$.

\begin{proposition}[One-limit rigidity]\label{fp:comparison}
Let $v$ be a free energy wave, let $t_n\to\infty$, and fix an integer
$J\ge0$ and nonzero data $\mathbf f_j\in\mathcal H$, $1\le j\le J$.
Suppose that the similarities $S_{j,n}$ satisfy
\eqref{main:profile-separation} and that
\begin{equation}\label{fp:core-representation}
 \vec v(t_n)+\sum_{j=1}^J S_{j,n}\mathbf f_j
 \longrightarrow\mathbf Q\ne0\quad\text{in }\mathcal H.
\end{equation}
Then $v=0$, $J=1$, and, after extraction,
$S_{1,n}\to S_*$ with $\mathbf Q=S_*\mathbf f_1$.
\end{proposition}
\begin{proof}
We pass to a subsequence on which each $S_{j,n}$ either converges or
escapes. By \eqref{main:profile-separation}, at most one converges,
and $S_{i,n}^{-1}S_{j,n}$ escapes whenever $i\ne j$. Expanding the
square and applying \eqref{fp:dislocation} and \eqref{fp:free-separation}
to the mixed terms, we obtain
\begin{equation}\label{fp:exhaustion}
 \begin{aligned}
 0&=\lim_{n\to\infty}
       \left\|\mathbf Q-\vec v(t_n)
              -\sum_{j=1}^J S_{j,n}\mathbf f_j\right\|_{\mathcal H}^2\\
  &=\left\|\mathbf Q-\sum_{S_{j,n}\to S_j}S_j\mathbf f_j\right\|_{\mathcal H}^2
       +\|\vec v(0)\|_{\mathcal H}^2
       +\sum_{S_{j,n}\text{ escapes}}\|\mathbf f_j\|_{\mathcal H}^2.
 \end{aligned}
\end{equation}
Since $\mathbf Q\ne0$ and $\mathbf f_j\ne0$, we conclude from
\eqref{fp:exhaustion} that
\[
 v=0,\qquad J=1,\qquad S_{1,n}\to S_*,\qquad
 \mathbf Q=S_*\mathbf f_1.\qedhere
\]
\end{proof}

\subsubsection{The two subsequences and the fixed-profile contradiction}\label{fp:energy-shape}
\begin{proof}[Proof of Theorem~\ref{main:counterexample}]
We choose the reference polygon affine on each $I_j$, with
\begin{equation}\label{end:eq:reference-vertices}
 \bar y(0)=0,\qquad
 \bar y(T_j)=\frac{1-(-1)^j}{2}de\quad(j\ge1).
\end{equation}
By Proposition~\ref{construction:response}, we have
\eqref{main:counterexample-bounds}--\eqref{main:scattering}. Since
$e_j=T_{j+1}$, we obtain from \eqref{construction:eq:realization}
\[
 p(T_{2k})\longrightarrow0,\qquad p(T_{2k+1})\longrightarrow(de,0).
\]
By continuity of the chart and Proposition~\ref{geo:prop:shape}, we have
\eqref{main:two-limits}. Suppose \eqref{main:FP}--\eqref{main:profile-separation}
hold on $T_n$. Writing $\mathbf f_j=\vec Q_{j,\ell_j}$, we apply
Proposition~\ref{fp:comparison} first at $T_{2k}$ and then at $T_{2k+1}$.
After extraction, there are $S_{\rm e},S_{\rm o}\in\mathrm{Sim}$ such that
\[
 v_L=v_+,\qquad J=1,\qquad
 (Q,0)=S_{\rm e}\mathbf f_1,\qquad
 (Q_{(de,0)},0)=S_{\rm o}\mathbf f_1 .
\]
Since $S_{\rm o}S_{\rm e}^{-1}\in\mathrm{Sim}$, we reach the contradiction
\[
 0=\operatorname{dist}_{\dot H^1}(Q_{(de,0)},\mathrm{Sim}\cdot Q)
       \gtrsim d>0.\qedhere
\]
\end{proof}

\begin{corollary}[No finite catalogue at all late times]\label{fp:catalogue}
Let $u$ be the solution of Theorem~\ref{main:counterexample}. Even if the
choice of list is allowed to depend on time, no finite catalogue of
fixed lists gives a resolution at all late times. Each list may have its
own fixed free wave and its own finite family of nonzero traveling-wave
profiles with fixed velocities. Within each list the similarities may
vary with time and the separation is as in \eqref{main:profile-separation}
along the times when that list is used.
\end{corollary}
\begin{proof}
The path $t\mapsto Q_{p(t)}$ is continuous and precompact, since its
parameters stay in a closed ball inside the stationary chart. Its limit set
\[
 \Omega=\bigcap_{T>0}\overline{\{Q_{p(t)}:t\ge T\}}^{\,\dot H^1}
\]
is compact and connected, since the sets in the intersection are
nested, nonempty, compact and connected. For the continuous function
\[
 D(q)=\operatorname{dist}_{\dot H^1}(q,\mathrm{Sim}\cdot Q),
\]
the isometry of the similarity action gives $D(Sq)=D(q)$.
By \eqref{main:two-limits}, we have
\[
 [0,D(Q_{(de,0)})]\subset D(\Omega),\qquad D(Q_{(de,0)})\gtrsim d.
\]
For $q\in\Omega$, we choose $t_n\to\infty$ with $Q_{p(t_n)}\to q$.
The compact chart contains only nonzero profiles. If a finite catalogue
sufficed, we could pass to a subsequence using one entry. After we
subtract $\vec v_+(t_n)$, \eqref{main:scattering} and
Proposition~\ref{fp:comparison} give the free wave $v_+$ and a single
profile with first component $g$ such that $q=Sg$. Hence $D(q)=D(g)$.
A finite catalogue would make $D(\Omega)$ finite, a contradiction.
\end{proof}

\subsection{Proof of the one-core sequential capture theorem}\label{proof:capture}
We use the estimates of Section~\ref{cap:start} at one sufficiently
late time in the given sequence.

\begin{proof}[Proof of Theorem~\ref{cap:main}]
We set
\[
 U=u-v,\qquad
 d(t)=\operatorname{dist}_{\mathcal H}(\vec U(t),\mathcal M_\ell(P)).
\]
The orthogonality \eqref{fp:free-separation} holds uniformly on the
compact set of normalized representatives. Testing at $t_n$ against
free waves, we obtain $v=v_{\rm can}$ as in
Section~\ref{cap:sec-radiation}. Let $\eta(t)\downarrow0$ be as in
\eqref{cap:eta}. Lemma~\ref{cap:linear} and the modulation equations
give \eqref{cap:modulation} and \eqref{cap:control} in a fixed tube.
For all sufficiently large $n$, the continuation estimate of
Section~\ref{cap:sec-convergence} gives, uniformly in $n$,
\[
 \sup_{t\ge t_n}d(t)\le C_*\,[d(t_n)+\eta(t_n)]\longrightarrow0 .
\]
The $C^1$ lift in \eqref{cap:gauge} and the estimate
\eqref{cap:modulation} give the decomposition and the three derivative
limits with the original $v$ and $\ell$. Integrating the limits for
$\lambda'$ and $a'$, we also obtain
\[
 \frac{\lambda(t)}{t}\longrightarrow0,\qquad
 \frac{a(t)}{t}\longrightarrow\ell .\qedhere
\]
\end{proof}

For Theorem~\ref{main:counterexample} take $P=Q$, $\ell=0$,
$v=v_+$ and $t_n=T_{2n}$. Theorem~\ref{cert:thm:gap} identifies the
full energy kernel with $T_Q(GQ)$, and \eqref{main:two-limits} supplies
the one-core sequence. The odd limit has positive distance from
$\mathrm{Sim}\cdot Q$. Thus this same solution satisfies the full-time
conclusion while visiting inequivalent similarity classes.

\section{A conformal family with distinct similarity classes}\label{geo:start}
We vary the stereographic pole of a Hopf lift, then normalize the
monopole and dipole coefficients to find a direction transverse to
the similarity orbit. Axial translation and dilation give the other
two chart directions. Our starting point is the balanced single-node
profile of Appendix~\ref{cert:start}, in the isoparametric family of
\cite{FP,NP}. It satisfies
\begin{equation}\label{geo:eq:base}
 w''+\cot r\,w'+\frac{3}{16}(w^5-w)=0,\qquad
 w'(0)=0,\quad w(\frac{\pi}{2})=0,\quad w>0\text{ on }[0,\frac{\pi}{2}).
\end{equation}
We extend it by $w(\pi-r)=-w(r)$ and set
\begin{equation}\label{geo:eq:spherical-profile}
 q(z_1,z_2)=(\frac{3}{4})^{\frac{1}{4}}w(2\chi),\qquad
 z_1=\cos\chi\,e^{i\alpha},\quad z_2=\sin\chi\,e^{i\beta}.
\end{equation}
Thus $q$ satisfies
\[
 \left(-\Delta_{\Sph^3}+\frac{3}{4}\right)q=q^5.
\]
Appendix~\ref{cert:start} gives the kernel classification at $k=2$.

\subsection{The stereographic tail and the stationary chart}\label{geo:chart-section}
We fix a stereographic pole $p\in\Sph^3$ such that $q(p)\ne0$, choose
an orthonormal frame at $p$, and denote the inverse stereographic
projection by $\sigma$. For a smooth function $v$ on the sphere, we have
\begin{equation}\label{geo:eq:tail}
 \begin{gathered}
 \left(\frac{2}{1+|x|^2}\right)^{\frac{1}{2}}v(\sigma(x))
 =\frac{a}{r}+\frac{d\cdot\theta}{r^2}+O(r^{-3}),\\
 a=\sqrt2v(p),\qquad d=2\sqrt2\nabla v(p),\qquad
 r=|x|,\quad\theta=\frac{x}{r}.
 \end{gathered}
\end{equation}
The expansion and its spatial and parameter derivatives are uniform
for smooth compact families.

\begin{proposition}\label{geo:prop:normalization}
Let $R_s$ be a smooth local family of conformal images of $q$, with
coefficients $a_s\ne0,d_s$ in \eqref{geo:eq:tail}. Fix $a_0\ne0$ of
the same sign and put
\begin{equation}\label{geo:eq:normalize}
 Q_s(x)=\frac{a_s}{a_0}
 R_s\left(\frac{a_s^2}{a_0^2}x+\frac{d_s}{a_s}\right).
\end{equation}
Then each $Q_s$ solves $-\Delta Q_s=Q_s^5$, and
\begin{equation}\label{geo:eq:normalized-tails}
\begin{aligned}
 Q_s(x)&=\frac{a_0}{r}+O(r^{-3}),\qquad
 \partial_s^jQ_s(x)=O(r^{-3})\quad(j\ge1),\\
 5Q_s^4(x)&=5a_0^4r^{-4}+O(r^{-6}),\qquad
 \partial_s^j(5Q_s^4)=O(r^{-6})\quad(j\ge1).
\end{aligned}
\end{equation}
In particular, chart differences and tangents belong to $H^2(\R^3)$,
although $Q_s\notin L^2$.
\end{proposition}
\begin{proof}
Substituting \eqref{geo:eq:normalize} into \eqref{geo:eq:tail}, we find
the monopole $a_0$ and the dipole
\[
 \left(\frac{a_0}{a_s}\right)^3\left(d_s-a_s\frac{d_s}{a_s}\right)=0 .
\]
We differentiate in $s$ and expand $5Q_s^4$ to obtain
\eqref{geo:eq:normalized-tails}. For $j\ge1$ and $|\alpha|\le2$, we have
\[
 |\partial_x^\alpha\partial_s^j Q_s(x)|
 \lesssim_{j,\alpha}\langle x\rangle^{-3-|\alpha|}.
\]
The parameter derivatives and their integrals in $s$ therefore belong to $H^2$.
\end{proof}

We choose
\[
 p=(\cos\chi_0,0,\sin\chi_0,0),\qquad 0<\chi_0<\frac{\pi}{4},
\]
and rotate the two real coordinates to move $p$ transversely to the
Clifford tori. The rotations commute with
reflection of either imaginary coordinate. After stereographic
projection, these become the fixed Euclidean reflections
\[
 (x_1,x_2,x_3)\mapsto(-x_1,x_2,x_3),\qquad
 (x_1,x_2,x_3)\mapsto(x_1,-x_2,x_3).
\]
The raw dipole is axial, so normalization uses only translation in
$x_3$. The entire family remains even under both reflections.

\begin{proposition}\label{geo:prop:shape}
We can choose the normalized path so that $K=\partial_sQ_s|_{s=0}$
is nonzero and transverse to the similarity orbit. On a sufficiently
small parameter interval, we have
\[
 \operatorname{dist}_{\dot H^1}(Q_s,\operatorname{Sim}\cdot Q_0)
       \gtrsim |s|.
\]
This assertion is geometric and does not assume kernel completeness.
\end{proposition}
\begin{proof}
Of the six rotations and four boosts, the two block rotations preserve
$q$. The remaining eight tangents displayed in
Section~\ref{cert:count-section} are independent by their angular modes.
Hence the two block rotations span the infinitesimal stabilizer.

At $p$, the two block rotations and the chosen real-plane rotation
have respective values
\[
 (0,\cos\chi_0,0,0),\qquad
 (0,0,0,\sin\chi_0),\qquad
 (-\sin\chi_0,0,\cos\chi_0,0).
\]
The three vectors are independent, and similarity fields vanish at $p$.
The similarity differential is therefore injective, and the chosen
rotation is transverse to its range. Normalization adds only
similarity directions, so $K$ remains nonzero and transverse.

By \eqref{fp:dislocation}, any sequence $S_nR_0\to R_0$ admits,
after extraction,
\[
 S_n\to S_*,\qquad S_*R_0=R_0,\qquad S_nS_*^{-1}\to I.
\]
Since $S_nS_*^{-1}R_0=S_nR_0$, we work near the identity.
By similarity invariance in \eqref{geo:eq:normalize} and a first-order
expansion transverse to the orbit, we obtain
\begin{equation}\label{geo:eq:orbit-distance}
 \operatorname{dist}_{\dot H^1}(Q_s,\operatorname{Sim}\cdot Q_0)
 =\operatorname{dist}_{\dot H^1}(R_s,\operatorname{Sim}\cdot R_0)
 \gtrsim |s|.\qedhere
\end{equation}
\end{proof}

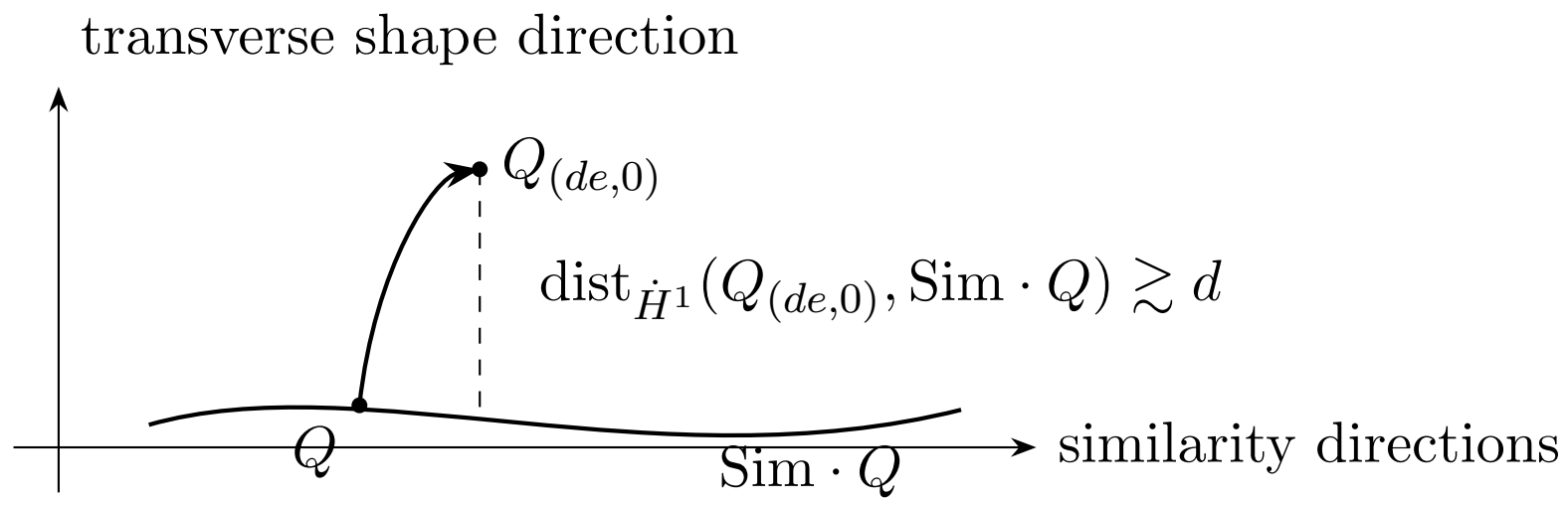
\begin{figure}[htbp]
\centering
\begin{tikzpicture}[>=Stealth,font=\small]
 \draw[->] (-.3,0)--(6.5,0) node[right] {similarity directions};
 \draw[->] (0,-.3)--(0,2.4) node[above right] {transverse shape direction};
 \draw[thick] (.6,.15)..controls(2,.55)and(4,-.25)..(6,.25);
 \node[below] at (5,.15) {$\mathrm{Sim}\cdot Q$};
 \fill (2,.28) circle(1.5pt) node[below left] {$Q$};
 \fill (2.8,1.85) circle(1.5pt) node[right] {$Q_{(de,0)}$};
 \draw[->,thick] (2,.28)..controls(2.1,1.2)and(2.5,1.8)..(2.8,1.85);
 \draw[dashed] (2.8,1.85)--(2.8,.2);
 \node[right] at (3.05,1.05) {$\operatorname{dist}_{\dot H^1}(Q_{(de,0)},\mathrm{Sim}\cdot Q)\gtrsim d$};
\end{tikzpicture}
\caption{Schematic transverse separation, proved in
Proposition~\ref{geo:prop:shape}. Translation, dilation and rotation
remain within one similarity orbit.}
\label{geo:shape-diagram}
\end{figure}

\subsection{The energy kernel in the double-even sector}\label{geo:validation-section}\label{geo:sec:kernel}
We set $Q=Q_0$ and adjoin axial translation and dilation to the shape
parameter. Since $a_0\ne0$, we may use the monopole and axial dipole
as coordinates. Before normalizing the resonance, we have
\begin{equation}\label{geo:eq:chart}
 Q_p=\frac{a_0+a}{r}+\frac{y_2\theta_3}{r^2}+O(r^{-3}),\qquad p=(y_1,y_2,a).
\end{equation}
The expansion is uniform under spatial and parameter differentiation.
At $a=0$, we take
\[
 Q_{(y_1,y_2,0)}(x)
 =Q_{y_1}\!\left(x-\frac{y_2}{a_0}e_3\right).
\]
Since the normalized family $Q_s$ has monopole $a_0$ and zero dipole,
this translation gives axial dipole $y_2$. We set
\[
 E_i=\partial_{y_i}Q_p,\qquad R=\partial_aQ_p,\qquad
 \partial_p^\alpha Q_p=O(r^{-3})\quad(|\alpha|\ge2).
\]
A fixed linear adjustment of $y$ by $a$ imposes the resonance
normalization of Section~\ref{lin:start}. The tail remains affine, so
\[
 |Z_p(x)-Z_*(x)|\lesssim |p|\langle x\rangle^{-3},
 \qquad Z_p=\partial_pQ_p.
\]

\begin{proposition}[Analytic stationary input]\label{geo:ass:stationary}
There is a smooth real sign-changing profile $Q$ and a three-dimensional
stationary chart $Q_p$ even separately in $x_1,x_2$ such that
\[
 \ker_{\dot H^1,\mathrm{even}}(-\Delta-5Q^4)
       =\operatorname{span}\{E_1,E_2,R\},\qquad
 \ker_{L^2,\mathrm{even}}(-\Delta-5Q^4)
       =\operatorname{span}\{E_1,E_2\}.
\]
The full spherical kernel has dimension eight and consists of conformal
tangents. At fixed monopole there is a unit finite-coordinate direction
$e\in\R^2$ such that
\[
 \operatorname{dist}_{\dot H^1}(Q_{(de,0)},\mathrm{Sim}\cdot Q)\gtrsim |d|.
\]
\end{proposition}
\begin{proof}
By \eqref{geo:eq:chart} and Proposition~\ref{geo:prop:shape}, the
vectors $E_1,E_2,R$ are independent, and
\[
 E_1,E_2\in L^2,\qquad R\notin L^2.
\]
By Theorem~\ref{cert:thm:gap} and the count in
Section~\ref{cert:count-section}, we have
\begin{align*}
 \dim\operatorname{span}\{E_1,E_2,R\}
 &=3=\dim\ker_{\dot H^1,\mathrm{even}}L,\\
 \dim\operatorname{span}\{E_1,E_2\}
 &=2=\dim\ker_{L^2,\mathrm{even}}L.\qedhere
\end{align*}
\end{proof}

\section{Threshold expansion and decay of the projected wave propagator}\label{lin:start}
We compute the large-time terms of the linearized wave. The zero modes
give a term linear in time, and the resonance and the dipole both
contribute to the constant term. After we project the forcing off these
terms, the response decays at an integrable rate.

\subsection{The selected response and its singular terms}\label{lin:green-section}
For threshold resolvent expansions, we refer to the work of Jensen
and Kato \cite{JK} and Jensen and Nenciu \cite{JN}. Beceanu
\cite{BeceanuThreshold} treats simultaneous eigenstates and resonances
for Schr\"odinger evolution. For wave estimates with a potential, see
Krieger and Schlag \cite{KS}, Beceanu \cite{Beceanu}, and Beceanu and
Goldberg \cite{BeceanuGoldberg}. Here we prove the weighted remainder
and free response estimates for a potential of order $|x|^{-4}$.

We work on the closed subspace of functions even in $x_1$ and $x_2$.
We complexify this space for the resolvent and use its Hilbert-space
operator norm. In the analytic finite-rank formulas, the transpose
denotes the complex bilinear pairing $f^tg=\int fg$. For projections
onto real finite-dimensional spaces, this convention gives the same
matrix as the Hilbert pairing.

We fix the centered stationary state of Proposition~\ref{geo:ass:stationary}.
In this section we suppress the subscript $*$ on the frozen operators.
We write
\begin{equation}\label{fr:eq:weighted-spaces}
 L=-\Delta-V,\quad V=5Q^4,\quad v=\sqrt5Q^2,\quad
 \rho=\langle x\rangle,\quad
 Y_2=L^2(\rho^4\dd x),\quad X_2=L^2(\rho^{-4}\dd x).
\end{equation}
Thus $X_2$ is the dual of $Y_2$ under the integral pairing. We initially choose
$R=r^{-1}+O(r^{-3})$, $T=-\partial_3Q$, and $K=\partial_sQ_s|_0$.
Then
\begin{equation}\label{fr:eq:kernel-basis}
 \ker_{\dot H^1}L=\operatorname{span}\{R,T,K\},\qquad
 E=(T,K),\qquad G_E=E^tE>0,\qquad P_E=EG_E^{-1}E^t.
\end{equation}
We set $f_Z=VZ$. By Gauss' formula and the differentiated stationary
equation, we have
\begin{equation}\label{fr:eq:charges}
 \int f_R=4\pi,\qquad \int f_E=0,\qquad
 d_j:=\int x_jf_E(x)\dd x=4\pi a_0\delta_{j3}(1,0)^t.
\end{equation}
Indeed,
\[
 f_T=-\partial_3(Q^5),\qquad \int Q^5=4\pi a_0,
\]
and the common-tail shape family has constant dipole coefficients.
These integrals converge absolutely by the density decay. We define
\begin{equation}\label{fr:eq:normalR}
 c_i=-\frac{1}{8\pi}\iint f_R(x)|x-y|f_{E_i}(y)\dd x\dd y,
 \qquad R_*=R-EG_E^{-1}c,\qquad D_j=EG_E^{-1}d_j.
\end{equation}
The integral defining $c_i$ is also absolutely convergent. After this
normalization $R_*=r^{-1}+O(r^{-2})$. Its leading charge is unchanged.

By \cite[Section~3.2.1]{DKMcompact}, the negative spectrum of $L$
consists of finitely many eigenvalues $-\nu^2$, counted with multiplicity.
Since the real bounded potential $-V$ lies in $L^{\frac{3}{2}}\cap L^2$,
the spectrum on $(0,\infty)$ is purely absolutely continuous
by \cite[Theorem~1]{GS}. We denote the corresponding projections by $P_\nu$
and $P_c=\mathbf1_{(0,\infty)}(L)$. The selected impulse response is
\begin{equation}\label{fr:eq:mixed}
 S(t)=\mathbf1_{t>0}\left(tP_E+
       \frac{\sin(t\sqrt L)}{\sqrt L}P_c\right)
       -\sum_\nu\frac{e^{-\nu|t|}}{2\nu}P_\nu.
\end{equation}
The positive and zero components are retarded. On the negative modes we use
\[
 (\partial_t^2-\nu^2)\left(-\frac{e^{-\nu|t|}}{2\nu}\right)=\delta_0,
 \qquad (\partial_t^2+L)S=\delta_0I.
\]
The negative-mode kernel decays at both time ends.

\begin{theorem}[Frozen response with all threshold moments]\label{fr:thm:full}
Let
\begin{equation}\label{fr:eq:residue}
 \mathcal A=\frac{1}{4\pi}R_*\otimes R_*
       +\frac{1}{12\pi}\sum_{j=1}^3D_j\otimes D_j.
\end{equation}
Then
\begin{equation}\label{fr:eq:main}
 S(t)=t\mathbf1_{t>0}P_E+\mathbf1_{t>0}\mathcal A+\mathcal K(t),
 \qquad
 \|\mathcal K(t)\|_{Y_2\to X_2}\lesssim \langle t\rangle^{-\frac{3}{2}}.
\end{equation}
Values at the single time $t=0$ may be assigned arbitrarily within the
bound in \eqref{fr:eq:main}.
\end{theorem}

By \eqref{fr:eq:charges} and \eqref{fr:eq:normalR}, we have
\[
 D_1=D_2=0,\qquad D_3=4\pi a_0EG_E^{-1}(1,0)^t\ne0 .
\]
Thus the dipole term in \eqref{fr:eq:residue} has rank one and takes
values in the $L^2$ kernel. The double pole acts on both finite modes,
and the retarded step acts on $R_*$ and $D_3$. We express these
operators in chart coordinates in \eqref{non:eq:residue-coordinates}.
When the source has a tangent component, the cancellation in
\eqref{non:eq:threshold-cancellation} uses both right factors of the step.

\subsection{Free kernels and a local convolution algebra}
We use the Fourier convention
\[
 \widehat f(k)=\int_{\R}e^{ikt}f(t)\dd t,\qquad z=-i(k+i0).
\]
For $\Re z>0$, the free resolvent has kernel
\begin{equation}\label{fr:eq:free-resolvent}
 G_z(x,y)=\frac{e^{-z|x-y|}}{4\pi|x-y|}.
\end{equation}
The boundary values $z^{-1}$ and $z^{-2}$ have inverse transforms
$\mathbf1_{t>0}$ and $t\mathbf1_{t>0}$, respectively.

\begin{lemma}[Weighted free estimates]\label{fr:lem:free}
If $|a(x)|+|b(x)|\lesssim \rho(x)^{-2}$, then
\begin{align}
 \|a\,\mathbf1_{t>0}\sin(tD)D^{-1}b\|_{2\to2}
   &\lesssim\langle t\rangle^{-2},\label{fr:eq:free2}\\
 \|a\cos(tD)b\|_{2\to2}&\lesssim\langle t\rangle^{-2}.
 \label{fr:eq:cosfree}
\end{align}
\end{lemma}
\begin{proof}
We reduce to $a=b=\rho^{-2}$. By Hardy's inequality, duality and
spectral calculus, we have
\[
 \sup_t\left(
 \left\|\frac{\sin(tD)}{D}\rho^{-2}\right\|_{2\to2}
 +\|\cos(tD)\rho^{-2}\|_{2\to2}\right)\lesssim1.
\]
Huygens' principle eliminates both propagators between two copies of
$\mathbf1_{|x|<\frac{|t|}{4}}$. Splitting at this radius and transposing the second exterior factor, we obtain
\begin{equation}\label{fr:eq:schur-free}
 \begin{aligned}
 &\left\|\rho^{-2}\frac{\sin(tD)}{D}\rho^{-2}\right\|_{2\to2}
       +\|\rho^{-2}\cos(tD)\rho^{-2}\|_{2\to2}\\
 &\qquad\lesssim
 \|\mathbf1_{|x|\ge\frac{|t|}{4}}\rho^{-2}\|_\infty
 \lesssim\langle t\rangle^{-2}.
 \end{aligned}
\end{equation}
\end{proof}

We call a multiplier $M$ \emph{regular near zero} if, for a smooth cutoff
$\chi$ equal to one near zero with sufficiently small support, we have
the operator norm estimate
\[
 \|\mathcal F^{-1}(\chi M)(t)\|\lesssim\langle t\rangle^{-\frac{3}{2}} .
\]
Constants are regular after localization. Products are regular because
\[
 \langle t\rangle^{-\frac{3}{2}}*
 \langle t\rangle^{-\frac{3}{2}}
 \lesssim\langle t\rangle^{-\frac{3}{2}}.
\]
Multiplication by a smooth localized
polynomial in $z$ also preserves regularity. Transposes interchange
the dual spaces and preserve these norm bounds.

\begin{lemma}[Local inversion of regular factors]\label{fr:lem:wiener}
If $B(k)$ is regular near zero and $B(0)$ is invertible, its inverse
is regular on a smaller neighborhood. The statement also holds centered
at any fixed real frequency, and with decay exponent $\sigma>1$
in place of $\frac{3}{2}$.
\end{lemma}
\begin{proof}
After localization, we write
\[
 B=\widehat b,\qquad \|b(t)\|\lesssim\langle t\rangle^{-\sigma},
 \qquad B(0)=\int_\R b(s)\dd s.
\]
We choose $\theta\in C_c^\infty$, equal to one near zero, and set
\[
 \widehat h(k)=\theta(\frac{k}{\epsilon})(B(k)-B(0))B(0)^{-1}.
\]
By translation continuity in $L^1$ and dominated convergence, we have
\begin{equation}\label{fr:eq:cutoff-smallness}
 \|h\|_1\lesssim\int_\R\|b(s)\|
 \|\check\theta(\cdot-\epsilon s)-\check\theta\|_1\dd s
 \longrightarrow0.
\end{equation}
We also have
\[
 \|\langle t\rangle^\sigma h\|_\infty<\infty.
\]
In an $n$-fold convolution, one delay has size at least $\frac{|t|}{n}$.
Choosing $\epsilon$ so that $\|h\|_1<1$, we sum the convolution powers to obtain
\begin{equation}\label{fr:eq:convolution-power}
 \sum_{n\ge1}\|\langle t\rangle^\sigma h^{*n}\|_\infty
 \le\|\langle t\rangle^\sigma h\|_\infty
 \sum_{n\ge1}n^{\sigma+1}\|h\|_1^{n-1}<\infty.
\end{equation}
On a smaller neighborhood, we have
\[
 B^{-1}=B(0)^{-1}\left(I+\sum_{n\ge1}(-\widehat h)^n\right).
\]
A cutoff turns the constant term into a Schwartz kernel, and
frequency translation preserves the bounds.
\end{proof}

\subsection{Density cancellations at the threshold}
We estimate the divided resolvent differences in time, using the zero
charges of the finite modes. Throughout the threshold calculation,
$R$ denotes the normalized resonance $R_*$. For $t\ge1$, we use
\begin{equation}\label{fr:eq:tailmoments}
 \int_{|x|>t}|x|^\beta|f_R|\lesssim t^{\beta-2}\ (0\le\beta<2),
 \qquad
 \int_{|x|>t}|x|^\beta|f_E|\lesssim t^{\beta-3}\ (0\le\beta<3).
\end{equation}
The densities are smooth. We use only the moments in
\eqref{fr:eq:tailmoments}.

For $r>0$ and $j\ge1$, we integrate explicitly to obtain
\begin{equation}\label{fr:eq:primitive}
 \frac{e^{-zr}-\sum_{n=0}^{j-1}\frac{(-zr)^n}{n!}}{4\pi r z^j}
 =\frac{(-1)^j}{4\pi(j-1)!r}
       \int_0^r e^{-zt}(r-t)^{j-1}\dd t.
\end{equation}
We omit the fixed constant and sign in the inverse kernel and write
\[
 K_{j,t}(r)=r^{-1}(r-t)^{j-1}\mathbf1_{r>t}.
\]

\begin{lemma}[Scalar and column remainders]\label{fr:lem:remainders}
For $t\ge1$, with or without a fixed smooth frequency cutoff, we have
\begin{equation}\label{fr:eq:remainder-table}
\begin{aligned}
 |f_R^tK_{2,t}f_R|+|f_R^tK_{3,t}f_E|+|f_E^tK_{4,t}f_E|
       &\lesssim t^{-2},\\
 \|K_{1,t}f_R\|_{X_2}+\|K_{2,t}f_E\|_{X_2}
       &\lesssim t^{-\frac{3}{2}}.
\end{aligned}
\end{equation}
The column bound remains true after multiplication by $v$, with target
$L^2$. All the kernels are bounded for $0<t\le1$.
\end{lemma}
\begin{proof}
We first take $t\ge1$. The kernels vanish on $r\le t$, and their
radial extensions satisfy
\[
 0\le K_{2,t}\le1,\qquad
 t\|\nabla K_{2,t}\|_\infty
 +\|\nabla K_{3,t}\|_\infty
 +\|\nabla^2K_{4,t}\|_\infty\lesssim1.
\]
By the mean value formula, we have
\begin{equation}\label{fr:eq:EE-subtraction}
 \begin{aligned}
 K_{2,t}(|x-y|)&\le\mathbf1_{|x|+|y|>t},\\
 |K_{3,t}(|x-y|)-K_{3,t}(|x|)|
 &\lesssim |y|\mathbf1_{|x|+|y|>t},\\
 |K_{4,t}(|x-y|)-K_{4,t}(|x|)-K_{4,t}(|y|)|
 &\lesssim |x||y|\mathbf1_{|x|+|y|>t}.
 \end{aligned}
\end{equation}
By $\int f_E=0$ and \eqref{fr:eq:tailmoments}, we obtain
\[
 \begin{aligned}
 &|f_R^tK_{2,t}f_R|+|f_R^tK_{3,t}f_E|+|f_E^tK_{4,t}f_E|\\
 &\qquad\lesssim\int_{|x|>\frac{t}{2}}
       (|f_R|+\rho|f_E|)\dd x\lesssim t^{-2}.
 \end{aligned}
\]
The subtractions converge absolutely by \eqref{fr:eq:tailmoments}.

For the columns, we use $K_{1,t}\le t^{-1}$ and subtract
$K_{2,t}(|x|)\int f_E=0$. By \eqref{fr:eq:tailmoments}, we obtain
\begin{equation}\label{fr:eq:column-pointwise}
 \begin{aligned}
 |K_{1,t}f_R(x)|+|K_{2,t}f_E(x)|&\lesssim t^{-1},\\
 |x|<\frac{t}{2}\quad\Longrightarrow\quad
 |K_{1,t}f_R(x)|+|K_{2,t}f_E(x)|
 &\lesssim\int_{|y|>\frac{t}{2}}(t^{-1}|f_R|+|f_E|)\dd y
 \lesssim t^{-3}.
 \end{aligned}
\end{equation}
Squaring \eqref{fr:eq:column-pointwise} and integrating against
$\rho^{-4}$, we obtain
\[
 \|K_{1,t}f_R\|_{X_2}^2+\|K_{2,t}f_E\|_{X_2}^2
 \lesssim t^{-6}\|\rho^{-2}\|_2^2
       +t^{-2}\|\mathbf1_{|x|>\frac{t}{2}}\rho^{-2}\|_2^2
 \lesssim t^{-3}.
\]
We use $v\lesssim\rho^{-2}$ for the $L^2$ columns. For $0<t\le1$,
we use local smoothness and
\[
 \int\rho|f_R|+\int\rho^2|f_E|<\infty.
\]
A smooth frequency cutoff convolves these kernels with a Schwartz
function and preserves \eqref{fr:eq:remainder-table}.
\end{proof}

\subsection{Schur reduction and the exact threshold coefficients}
We compute the three-dimensional Schur complement to identify the
coefficients of both poles. We set
\[
 G(k)=\lim_{\epsilon\downarrow0}G_{\epsilon-ik},
 \qquad A(k)=I-vG(k)v,\qquad \phi_Z=vZ.
\]
Thus $G(0)$ is the Newton kernel in \eqref{fr:eq:free-resolvent}.
For real $k\ne0$, the symmetric resolvent formula reads
\begin{equation}\label{fr:eq:BS}
 \mathcal R(k)=\lim_{\epsilon\downarrow0}(L-(k+i\epsilon)^2)^{-1}
       =G(k)+G(k)vA(k)^{-1}vG(k).
\end{equation}
Since $vf\in L^{\frac{6}{5}}$, we have
\[
 A(0)f=0
 \quad\Longleftrightarrow\quad
 Z=G(0)vf\in\dot H^1,\quad LZ=0,\quad f=vZ.
\]
Thus
\[
 \ker A(0)=\operatorname{span}\{\phi_R,\phi_{E_1},\phi_{E_2}\}.
\]
Moreover,
\[
 \iint\frac{v(x)^2v(y)^2}{|x-y|^2}\dd x\dd y
 \lesssim\|v\|_\infty^2\|v\|_2^2+\|v\|_2^4<\infty.
\]
We let $S_0$ project onto this kernel and set $S_0^{\perp}=I-S_0$.
By the Fredholm alternative and Lemmas~\ref{fr:lem:free},
\ref{fr:lem:wiener}, we have
\[
 D(k)=S_0^{\perp}\bigl(S_0^{\perp}A(k)S_0^{\perp}|_{S_0^{\perp}L^2}\bigr)^{-1}S_0^{\perp}
 \quad\text{is regular near zero}.
\]

In the not necessarily orthonormal columns $\phi=(\phi_R,\phi_E)$,
we set
\begin{equation}\label{fr:eq:schur-data}
 m_{ij}=\phi_i^tA\phi_j-(A\phi_i)^tD(A\phi_j),\qquad
 U_j=Gv(\phi_j-DA\phi_j).
\end{equation}
By block Gaussian elimination, we have
\begin{equation}\label{fr:eq:schur}
 \mathcal R=G+GvDvG+Um^{-1}U^t.
\end{equation}
We use the same columns for $m$ and $U$, without orthonormalization.

We expand
\[
 \phi_i^tA(k)\phi_j=f_i^t(G(0)-G(k))f_j.
\]
By \eqref{fr:eq:charges}, \eqref{fr:eq:normalR} and the kernel of
$(-\Delta)^{-2}$, the linear coefficient and the quadratic and cubic
$EE$ coefficients are, respectively,
\[
 \begin{gathered}
 \frac{1}{4\pi}\left(\int f_i\right)\left(\int f_j\right),\\
 -\frac{1}{8\pi}\iint f_E(x)|x-y|f_E(y)^t\dd x\dd y=G_E,\\
 \frac{1}{24\pi}\iint |x-y|^2f_E(x)f_E(y)^t\dd x\dd y
 =-\frac{1}{12\pi}\sum_jd_jd_j^t=:N_1.
 \end{gathered}
\]
The mixed quadratic coefficient vanishes by \eqref{fr:eq:normalR}.
We cancel the $|x|^2,|y|^2$ terms using $\int f_E=0$,
with absolute convergence by \eqref{fr:eq:tailmoments}.

Since $G(0)f_Z=Z$, we have, for $Z=R,E_1,E_2$,
\begin{equation}\label{fr:eq:column-identities}
 A\phi_Z=-v(G-G(0))f_Z,\qquad
 U_Z-Z=(I+GvDv)(G-G(0))f_Z.
\end{equation}
By Lemma~\ref{fr:lem:remainders}, the columns
\begin{equation}\label{fr:eq:divided-columns}
 \frac{A\phi_R}{z},\qquad \frac{A\phi_E}{z^2}
\end{equation}
are regular. Consequently the products $(A\phi_i)^tD(A\phi_j)$ start
at orders $z^2$, $z^3$ and $z^4$ in the $RR$, $RE$ and $EE$ blocks,
respectively. Thus \eqref{fr:eq:remainder-table} gives regularity near
zero for each of the five factors
\begin{equation}\label{fr:eq:factors}
 \frac{m_{RR}-4\pi z}{z^2},\qquad \frac{m_{RE}}{z^3},\qquad
 \frac{m_{EE}-z^2G_E-z^3N_1}{z^4},\qquad
 \frac{U_R-R_*}{z},\qquad \frac{U_E-E}{z^2}.
\end{equation}
The free terms in the first three quotients have time kernels
\[
 f_R^tK_{2,t}f_R,\qquad f_R^tK_{3,t}f_E,\qquad
 f_E^tK_{4,t}f_E,
\]
respectively. We use \eqref{fr:eq:column-identities} and the column
bounds in \eqref{fr:eq:remainder-table} for the last two quotients.
The complement terms are products of $D$ and the columns in
\eqref{fr:eq:divided-columns}.

After scaling the rows of $m$, we have
\begin{equation}\label{fr:eq:schur-inverse}
 \begin{pmatrix}z^{-1}&0\\0&z^{-2}I_2\end{pmatrix}m
 \longrightarrow\begin{pmatrix}4\pi&0\\0&G_E\end{pmatrix}.
\end{equation}
By \eqref{fr:eq:factors}, the scaled matrix is regular and its limit
is invertible. We apply Lemma~\ref{fr:lem:wiener} once and obtain
\begin{equation}\label{fr:eq:inverse-expansion}
 m^{-1}-
 \begin{pmatrix}
 \frac{1}{4\pi z}&0\\
 0&\frac{G_E^{-1}}{z^2}-\frac{G_E^{-1}N_1G_E^{-1}}{z}
 \end{pmatrix}
 \quad\text{is regular near zero}.
\end{equation}
To verify \eqref{fr:eq:inverse-expansion}, we multiply its difference
on the left by the scaled matrix in \eqref{fr:eq:schur-inverse}.
The $z^{-2}$ and $z^{-1}$ terms in the $EE$ block cancel,
and all remaining entries are regular by \eqref{fr:eq:factors}.
Substituting \eqref{fr:eq:inverse-expansion} into \eqref{fr:eq:schur}
and using the last two factors in \eqref{fr:eq:factors}, we conclude that
\begin{equation}\label{fr:eq:laurent}
 \mathcal R(k)-z^{-2}P_E-z^{-1}\mathcal A
 \quad\text{is regular near zero}.
\end{equation}
Indeed, $-EG_E^{-1}N_1G_E^{-1}E^t$ is exactly the dipole term in
\eqref{fr:eq:residue}.

The divided free kernels have continuous Laplace boundary values
by \eqref{fr:eq:remainder-table}. Their products and the inverses with
nonsingular zero values remain continuous up to $\Re z=0$.
Thus the remainder in \eqref{fr:eq:laurent} introduces no distribution
supported at zero. With the prescribed boundary values, we have
\[
 \mathcal F^{-1}(z^{-1})=\mathbf1_{t>0},\qquad
 \mathcal F^{-1}(z^{-2})=t\mathbf1_{t>0},\qquad
 \mathcal F^{-1}\left(-\frac{1}{k^2+\nu^2}\right)
 =-\frac{e^{-\nu|t|}}{2\nu}.
\]
Spectral calculus therefore identifies the boundary resolvent with
the transform of \eqref{fr:eq:mixed}.

\subsection{Middle and high frequencies}\label{lin:frequency-section}
We apply the Wiener argument of \cite[Lemma~2.4]{Beceanu} in the
weighted convolution class of Lemma~\ref{fr:lem:wiener}, first on
compact intervals away from zero, then at high frequencies.

On a compact frequency set away from zero, we apply
\cite[Theorem~1]{GS} and H\"older with $v\in L^4$ to obtain
\[
 \|v\mathcal R(k)v\|_{2\to2}
 \le\|v\|_4^2\|\mathcal R(k)\|_{L^{\frac{4}{3}}\to L^4}\lesssim1,
 \qquad A(k)^{-1}=I+v\mathcal R(k)v.
\]
We pass to the boundary by Hilbert--Schmidt continuity of $vG(k)v$,
then apply Lemma~\ref{fr:lem:wiener} with decay exponent $2$
on a finite frequency cover.
The resolvent formula \eqref{fr:eq:BS} then has a localized
$O(\langle t\rangle^{-2})$ kernel.

For high frequencies, we use
\[
 \widehat b(k)=vG(k)v,\qquad
 b(t)=\mathbf1_{t>0}v\frac{\sin(tD)}{D}v.
\]
By Lemma~\ref{fr:lem:free}, we have
\[
 \|b(t)\|_{2\to2}\lesssim\langle t\rangle^{-2}.
\]
We choose a real even $\chi\in C_c^\infty$, equal to one on $[-1,1]$
and zero outside $[-2,2]$, and set
\[
 b_\Lambda=b-\Lambda\check\chi(\Lambda\cdot)*b,
 \qquad
 \widehat b_\Lambda(k)=(1-\chi(\frac{k}{\Lambda}))vG(k)v.
\]
By translation continuity in $L^1$ and dominated convergence, we have
\begin{equation}\label{fr:eq:highpass-smallness}
 \|b_\Lambda\|_1
 \le\int_\R|\check\chi(s)|
       \|b(\cdot)-b(\cdot-\frac{s}{\Lambda})\|_1\dd s
 \longrightarrow0\qquad(\Lambda\to\infty).
\end{equation}
We fix $\Lambda$ large enough that
\[
 \|b_\Lambda\|_1<\frac{1}{2},\qquad
 \|\langle t\rangle^2b_\Lambda\|_\infty<\infty.
\]
Applying \eqref{fr:eq:convolution-power} with $\sigma=2$, we obtain
\begin{equation}\label{fr:eq:highpass-inverse}
 (I-\widehat b_\Lambda)^{-1}
       =I+\sum_{n\ge1}\widehat b_\Lambda^n,
 \qquad
 \left\|\sum_{n\ge1}b_\Lambda^{*n}(t)\right\|
       \lesssim_\Lambda\langle t\rangle^{-2}.
\end{equation}
With $\psi(k)=1-\chi(\frac{k}{2\Lambda})$, we have
\[
 \widehat b_\Lambda(k)=vG(k)v\quad\text{on }\operatorname{supp}\psi.
\]
From \eqref{fr:eq:BS}, we obtain
\begin{equation}\label{fr:eq:highfrequency-inverse}
 \psi\mathcal R
 =\psi\bigl[G+Gv(I-\widehat b_\Lambda)^{-1}vG\bigr].
\end{equation}
We convolve the free kernels in Lemma~\ref{fr:lem:free} with
\[
 \mathcal F^{-1}\psi=\delta_0-2\Lambda\check\chi(2\Lambda\cdot)
\]
and apply \eqref{fr:eq:highpass-inverse} to obtain
\begin{equation}\label{fr:eq:highfrequency-decay}
 \|\mathcal F^{-1}(\psi\mathcal R)(t)\|_{Y_2\to X_2}
 \lesssim_\Lambda\langle t\rangle^{-2}.
\end{equation}

For the packet construction, we use the high-energy bound
\begin{equation}\label{fr:eq:high-resolvent}
 \|G(k)\|_{Y_2\to X_2}+\|\mathcal R(k)\|_{Y_2\to X_2}
 \lesssim |k|^{-1}\qquad(|k|\gg1).
\end{equation}
The free estimate is the case $\alpha=0$, $\sigma=2$ of
\cite[Eq.~(1.8)]{EGS}. It also follows by integration by parts in
time using \eqref{fr:eq:cosfree}. Since $v\rho^2\in L^\infty$, the same
bound holds for $vGv,Gv,vG$. We then apply \eqref{fr:eq:BS}
and the uniformly bounded inverse of $I-vGv$.

For the high-pass transforms of $z^{-1},z^{-2}$, we subtract a
Schwartz convolution of mass one and vanishing higher moments from
$\mathbf1_{t>0},t\mathbf1_{t>0}$. The remainders are bounded and
rapidly decreasing away from zero.
We combine these bounds with \eqref{fr:eq:laurent} and the
regular-frequency estimates to conclude \eqref{fr:eq:main}.

\subsection{Compatible forcing, spatial derivatives and scalar observations}
We use the chart normalized by \eqref{fr:eq:normalR}. Since its first
two tail coefficients are affine in $p$, we have
\[
 \begin{gathered}
 Z_p=(E_{p,1},E_{p,2},R_p)=\partial_pQ_p,\qquad D_p=Z_p-Z_*,\\
 |D_p(x)|\lesssim |p|\rho^{-3},\qquad
 |\partial_p^\alpha D_p(x)|\lesssim_\alpha\rho^{-3}
       \quad(|\alpha|\ge1).
 \end{gathered}
\]
By the dual-test construction in \cite[(3.23)]{DKMcompact}, followed
by averaging over the reflections, we choose real functions with
\[
 \mu_j\in C_c^\infty,\qquad
 \langle Z_{*,i},\mu_j\rangle=\delta_{ij}.
\]
We set
\begin{equation}\label{lin:eq:projection}
 \Pi f=f-\mu\langle Z_*,f\rangle,
 \qquad \langle Z_p,\mu\rangle=I+O(|p|).
\end{equation}

We define
\begin{equation}\label{lin:eq:spaces}
 \|u\|_X=\sum_{|\alpha|\le2}\|\rho^{-2}\partial^\alpha u\|_2,
 \quad
 \|F\|_Y=\sum_{|\alpha|\le2}\|\rho^2\partial^\alpha F\|_2,
 \quad \|u\|_{\mathcal O}=\|\rho^{-3}u\|_1.
\end{equation}
The projection $\Pi$ is bounded on $Y_2$ and $Y$. Since each right factor in \eqref{fr:eq:residue} lies in
$\operatorname{span}\{Z_*\}$, we have
\begin{equation}\label{lin:eq:annihilation}
 P_E\Pi=0,\qquad\mathcal A\Pi=0,\qquad S_*(t)\Pi=\mathcal K_*(t)\Pi.
\end{equation}
Here
\[
 \Pi^2=\Pi,\qquad
 \operatorname{ran}\Pi=\{f\in Y_2:\langle Z_*,f\rangle=0\}.
\]
We impose the resonance condition by this dual pairing, since
$R_*\notin L^2$. The projection $\Pi$ is not the spectral projection $P_c$.

\begin{proposition}[Compatible estimates]\label{lin:prop:compatible}
For $f\in Y$ and almost every $t$, we have
\begin{equation}\label{lin:eq:regular-X}
 \|\mathcal K_*(t)f\|_X
 \lesssim\langle t\rangle^{-\frac{3}{2}}\|f\|_Y.
\end{equation}
For the compatible response, we also have
\begin{equation}\label{fix:eq:linear}
 \|S_*(t)\Pi f\|_X\lesssim\langle t\rangle^{-\frac{3}{2}}\|f\|_Y.
\end{equation}
\end{proposition}
\begin{proof}
By local elliptic estimates on unit balls with bounded overlap,
we have
\begin{equation}\label{lin:eq:elliptic}
 \|u\|_X\lesssim\|u\|_{X_2}+\|L_*u\|_{X_2}.
\end{equation}
The weight is comparable on each enlarged ball, and $V_*$ is bounded.
Moreover,
\[
 \|L_*f\|_{Y_2}\lesssim\|f\|_Y,\qquad
 \langle Z_*,L_*f\rangle=0,\qquad
 L_*\mathcal K_*(t)f=\mathcal K_*(t)L_*f.
\]
We use cutoff integration by parts for the pairing. The threshold
residues annihilate $L_*f$ and take values in $\ker L_*$, so they
commute with $L_*$. By Theorem~\ref{fr:thm:full} and
\eqref{lin:eq:elliptic}, we have
\[
 \|\mathcal K_*(t)f\|_X
 \lesssim\|\mathcal K_*(t)f\|_{X_2}
       +\|\mathcal K_*(t)L_*f\|_{X_2}
 \lesssim\langle t\rangle^{-\frac{3}{2}}\|f\|_Y.
\]
We use \eqref{lin:eq:annihilation} for \eqref{fix:eq:linear}.
\end{proof}

\begin{lemma}[Free response estimates]\label{lin:pointwise}
Let $z=S_{\rm free}*(F+V_*h)$, with $F\in L^1_tY$
and $h\in L^1_tX$.
Then $z$ has continuous energy, scatters freely, and
\begin{equation}\label{lin:eq:free-response}
 \begin{gathered}
 \|\vec z\|_{L^\infty\mathcal H}+\|z\|_{L^5L^{10}}
       +\|z\|_{L^\infty_{t,x}}
 \lesssim\|F\|_{L^1Y}+\|h\|_{L^1X},\\
 \|\rho z\|_\infty\lesssim\|\rho^5F\|_\infty+\|h\|_\infty.
 \end{gathered}
\end{equation}
The weighted bound holds whenever its right-hand side is finite.
With the time norms of \eqref{time:eq:norms}, we also have
\begin{equation}\label{lin:eq:free-observation}
 \|z_t\|_{\mathcal T(\mathcal O)}
 \lesssim\|F\|_{\mathcal T(Y)}+\|h\|_{\mathcal T(X)}.
\end{equation}
\end{lemma}
\begin{proof}
We apply \eqref{fp:eq:free-strichartz}, also to one spatial derivative,
and use
\[
 \dot H^1\cap\dot H^2\hookrightarrow L^\infty,\qquad
 \|F+V_*h\|_{L^1H^1}
 \lesssim\|F\|_{L^1Y}+\|h\|_{L^1X}
\]
for the first line of \eqref{lin:eq:free-response} and energy continuity.
The Duhamel tail tends to zero in $\mathcal H$, since the source is in
$L^1L^2$.
For the weighted supremum, we use positivity of the free sine
kernel. By Newton's formula, we have
\begin{equation}\label{lin:eq:newton-potential}
 (-\Delta)^{-1}[\rho^{-4}](x)
 =\frac{\arctan|x|}{2|x|}\lesssim\rho(x)^{-1}.
\end{equation}
The quotient has value $\frac{1}{2}$ at zero. Since
$|V_*|\lesssim\rho^{-4}$, we obtain
\[
 |z(t,x)|\lesssim
 (\|\rho^5F\|_\infty+\|h\|_\infty)
                   (-\Delta)^{-1}[\rho^{-4}](x).
\]
For $g\in L^2$, we apply Cauchy--Schwarz on $|x|\le r$ and $|x|>r$, then minimize in $r$ to obtain
\begin{equation}\label{lin:eq:observation-split}
 \|g\|_{\mathcal O}
 \lesssim\inf_{r\ge1}
       \left(r^{\frac{1}{2}}\|g\|_{X_2}
             +r^{-\frac{3}{2}}\|g\|_2\right)
 \lesssim\|g\|_{X_2}^{\frac{3}{4}}\|g\|_2^{\frac{1}{4}}.
\end{equation}
The minimizing radius is at least one, since $\|g\|_{X_2}\le\|g\|_2$.
By \eqref{fr:eq:cosfree} and the unweighted cosine bound, we have
\[
 \|\cos(tD)g\|_{\mathcal O}
 \lesssim k_{\frac{3}{2}}(t)\|g\|_{Y_2}.
\]
We differentiate the free Duhamel formula, use
$\|V_*h\|_{Y_2}\lesssim\|h\|_X$, and apply
Lemma~\ref{time:convolution}. This proves
\eqref{lin:eq:free-observation}.
\end{proof}

\section{Time estimates on separated intervals}\label{lin:old-start}
We use the scales \eqref{main:scales} and the intervals $I_j$ of
Section~\ref{construction:controls-section}. To estimate convolutions
and products on the full time axis, we set
\[
 J_0=[-4T_0,4T_0],\qquad
 J_j=\{T_j\le|t|\le4T_j\}\quad(j\ge1).
\]
For a Banach-valued function, we define
\begin{equation}\label{time:eq:norms}
 \|f\|_{\mathcal T(B)}=\|f\|_{L^\infty B}
 +\|f\|_{L^1B}
 +\eta_0\sup_{j\ge0}\eta_j^{-1}\|f\|_{L^2(J_j;B)}.
\end{equation}
For finite-dimensional targets, we omit $B$. We identify strongly
measurable functions that agree almost everywhere and use a continuous
representative whenever we have proved continuity. We write
$k_s(t)=\langle t\rangle^{-s}$.

\begin{lemma}[Block convolution estimates]\label{time:convolution}
For each $\sigma>1$, convolution by $k_\sigma$ is bounded on
$\mathcal T(B)$ uniformly for $T_0\ge1$.
The same holds for kernels dominated by $k_\sigma$.
If $\|K(t)\|\lesssim k_3(t)$ and
$\|f\|_{L^1(J_j;B)}\le A\eta_j^2$ for every $j$, then
\begin{equation}\label{time:eq:Dconv}
 \|K*f\|_{L^1(J_j;B)}\lesssim A\eta_j^2.
\end{equation}
Bounded instantaneous multipliers may be added. The estimates also
hold for two-time kernels bounded pointwise by the indicated scalar
kernels.
\end{lemma}
\begin{proof}
We apply Young's inequality to the global $L^1$ and $L^\infty$ norms
and to the terms with $|i-j|\le1$. On the other blocks,
$\operatorname{dist}(J_i,J_j)\gtrsim T_j$, so for $p=1,2$ we have
\begin{equation}\label{time:eq:offdiagonal}
 \|K*f\|_{L^p(J_j;B)}
 \lesssim\sum_{|i-j|\le1}\|f\|_{L^p(J_i;B)}
       +T_j^{\frac{1}{p}-\sigma}\|f\|_{L^1B}.
\end{equation}
We take $p=2$ for the $\mathcal T(B)$ bound, and $p=1$, $\sigma=3$
for \eqref{time:eq:Dconv}. Indeed, the neighboring weights are
comparable, $\sum_i\eta_i^2\lesssim\eta_0^2$, and
\[
 \frac{\eta_0}{\eta_j}T_j^{\frac{1}{2}-\sigma}
 =\frac{T_j^{1-\sigma}}{\sqrt{T_0}}\le1,\qquad
 \frac{\eta_0^2}{\eta_j^2}T_j^{-2}
 =\frac{1}{T_0T_j}\le1.
\]
The same bounds hold for two-time kernels and bounded instantaneous
multipliers.
\end{proof}

\begin{lemma}[Product estimates]\label{time:products}
Suppose $f$ and $g$ satisfy, for every $j$,
\[
 \|f\|_{L^2(J_j)}\le A\eta_j,\qquad \|g\|_{L^2(J_j)}\le B\eta_j,\qquad
 \|f\|_\infty\le a,\qquad \|g\|_\infty\le b .
\]
Then
\begin{align}
 \|fg\|_{L^1(J_j)}&\le AB\eta_j^2,\label{time:eq:prodD}\\
 \|fg\|_{\mathcal T}
 &\lesssim ab+AB\eta_0^2+\eta_0\min(aB,bA).
 \label{time:eq:prodT}
\end{align}
\end{lemma}
\begin{proof}
By Cauchy--Schwarz and \eqref{main:scales}, we have
\[
 \begin{aligned}
 \|fg\|_{L^1(J_j)}&\le AB\eta_j^2,&
 \|fg\|_1&\le AB\sum_j\eta_j^2\lesssim AB\eta_0^2,\\
 \|fg\|_{L^2(J_j)}&\le\min(aB,bA)\eta_j,& \|fg\|_\infty&\le ab .
 \end{aligned}
\]
We substitute these bounds in \eqref{time:eq:norms}.
\end{proof}

\section{Radiation packets and their quadratic impulses}\label{cat:start}
To prove Proposition~\ref{construction:controls}, we find finitely
many rays whose impulses generate the required positive cone. We
realize the impulses by exact wave packets and estimate the time
integral of the error in their superposition.

\subsection{The cone of ray impulses}\label{cat:sec:ray}
For an oriented line $y+\R\nu$, we define its line transform and vector
of moments by
\[
 \begin{aligned}
 Xf(\nu,y)&=\int_\R f(y+t\nu)\,dt,\qquad \nu\in\Sph^2,\quad y\in\nu^\perp,\\
 F(\nu,y)&=20X\bigl(Q^3(R,T,K)\bigr)(\nu,y)\in\R^3 .
 \end{aligned}
\]

\begin{theorem}\label{ray:thm:ray}
The positive cone generated by all $F(\nu,y)$ is $\R^3$.
There is a finite family of lines whose moment vectors span $\R^3$
and admit a strictly positive linear combination equal to zero.
The lines may be symmetrized under the two reflections without
changing their vectors.
\end{theorem}
\begin{proof}
We take surface area measure $d\nu$ on $\Sph^2$ and Lebesgue measure
$dy$ on $\nu^\perp$. By Fubini's theorem, we obtain
\begin{equation}\label{ray:eq:riesz}
 \int(XV)^2\,d\nu\,dy
   =2\iint\frac{V(x)V(z)}{|x-z|^2}\,dx\,dz.
\end{equation}
We verify conformal invariance by changes of variables, using
for inversion
\[
 V_I(x)=|x|^{-4}V\!\left(\frac{x}{|x|^2}\right),\qquad
 |Iu-Iv|=\frac{|u-v|}{|u||v|},\qquad
 |\det DI(u)|=|u|^{-6}.
\]
Translations, rotations, dilations and inversion generate the
conformal action. Differentiating \eqref{ray:eq:riesz} along the
conformal tangent $Z$, we obtain
\begin{equation}\label{ray:eq:balance}
 \int (XV)X(Q^3Z)\,d\nu\,dy=0,\qquad Z=R,T,K.
\end{equation}
The differentiated integrals converge absolutely, since
\[
 |V|+|Q^3R|\lesssim\rho^{-4},\qquad
 |Q^3T|\lesssim\rho^{-5},\qquad |Q^3K|\lesssim\rho^{-6}.
\]
Also $XV>0$ on every line, since $Q$ is nonzero for large $|x|$.

For a covector $a$ with $a\cdot F\ge0$, we deduce from
\eqref{ray:eq:balance} and the Fourier slice identity that
\[
 \begin{aligned}
 0=\int(XV)(a\cdot F)\,d\nu\,dy
 &\quad\Longrightarrow\quad X\Bigl(Q^3\sum_Za_ZZ\Bigr)=0\\
 &\quad\Longrightarrow\quad Q^3\sum_Za_ZZ=0
 \quad\Longrightarrow\quad a=0.
 \end{aligned}
\]
Here we use $XV>0$ and continuity and independence of the tangents.
By \eqref{geo:eq:base}--\eqref{geo:eq:spherical-profile}, the nodal set
before conformal transformation is the Clifford torus. Hence $Q\ne0$
on a dense set.

By finite-dimensional separation, we have
\[
 \operatorname{cone}\{F(\nu,y):\nu\in\Sph^2,\ y\in\nu^\perp\}=\R^3.
\]
We add finite positive representations of $\pm e_i$ and discard
zero coefficients. The resulting vectors span $\R^3$ with a positive
zero balance. Reflections preserve their moments because the weights
are double-even.
\end{proof}

\subsection{Exact wave packets}\label{cat:packets-section}
We fix an oriented line $b+\R\omega$, with $b\perp\omega$, and choose
a real Schwartz envelope $g$ with $\|g\|_2=1$ and smooth compact
Fourier support. We first take the high-frequency limit, then
concentrate $g$. We define the real free wave
\[
 v_h^0(t,x)=h\operatorname{Re}\left[e^{-itD}
                   (e^{i\frac{\omega\cdot x}{h}}g(x-b))\right].
\]
For a fixed envelope and sufficiently small $h$, its time frequencies
lie in two compact bands centered at $\pm h^{-1}$, disjoint from zero.

\begin{lemma}[Exact packets realizing ray impulses]\label{cat:lem:packets}
There are real frozen waves $v_h$ satisfying
$(\partial_t^2+L_*)v_h=0$ with the same time-frequency bands and
\[
 \|v_h-v_h^0\|_{L^2_tX_2}\lesssim_g h^2,
 \qquad \|v_h^0\|_{L^2_tX_2}\lesssim_g h.
\]
Each fixed $v_h$ belongs to $L^5_tL^{10}_x$, has uniformly bounded
energy, scatters freely at both time ends, and satisfies
\begin{equation}\label{fix:eq:packetdecay}
 \|v_h(t)\|_X+\sum_{m=1}^2\|\partial_t^m v_h(t)\|_{X_2}
 \lesssim_h k_2(t),\qquad
 \|\rho v_h\|_{L^\infty_{t,x}}<\infty.
\end{equation}
Its time derivative is $O_h(k_{\frac{3}{2}}(t))$ in
$\mathcal O$. If
\begin{equation}\label{cat:eq:generic-direction}
 \omega_1\omega_2\omega_3\ne0,
\end{equation}
we may symmetrize both waves under the two reflections and retain
these estimates in the double-even sector, with
\begin{equation}\label{cat:eq:normalizedray}
 \lim_{g\ {\rm concentrates}}\ \lim_{h\downarrow0}
 h^{-2}\int_{\mathbb R}\mathsf M^{-1}m_0(v_h)(t)\dd t
 =\frac{1}{4}\mathsf M^{-1}X(20Q^3Z_*)(\omega,b).
\end{equation}
In this formula the tangent vector is expressed in the order of
the chart coordinates.
\end{lemma}
\begin{proof}
We recover the half-wave from twice the cosine propagator by a fixed
smooth cutoff in time frequency, translated by $h^{-1}$ and equal
to one on the positive packet band. Its inverse transform has an
$h$-independent Schwartz envelope. By \eqref{fr:eq:cosfree}, we have
\[
 \|h^{-1}v_h^0(t)\|_{X_2}
 \lesssim_g(k_4*k_2)(t)\lesssim_g k_2(t),
\]
uniformly for small $h$.

For this unsymmetrized packet, we use the full-space resolvent.
We apply Lemma~\ref{fr:lem:free} and \eqref{fr:eq:high-resolvent}
on the full space, as in Section~\ref{lin:frequency-section}, and set
\begin{equation}\label{cat:eq:distorted-packet}
 \widehat v_h=\widehat v_h^0+\mathcal R_*(k)V_*\widehat v_h^0.
\end{equation}
A symmetric band cutoff equals one on the source bands. By the
resolvent identity, the high-frequency bound, and Plancherel, we have
\[
 (\partial_t^2+L_*)v_h=0,\qquad
 \|v_h-v_h^0\|_{L^2X_2}
 \lesssim_g h\|v_h^0\|_{L^2X_2}\lesssim_g h^2.
\]
For the kernel $K_h$ on the packet bands, we have
$\|K_h(t)\|_{Y_2\to X_2}\lesssim_h k_2(t)$ by \eqref{fr:eq:highfrequency-decay}.
We convolve with the free packet bound to obtain
\[
 \begin{gathered}
 \|V_*v_h^0(t)\|_{Y_2}\lesssim_g hk_2(t),\\
 \|v_h(t)\|_{X_2}
 \lesssim_h k_2(t)+(k_2*k_2)(t)\lesssim_h k_2(t).
 \end{gathered}
\]
Localized multipliers $k,k^2$ preserve this bound. We use
\eqref{lin:eq:elliptic} and $L_*v_h=-\partial_t^2v_h$ for
\eqref{fix:eq:packetdecay}.

By the free resolvent identity and the outgoing prescription,
we have
\[
 v_h-v_h^0=S_{\rm free}*(V_*v_h),\qquad V_*v_h\in L^1_tY.
\]
By Lemma~\ref{lin:pointwise} and \eqref{fp:eq:free-strichartz},
$v_h$ scatters freely at both time ends, and
\[
 \|\vec v_h\|_{L^\infty\mathcal H}
 +\|v_h\|_{L^5L^{10}}+\|v_h\|_\infty<\infty.
\]
For the free packet,
Kirchhoff's formula controls $|x|\le2|t|$ by free
$\langle t\rangle^{-1}$ decay and $|x|>2|t|$ by Schwartz tails.
From the weighted estimate in the same lemma we obtain
\begin{equation}\label{cat:eq:packet-sup}
 \|\rho v_h\|_\infty
 \lesssim\|\rho v_h^0\|_\infty+\|v_h\|_\infty<\infty.
\end{equation}
By \eqref{lin:eq:observation-split} and \eqref{fix:eq:packetdecay}, we have
\[
 \|\partial_tv_h(t)\|_{\mathcal O}
 \lesssim\|\partial_tv_h(t)\|_{X_2}^{\frac{3}{4}}
             \|\partial_tv_h(t)\|_2^{\frac{1}{4}}
 \lesssim_h k_{\frac{3}{2}}(t).
\]
Since the potential and the free packet are real, we have
\[
 \mathcal R_*(-k)=\overline{\mathcal R_*(k)},\qquad
 \widehat v_h(-k)=\overline{\widehat v_h(k)}.
\]

For the impulse, we factor the difference of the squares and estimate
\[
 \left|\int\bigl[m_0(v_h)-m_0(v_h^0)\bigr]\dd t\right|
 \lesssim
 \|v_h-v_h^0\|_{L^2X_2}
       \bigl(\|v_h\|_{L^2X_2}+\|v_h^0\|_{L^2X_2}\bigr)
 \lesssim_g h^3.
\]
On the fixed Fourier support of the envelope, we expand
\[
 \begin{gathered}
 \left|\frac{\omega}{h}+\zeta\right|
 =h^{-1}+\omega\cdot\zeta+O_g(h),\\
 h^{-1}v_h^0(t,x)
 =g(x-b-t\omega)\cos\!\left(\frac{\omega\cdot x-t}{h}\right)
       +o_{L_x^2}(1),
 \end{gathered}
\]
uniformly for bounded $t$.
Averaging the oscillation on bounded time intervals and using the uniform
$k_2$ bound on the tails, we obtain
\[
 \lim_{h\downarrow0}h^{-2}\int\mathsf M^{-1}m_0(v_h^0)\dd t
 =\frac{1}{4}\mathsf M^{-1}\int |g(z)|^2
 X(20Q^3Z_*)(\omega,b+\operatorname{proj}_{\omega^\perp}z)\dd z.
\]
The factor is $\frac{1}{2}\cdot\frac{10}{20}=\frac{1}{4}$.
We concentrate $|g|^2$ at zero for \eqref{cat:eq:normalizedray}.

Under \eqref{cat:eq:generic-direction}, the four reflected carriers
are distinct up to sign. We take one half the sum of the four
reflections. The cross impulses vanish by oscillatory averaging,
and the diagonal factor is $4(\frac{1}{2})^2=1$.
Distortion commutes with reflections. By continuity of the ray
moments, we may perturb the finite ray catalogue to satisfy
\eqref{cat:eq:generic-direction}, preserving rank and positive balance.
\end{proof}

\subsection{Packet impulses and their integrated tails}\label{cat:paired-section}
We choose disjoint symmetric frequency bands in
Lemma~\ref{cat:lem:packets}, retaining the spanning property and the
strictly positive zero balance of Theorem~\ref{ray:thm:ray}.
For a prescribed impulse $(c,0)$ with $c$ in a bounded subset of $\R^2$,
we use a linear right inverse and add a sufficiently large multiple
of the zero balance to obtain smooth positive coefficients. We denote
by $\mathcal V_c$ the sum of the corresponding packets weighted by the
square roots of these coefficients. By Plancherel, disjoint time bands
give zero integrated cross terms, so
\begin{equation}\label{cat:eq:catalogue}
 \int_\R[\mathsf M^{-1}m_0(\mathcal V_c)(t)]_y\dd t=c,
 \qquad
 \sup_c\sum_{|\alpha|\le1}
       \|\partial_c^\alpha\mathcal V_c(t)\|_X\lesssim k_2(t).
\end{equation}
All constants are fixed before $d$.

For the train in \eqref{cat:eq:train}, the neighboring packet centers
contribute $O(\eta_j)$ to $L^2(J_j)$. Every other center is at
distance $\gtrsim T_j$. For $\sigma=2,\frac{3}{2}$, we obtain
\begin{equation}\label{cat:eq:train-block}
 \left\|\sum_i\eta_i k_\sigma(\cdot-e_i)\right\|_{L^2(J_j)}
 \lesssim\eta_j+\eta_0T_j^{\frac{1}{2}-\sigma}\lesssim\eta_j.
\end{equation}
We sum the outgoing packet states in $\mathcal H$, using
$\sum_i\eta_i=2\eta_0$. By Lemmas~\ref{cat:lem:packets} and
\ref{time:products}, we obtain the wave and force bounds in
Proposition~\ref{construction:controls}.

For one packet, we use \eqref{cat:eq:catalogue} and write
\[
 a(s)=[\mathsf M^{-1}m_0(\mathcal V_c)(s)]_y,\qquad
 \int_\R a(s)\dd s=c,\qquad |a(s)|\lesssim k_2(s)^2.
\]
Translating a half-line by $s$ changes its indicator on an interval
of length $|s|$. Hence
\begin{equation}\label{cat:eq:tail-translation}
 \int_0^\infty\left|
 \int_\R a(s)\bigl(\mathbf1_{t<e_j+s}-\mathbf1_{t<e_j}\bigr)\dd s
 \right|\dd t
 \le\int_\R|s|\,|a(s)|\dd s\lesssim1.
\end{equation}
We sum the diagonal errors using $\sum_j\eta_j^2\lesssim\eta_0^2$.
For $i<k$, we write $E_{ik}(t)$ for the contribution of that pair to
the integrated force in \eqref{cat:eq:integrated-error}. By Tonelli's theorem, we have
\[
 \|E_{ik}\|_{L^1(0,\infty)}
 \lesssim\eta_i\eta_k\int_0^\infty s\,k_2(s-e_i)k_2(s-e_k)\dd s .
\]
Since $e_i\le\frac{e_k}{4}$, we have
\begin{equation}\label{cat:eq:cross-tail}
 \begin{aligned}
 s\,k_2(s-e_i)k_2(s-e_k)
 &\lesssim\frac{k_2(s-e_i)+k_2(s-e_k)}{e_k}\quad(s\ge0),\\
 \int_0^\infty s\,k_2(s-e_i)k_2(s-e_k)\dd s
 &\lesssim e_k^{-1}.
 \end{aligned}
\end{equation}
We sum the cross errors by
\[
 \sum_{i<k}\frac{\eta_i\eta_k}{e_k}
 \lesssim\frac{\eta_0^2}{T_0}.
\]
Together with the diagonal bound, this proves
\eqref{cat:eq:integrated-error}.
For differences, we vary one factor at a time. The packet derivatives
obey the same envelopes, and their integrated mass difference is
exactly $c_j-\widetilde c_j$. These summable bounds also justify the
infinite train.

\section{The modulated core and its radiative correction}\label{non:start}
We fix the polygon in \eqref{construction:eq:polygon} and choose the
packet controls so that the core reaches its prescribed vertices. We
then cancel the threshold terms and estimate the corrected wave and
its quintic remainder for the contraction in Section~\ref{end:start}.

\subsection{Reference path and packet controls}\label{non:sec:field}
We determine the controls from the trial source and velocity by
prescribing the average velocity on each interval. We use
the norms
\begin{equation}\label{non:eq:contraction-spaces}
 \|w\|_{\mathcal B}=\|w\|_{\mathcal T(X)}+\|\rho w\|_\infty,
 \qquad
 \|f\|_{\mathcal Y}=\|f\|_{\mathcal T(Y)}+\|\rho^5f\|_\infty.
\end{equation}
For $I=[0,\infty)$, we take a real double-even source $f=\Pi f$
on $\R$ and a locally absolutely continuous $b:I\to\R^3$, with
\begin{equation}\label{non:eq:trial-ball}
 \|x\|:=\max\left\{
 \eta_0^{-1}\|f\|_{\mathcal Y},\quad
 \eta_0^{-1}\|b'\|_{\mathcal T(I)},\quad
 \sup_j\eta_j^{-1}\|b\|_{L^2(I_j)}\right\}\le\sqrt d,
 \qquad x=(f,b).
\end{equation}
Here $\mathcal T(I)$ is \eqref{time:eq:norms} restricted to $I$.
Since $b\in L^2(I)$ and $b'\in L^1(I)$, we have
\begin{equation}\label{non:eq:velocity-tail}
 b(t)=-\int_t^\infty b'(s)\dd s,\qquad
 \|b\|_\infty\le\eta_0\|x\|.
\end{equation}
This identity also proves completeness of the ball.

We pass from the spectral basis to the chart basis,
\[
 E=(T,K),\qquad
 E_*=(\partial_{y_1}Q_p,\partial_{y_2}Q_p)\big|_{p=0}.
\]
Since $T=-\partial_3Q$ has axial dipole $a_0$, \eqref{geo:eq:chart} gives
\[
 B=\begin{pmatrix}0&a_0^{-1}\\1&0\end{pmatrix},\qquad E_*=EB,\qquad
 G_y=B^tG_E B,\qquad \mathsf M=\operatorname{diag}(G_y,4\pi).
\]
At $a=0$, the adjustment of $y$ by $a$ preserves both columns.
With $d_j^y=B^td_j$, we obtain
\[
 P_E=E_*G_y^{-1}E_*^t,\qquad D_j=E_*G_y^{-1}d_j^y,\qquad
 d_3^y=4\pi(0,1)^t .
\]
The change of basis satisfies
\[
 E_*G_y^{-1}B^tc=EG_E^{-1}c,\qquad
 E_*^tf=0\ \Longleftrightarrow\ E^tf=0 .
\]
Hence both frames give the same normalized resonance in
\eqref{fr:eq:normalR}. For $z=(z_y,z_a)$, the duality
$\langle Z_*,\mu\rangle=I$ gives
\[
 \mathcal A_*\mu\mathsf M z
 =R_*z_a+\frac{1}{12\pi}D_3(d_3^y)^tz_y .
\]
We therefore have a unique constant matrix $A$ satisfying
\begin{equation}\label{non:eq:residue-coordinates}
 Z_*A=\mathcal A_*\mu\mathsf M,
 \qquad P_E\mu\mathsf M z=E_*z_y\quad(z\in\R^3).
\end{equation}
For the force $P(c)$ in Proposition~\ref{construction:controls}, we set
\begin{equation}\label{non:eq:path}
 r(t)=b(t)-\int_t^\infty P(c)(s)\dd s,
 \qquad y(t)=\int_0^t r_y(s)\dd s,
 \qquad p(t)=(y(t),0)+Ar(t).
\end{equation}
We determine the controls from
\begin{equation}\label{non:eq:control-map}
 \eta_j^2c_j=
 \bigl(\operatorname{avg}_{I_j}-\operatorname{avg}_{I_{j+1}}\bigr)
             (b_y-\bar y'-e_c),\qquad j\ge0,
\end{equation}
where $e_c$ is defined in \eqref{cat:eq:integrated-error} and
\[
 \operatorname{avg}_{J}g=\frac{1}{|J|}\int_Jg.
\]
For differences with the same reference polygon, we write
$\Delta$ for subtraction and $\delta=\|x-\widetilde x\|$.

We denote the right side of \eqref{non:eq:control-map}, divided by
$\eta_j^2$, by $\mathcal C_x(c)_j$. The polygon bound and
Cauchy--Schwarz give, independently of the control radius,
\[
 |\operatorname{avg}_{I_j}\bar y'|\lesssim\eta_j^2,\qquad
 |\operatorname{avg}_{I_j} b_y|
 \lesssim |I_j|^{-1/2}\eta_j\|x\| .
\]
For $\|c\|_\infty\le R$, we use \eqref{cat:eq:integrated-error} and
$|I_j|\eta_j^2\simeq d$ to estimate
\[
 \frac{|\operatorname{avg}_{I_j}e_c|}{\eta_j^2}
 \le\frac{\|e_c\|_{L^1}}{|I_j|\eta_j^2}
 \lesssim C_R\frac{\eta_0^2}{d}=C_Rd^5 .
\]
The adjacent block obeys the same bounds up to fixed factors. Hence
\begin{equation}\label{non:eq:control-estimate}
 \begin{aligned}
 \|\mathcal C_x(c)\|_\infty
 &\le C_0(1+d^{-\frac{1}{2}}\|x\|)+C_Rd^5,\\
 \|\mathcal C_x(c)-\mathcal C_{\widetilde x}(\widetilde c)\|_\infty
 &\le C_0d^{-\frac{1}{2}}\delta+C_Rd^5\|c-\widetilde c\|_\infty .
 \end{aligned}
\end{equation}
We first fix $R>4C_0$, then choose $d$ such that
$C_Rd^5\le\min\{R/4,1/2\}$. Since $\|x\|\le\sqrt d$,
$\mathcal C_x$ is a contraction of the closed ball into itself.
The block averages tend to zero, so we sum \eqref{non:eq:control-map}
and use \eqref{cat:eq:integrated-error} to obtain
\[
 \operatorname{avg}_{I_j}(r_y-\bar y')
 =\operatorname{avg}_{I_j}(b_y-\bar y'-e_c)
       -\sum_{i\ge j}\eta_i^2c_i=0.
\]
Since $y(0)=\bar y(0)=0$, for $t\in I_j$ we obtain
\[
 (y-\bar y)(t)=\int_{\inf I_j}^t
 \bigl(b_y-e_c-\operatorname{avg}_{I_j}(b_y-e_c)\bigr)(s)\dd s.
\]
In particular,
\begin{equation}\label{non:eq:vertices}
 y(e_j)=\bar y(e_j),\qquad
 \|y-\bar y\|_{L^\infty(I_j)}
 \lesssim\|b\|_{L^1(I_j)}+\|e_c\|_{L^1(I_j)}.
\end{equation}

For $t<0$, we extend $r$ by $\chi(t)r(0)$, with a fixed smooth
$\chi=1$ near zero and $\chi=0$ on $(-\infty,-1]$, and continue
$y'=r_y$. By \eqref{cat:eq:train-bounds}, we have
\[
 \sup_{t\in I_j}\left|\int_t^\infty P(c)(s)\dd s\right|
 \le\sum_{i\ge j}\|P(c)\|_{L^1(I_i)}\lesssim\eta_j^2.
\]
Using \eqref{non:eq:velocity-tail} and \eqref{non:eq:vertices}, we obtain
\begin{equation}\label{non:eq:path-bounds}
 \begin{gathered}
 \|p\|_\infty\lesssim d,\qquad
 \|(p',r)\|_{L^2(J_j)}\lesssim\sqrt d\,\eta_j,\qquad
 \|(p',r)\|_\infty\lesssim\sqrt d\,\eta_0,\\
 \|r'\|_{\mathcal T}\lesssim(\|x\|+\eta_0)\eta_0,
 \qquad r'=b'+P(c)\quad\text{on }I.
 \end{gathered}
\end{equation}
For the fixed controls, we subtract \eqref{non:eq:control-map},
use \eqref{non:eq:control-estimate}, and absorb the $d^5$ term.
From \eqref{non:eq:vertices} we then obtain
\begin{equation}\label{non:eq:path-differences}
 \begin{gathered}
 \|\Delta c\|_{\ell^\infty}+d^{-1}\|\Delta p\|_\infty
       \lesssim d^{-\frac{1}{2}}\delta,\qquad
 \|\Delta r'\|_{\mathcal T}\lesssim\eta_0\delta,\\
 \|(\Delta p',\Delta r)\|_{L^2(J_j)}\lesssim\eta_j\delta,
 \qquad\|(\Delta p',\Delta r)\|_\infty\lesssim\eta_0\delta.
 \end{gathered}
\end{equation}

\subsection{The retarded core and the threshold cancellation}
\label{non:retarded-section}
We cancel the threshold terms and estimate the remaining wave
by the free estimates of Lemma~\ref{lin:pointwise}.
We define
\begin{equation}\label{non:eq:core-pair}
 \binom{Q_p+H}{\partial_t(Q_p+H)}(t)
 =\binom{Q_{p(t)}}{E_{p(t)}r_y(t)}
 -\int_{-\infty}^t U_0(t-s)
       \binom{Z_{p(s)}Ar'(s)}{(E_{p(s)}r_y(s))'}\dd s.
\end{equation}
By \eqref{non:eq:path-bounds} and Lemma~\ref{time:products}, we have
\begin{equation}\label{non:eq:core-sources}
 \|Z_pAr'\|_{L^1\dot H^1}
 +\|(E_pr_y)'\|_{L^1L^2}
 \lesssim\|r'\|_1+\|p'r\|_1<\infty.
\end{equation}

We denote the negative Duhamel integral in \eqref{non:eq:core-pair}
by $(H,J)$. By \eqref{non:eq:core-sources}, $(H,J)\in C_t\mathcal H$.
Differentiating in distributions, we obtain
\[
 \partial_t H=J-Z_pAr',\qquad
 \partial_t J=\Delta H-(E_pr_y)' .
\]
Although $\partial_aQ_p$ has a $|x|^{-1}$ tail, the two terms containing
$Z_pAr'$ cancel in the physical velocity,
\[
 \partial_t(Q_p+H)=E_pr_y+Z_pAr'+J-Z_pAr'=E_pr_y+J\in L^2 .
\]
Freezing both columns at $p=0$ defines a second energy pair
$(H_*,J_*)$. Using $-\Delta Q_p=Q_p^5$ and the two equations for
these pairs, we obtain
\begin{equation}\label{non:eq:core-equations}
 \square(Q_p+H)=Q_p^5,\qquad
 \square(Z_*p+H_*)=V_*Z_*p.
\end{equation}
We set
\begin{equation}\label{non:eq:wave-definition}
 F=f+\mu\mathsf M r',\qquad
 w=v_c+\mathcal K_* *F-H_*,\qquad u=Q_p+H+w.
\end{equation}
Since $f=\Pi f$ and $r(-\infty)=0$, we use
\eqref{non:eq:residue-coordinates} to compute
\begin{equation}\label{non:eq:threshold-cancellation}
 (t\mathbf1_{t>0}P_E+\mathbf1_{t>0}\mathcal A_*)*F
 =E_*[y-y(-\infty)]+Z_*Ar
 =Z_*p-E_*y(-\infty).
\end{equation}
Subtracting \eqref{non:eq:core-equations} from the full Green response,
we obtain
\begin{equation}\label{non:eq:wave-identities}
 \begin{gathered}
 (\partial_t^2+L_*)w=F+V_*H_*,\\
 w-v_c=S_{\rm free}*\bigl(F+V_*(\mathcal K_* *F)\bigr).
 \end{gathered}
\end{equation}

To justify the second identity on the full line, we truncate the source,
\[
 \chi_N=\mathbf 1_{[-N,N]},\qquad f_N=\chi_N f,\qquad
 g_N=\chi_N r',\qquad F_N=f_N+\mu\mathsf M g_N .
\]
Then $\Pi f_N=f_N$. In the frozen pair $(H_{*,N},J_{*,N})$ we use
$g_N$ in place of $r'$, and in \eqref{non:eq:threshold-cancellation} we
integrate $g_N$ from $-\infty$ and take a second primitive of its finite
components from the same lower limit. Before the support of
$F_N$, the retarded terms vanish, and the negative modes decay as
$e^{\nu t}$ by \eqref{fr:eq:mixed}. Since their past energy tends to zero,
the free Duhamel formula for $z_N=\mathcal K_* *F_N-H_{*,N}$ gives
\[
 \vec z_N(t)=\int_{-\infty}^t U_0(t-s)(0,G_N(s))\,\dd s,\qquad
 G_N=F_N+V_*(\mathcal K_* *F_N).
\]
We put $G=F+V_*(\mathcal K_* *F)$. Since $F_N\to F$ in $L^1_tY$
and $g_N\to r'$ in $L^1_t$, we have
\[
 \begin{aligned}
 \|\mathcal K_* *(F_N-F)\|_{L^1_tX}
   &\lesssim \|F_N-F\|_{L^1_tY},\\
 \|(H_{*,N}-H_*,J_{*,N}-J_*)\|_{L^\infty_t\mathcal H}
   &\lesssim \|g_N-r'\|_{L^1_t},\\
 \|G_N-G\|_{L^1_tL^2_x}&\lesssim \|F_N-F\|_{L^1_tY}.
 \end{aligned}
\]
Here we used \eqref{lin:eq:regular-X}, the frozen estimate
\eqref{non:eq:core-sources}, and $V_*:X\to L^2$. The free integrals
converge uniformly in $\mathcal H$. Their first components converge
in distributions to $\mathcal K_* *F-H_*$, so the energy limit is
$(w-v_c,\partial_t(w-v_c))$. Its outgoing free data are
\[
 \int_{\R}U_0(-s)(0,G(s))\dd s .
\]
Together with the scattering of $v_c$, this proves the assertion for $w$.
\begin{lemma}[Core and wave bounds]\label{non:lem:retarded}
For trials in \eqref{non:eq:trial-ball}, we have
\begin{align}
 \|H_*\|_{\mathcal B}&\lesssim\|r'\|_{\mathcal T},\notag\\
 \|H-H_*\|_{\mathcal B}
 &\lesssim d\|r'\|_{\mathcal T}+d\eta_0^2,
 \label{non:eq:core-bounds}\\
 \|w-v_c\|_{\mathcal B}+\|\partial_t(w-v_c)\|_{\mathcal T(\mathcal O)}
 +\|\vec w-\vec v_c\|_{L^\infty\mathcal H}
 +\|w-v_c\|_{L^5L^{10}}
 &\lesssim\|f\|_{\mathcal Y}+\|r'\|_{\mathcal T}.
 \label{non:eq:wave-bounds}
\end{align}
The field $w$ has continuous energy and scatters to a free wave.
For $n=H+w$, we also have
\begin{equation}\label{non:eq:field-differences}
 \begin{gathered}
 \|n\|_{\mathcal B}+\|w_t\|_{\mathcal T(\mathcal O)}\lesssim\eta_0,
 \qquad \|n-v_c\|_{\mathcal B}\lesssim(\|x\|+\eta_0)\eta_0,\\
 \|\Delta n\|_{\mathcal B}+\|\Delta w_t\|_{\mathcal T(\mathcal O)}
       \lesssim\eta_0d^{-\frac{1}{2}}\delta,
 \qquad\|\Delta(n-v_c)\|_{\mathcal B}\lesssim\eta_0\delta.
 \end{gathered}
\end{equation}
\end{lemma}
\begin{proof}
We first estimate the core corrections. For $f_{Z,p}=V_pZ_p$,
we use $\int f_{E,p}=0$ in Newton's formula to write
\begin{equation}\label{non:eq:free-columns}
 \begin{aligned}
 \cos(sD)Z_p(x)
 &=\frac{1}{4\pi}\int_{|x-x'|>s}
           \frac{f_{Z,p}(x')}{|x-x'|}\dd x',\\
 S_{\rm free}(s)E_p(x)
 &=-\frac{1}{4\pi}\int
       \left(1-\frac{s}{|x-x'|}\right)_+f_{E,p}(x')\dd x'.
 \end{aligned}
\end{equation}
We transfer spatial derivatives to the densities and apply
Lemma~\ref{fr:lem:remainders}. For $|\alpha|\le2$, we obtain uniformly in $p$
\begin{equation}\label{non:eq:column-bounds}
 \begin{gathered}
 \|\partial_p^\alpha\cos(sD)Z_p\|_X
       +\|\partial_p^\alpha S_{\rm free}(s)E_p\|_X
       \lesssim k_{\frac{3}{2}}(s),\\
 \|\rho\,\partial_p^\alpha\cos(sD)Z_p\|_\infty
       +\|\rho\,\partial_p^\alpha S_{\rm free}(s)E_p\|_\infty
       \lesssim1.
 \end{gathered}
\end{equation}
For the second column, we subtract $(1-\frac{s}{|x|})_+$ in
\eqref{non:eq:free-columns}. The required moments are
\[
 \int\partial_p^\alpha f_{E,p}=0,\qquad
 \int|x|\,|\partial_p^\alpha f_{E,p}|\dd x\lesssim1.
\]
They follow from the affine tails and integration by parts on balls.
Pairing \eqref{fr:eq:column-pointwise} with $f_{Z,q}$ and using
\eqref{fr:eq:tailmoments}, we also have, uniformly in $p,q$,
\begin{equation}\label{non:eq:scalar-columns}
 |\langle f_{Z,q},\cos(sD)Z_p\rangle|
 +|\langle f_{Z,q},S_{\rm free}(s)E_p\rangle|
 \lesssim k_3(s).
\end{equation}
For $s\ge1$, we split at $|x|=\frac{s}{2}$ and bound both pairings by
\[
 s^{-3}\|f_{Z,q}\|_1
 +s^{-1}\|f_{Z,q}\|_{L^1(|x|>\frac{s}{2})}\lesssim s^{-3}.
\]
For $s\le1$, we use
Lemma~\ref{fr:lem:remainders}. The same argument applies after two
parameter derivatives.

By \eqref{non:eq:column-bounds} and Lemma~\ref{time:convolution},
the frozen correction satisfies
the first bound in \eqref{non:eq:core-bounds}. For the difference, the mean value formula and Lemma~\ref{time:products} imply
\[
 \|H-H_*\|_{\mathcal B}
 \lesssim\|p\|_\infty\|r'\|_{\mathcal T}+\|p'r\|_{\mathcal T}
 \lesssim d\|r'\|_{\mathcal T}+d\eta_0^2.
\]
Subtracting \eqref{non:eq:core-pair} for two trials and using
\eqref{non:eq:path-differences}, we obtain
\[
 \begin{aligned}
 \|\Delta(H-H_*)\|_{\mathcal B}
 &\lesssim d\|\Delta r'\|_{\mathcal T}
 +\|\Delta p\|_\infty
       (\|\widetilde r'\|_{\mathcal T}
        +\|\widetilde p'\widetilde r\|_{\mathcal T})
 +\|\Delta(p'r)\|_{\mathcal T}\\
 &\lesssim d\eta_0\delta.
 \end{aligned}
\]

Since $\mathsf M r'=\langle Z_*,F\rangle$, both $H_*$ and $w-v_c$
are linear in $F$. From \eqref{lin:eq:regular-X}, we obtain
\[
 \|\mathcal K_* *F\|_{\mathcal T(X)}+\|H_*\|_{\mathcal B}
 \lesssim\|F\|_{\mathcal Y}
 \lesssim\|f\|_{\mathcal Y}+\|r'\|_{\mathcal T}.
\]

We apply Lemma~\ref{lin:pointwise} with $z=w-v_c$ and $h=\mathcal K_* *F$.
The unweighted estimate gives
\[
 \|z\|_\infty\lesssim\|F\|_{L^1Y}+\|h\|_{L^1X}.
\]
Since $h=z+H_*$ and $\|H_*\|_\infty\le\|H_*\|_{\mathcal B}$,
the weighted estimate gives
\[
 \|\rho z\|_\infty
 \lesssim\|\rho^5 F\|_\infty+\|z\|_\infty+\|H_*\|_\infty
 \lesssim\|F\|_{\mathcal Y}.
\]
Both estimates apply to the limiting field itself.
This proves \eqref{non:eq:wave-bounds}, continuity, scattering, and
\begin{equation}\label{non:eq:linear-field}
 \|\mathcal K_* *F\|_{\mathcal B}\lesssim\|F\|_{\mathcal Y}.
\end{equation}
We use the cancellation
\begin{equation}\label{non:eq:field-cancellation}
 n-v_c=\mathcal K_* *F+(H-H_*).
\end{equation}
For differences, we substitute \eqref{non:eq:field-cancellation} in
\eqref{non:eq:linear-field} and use the wave estimates to obtain
\[
 \begin{gathered}
 \|\Delta F\|_{\mathcal Y}\lesssim\eta_0\delta,\\
 \|\Delta(n-v_c)\|_{\mathcal B}
 +\|\partial_t\Delta(w-v_c)\|_{\mathcal T(\mathcal O)}
 \lesssim\eta_0\delta,\\
 \|\Delta v_c\|_{\mathcal B}
 +\|\partial_t\Delta v_c\|_{\mathcal T(\mathcal O)}
 \lesssim\eta_0\|\Delta c\|_{\ell^\infty}
 \lesssim\eta_0d^{-\frac{1}{2}}\delta.
 \end{gathered}
\]
This proves \eqref{non:eq:field-differences}.
\end{proof}

\subsection{Weighted estimates for the quintic remainder}\label{non:contraction-section}
We expand the quintic nonlinearity about the stationary state and use
the Sobolev product estimate on unit balls. We write
\begin{equation}\label{non:eq:quintic-remainder}
 N_Q(n)=(Q+n)^5-Q^5-5Q^4n
       =10Q^3n^2+n^3(10Q^2+5Qn+n^2).
\end{equation}

\begin{lemma}[Weighted polynomial estimates]\label{non:lem:tame}
Uniformly for $p,\widetilde p$ in the compact stationary chart and
$\|\rho n\|_\infty,\|\rho\widetilde n\|_\infty\le1$, we have
\begin{equation}\label{non:eq:tame}
\begin{aligned}
 \|N_{Q_p}(n)\|_Y&\lesssim \|\rho n\|_\infty\|n\|_X,
 &\|\rho^5N_{Q_p}(n)\|_\infty&\lesssim \|\rho n\|_\infty^2,\\
 \|N_{Q_p}(n)-N_{Q_p}(\widetilde n)\|_Y
 &\lesssim (\|\rho n\|_\infty+\|\rho\widetilde n\|_\infty)
               \|n-\widetilde n\|_X\\
 &\quad+\|\rho(n-\widetilde n)\|_\infty
                  (\|n\|_X+\|\widetilde n\|_X).
\end{aligned}
\end{equation}
We also have
\[
 \|V_p-V_{\widetilde p}\|_{X\to Y}\lesssim|p-\widetilde p|.
\]
If $n_m\to n$ in $X$ and $\sup_m\|\rho n_m\|_\infty<\infty$, then
\[
 N_{Q_p}(n_m)\longrightarrow N_{Q_p}(n)\quad\text{in }Y.
\]
The convergence is uniform on compact sets of strongly converging
$X$-valued paths.
\end{lemma}
\begin{proof}
By \eqref{non:eq:quintic-remainder} and the Sobolev product estimate
\cite[Lemma~A.8]{Tao}, we have, uniformly for $a\in\mathbb Z^3$,
\[
 \|N_{Q_p}(n)\|_{H^2(B(a,1))}
 \lesssim\langle a\rangle^{-4}\|\rho n\|_\infty
                    \|n\|_{H^2(B(a,2))}.
\]
Here $\rho\simeq\langle a\rangle$ on $B(a,2)$, and the first two
spatial derivatives of $Q_p$ are $O(\rho^{-1})$.
We square, multiply by $\langle a\rangle^4$, and sum in $a$ for
the first estimate in \eqref{non:eq:tame}. For differences, we use
\[
 n^j-\widetilde n^j
 =(n-\widetilde n)\sum_{\ell=0}^{j-1}
                  n^\ell\widetilde n^{j-1-\ell},\qquad 2\le j\le5,
\]
and apply the same product estimate. The pointwise bound is immediate
from \eqref{non:eq:quintic-remainder}.
The potential estimate follows from
\[
 |\partial_x^\alpha(V_p-V_{\widetilde p})|
       \lesssim|p-\widetilde p|\rho^{-4},\qquad |\alpha|\le2.
\]
For continuity, we use local $H^2$ convergence and the same bounds
with cutoffs. Uniformly on compact subsets of $X$, we obtain
\[
 \limsup_{m\to\infty}\|N_{Q_p}(n_m)-N_{Q_p}(n)\|_Y
 \lesssim\sum_{|\alpha|\le2}
       \|\rho^{-2}\partial^\alpha n\|_{L^2(|x|>R)}
 \longrightarrow0\quad(R\to\infty).\qedhere
\]
\end{proof}

\section{A contraction for the source and the core velocity}\label{end:start}
The block averages fix the reference vertices. We now solve simultaneously
for the compatible source and the velocity correction. We split the
residual into its projected and tangential equations and use a momentum
identity to control the accumulated velocity error on each interval.

\subsection{The exact residual and its update}\label{end:update-section}
The projected residual determines the new source. In the tangential
equation, we isolate a total time derivative before integrating over
each interval. For the fields in \eqref{non:eq:wave-definition}, we have
\begin{equation}\label{end:eq:residual}
 \begin{gathered}
 S=(V_p-V_*)n+V_*(H-H_*)+N_{Q_p}(n),\\
 \square u-u^5=F-S.
 \end{gathered}
\end{equation}
The identity follows directly from
\eqref{non:eq:core-equations} and \eqref{non:eq:wave-identities}.
We define $\Phi(f,b)=(f^+,b^+)$ by
\begin{equation}\label{end:eq:update}
 f^+=\Pi S\quad\text{on }\R,\qquad
 b^+(t)=b(t)+\int_t^\infty\mathsf M^{-1}
                  \langle Z_{p(s)},F(s)-S(s)\rangle\dd s
 \quad(t\ge0).
\end{equation}
The source belongs to $L^1Y$, so the integral is defined.
At a fixed point, we have
\begin{equation}\label{end:eq:fixed-residual}
 \Pi(F-S)=0,\qquad \langle Z_p,F-S\rangle=0,
 \qquad \langle Z_p,\mu\rangle=I+O(d).
\end{equation}
The first identity gives $F-S\in\operatorname{ran}\mu$. Since
$\langle Z_p,\mu\rangle$ is invertible, the second yields $F=S$,
and \eqref{end:eq:residual} gives $\square u=u^5$ on $I$.

\begin{lemma}[Bounds for the update]\label{non:lem:comparison}
On the ball \eqref{non:eq:trial-ball}, we have
\begin{align}
 \|\Phi(x)\|&\lesssim d+\sqrt d\,\|x\|,
 \label{end:eq:update-size}\\
 \|\Phi(x)-\Phi(\widetilde x)\|
 &\lesssim\sqrt d\,\|x-\widetilde x\|.
 \label{end:eq:update-difference}
\end{align}
We also have
\begin{equation}\label{end:eq:strong-update}
 \|f^+\|_{\mathcal Y}+\|(b^+)'\|_{\mathcal T(I)}\lesssim d\eta_0,
 \qquad
 \sup_j\|b^+\|_{L^1(I_j)}\lesssim d(d+\|x\|).
\end{equation}
\end{lemma}
\begin{proof}
\noindent\textit{Step 1. The source and the differentiated velocity.}
By the core and field estimates in Lemma~\ref{non:lem:retarded}
and the product estimates in Lemma~\ref{non:lem:tame}, we have
\[
 \begin{aligned}
 \|S\|_{\mathcal Y}
 &\lesssim d\|n\|_{\mathcal B}
       +\|H-H_*\|_{\mathcal B}+\|n\|_{\mathcal B}^2
 \lesssim d\eta_0,\\
 \|\Delta S\|_{\mathcal Y}&\lesssim\sqrt d\,\eta_0\delta.
 \end{aligned}
\]
On $I$, we substitute
\[
 \langle Z_*,F\rangle=\mathsf M r',\qquad
 \mathsf M(r'-b')=m_0(v_c)
\]
into the derivative of \eqref{end:eq:update} to obtain
\begin{equation}\label{end:eq:velocity-derivative}
 \mathsf M(b^+)'
 =-m_0(v_c)-\langle D_p,F\rangle+\langle Z_p,S\rangle.
\end{equation}
By $|D_p|\lesssim d\rho^{-3}$, \eqref{cat:eq:train-bounds} and
\eqref{non:eq:path-differences}, we obtain
\[
 \|f^+\|_{\mathcal Y}+\|(b^+)'\|_{\mathcal T(I)}
 \lesssim d\eta_0,
\]
\begin{equation}\label{end:eq:source-difference}
 \|\Delta f^+\|_{\mathcal Y}
 +\|(\Delta b^+)'\|_{\mathcal T(I)}
 \lesssim\sqrt d\,\eta_0\delta.
\end{equation}

\smallskip
\noindent\textit{Step 2. The undifferentiated velocity.}
We test the first equation in \eqref{non:eq:wave-identities}
against $D_p$. Since
\[
 L_*D_p=(V_p-V_*)Z_p,
\]
we have
\[
 \langle(V_p-V_*)Z_p,w\rangle
 =\langle D_p,F+V_*H_*\rangle+\langle\dot D_p,w_t\rangle
                 -\partial_t\langle D_p,w_t\rangle.
\]
Cutoff integration by parts is valid because
$D_p=O(\rho^{-3})$ with differentiated bounds.
After substitution in \eqref{end:eq:velocity-derivative}, we obtain
\begin{equation}\label{end:eq:momentum}
 \partial_t\bigl(\mathsf M b^++\langle D_p,w_t\rangle\bigr)
 =\langle Z_p,N_{Q_p}(n)\rangle-m_0(v_c)
       +\langle\dot D_p,w_t\rangle+\tau,
\end{equation}
where
\[
 \tau=\langle V_pZ_p,H\rangle-\langle V_*Z_*,H_*\rangle.
\]
For the nonlinear term, we use the exact factorization
\begin{equation}\label{non:eq:moment-factorization}
 \begin{aligned}
 \langle Z_p,N_{Q_p}(n)\rangle-m_0(v_c)
 &=10\langle Z_pQ_p^3,(n-v_c)(n+v_c)\rangle\\
 &\quad+10\langle Z_pQ_p^3-Z_*Q^3,v_c^2\rangle\\
 &\quad+\langle Z_p,n^3(10Q_p^2+5Q_pn+n^2)\rangle.
 \end{aligned}
\end{equation}
By the $O(\rho^{-4})$ quadratic weights, Cauchy--Schwarz and
\eqref{non:eq:field-differences}, we obtain both bounds
\[
 \begin{aligned}
 \|\langle Z_p,N_{Q_p}(n)\rangle-m_0(v_c)\|_{L^1(I_j)}
 &\lesssim(d+\|x\|)\eta_j^2,\\
 \|\Delta(\langle Z_p,N_{Q_p}(n)\rangle-m_0(v_c))\|_{L^1(I_j)}
 &\lesssim\eta_j^2\delta.
 \end{aligned}
\]
In particular, when the difference falls on the packet in the first
product of \eqref{non:eq:moment-factorization}, the block bound is
\[
 \|n-v_c\|_{L^2(I_j;X)}
 \|\Delta v_c\|_{L^2(I_j;X)}
 \lesssim (\|x\|+\eta_0)\eta_j\,
             d^{-\frac{1}{2}}\eta_j\delta
 \lesssim \eta_j^2\delta .
\]
For the moving test, we substitute
\[
 p'=(r_y,0)+Ar'
\]
in the force and tail bounds of \eqref{non:eq:path} to obtain
\[
 \begin{aligned}
 \|\langle\dot D_p,w_t\rangle\|_{L^1(I_j)}
 &\lesssim\|(b,b')\|_{L^2(I_j)}\|w_t\|_{L^2(I_j;\mathcal O)}\\
 &\quad+\eta_j^2
       (\|w_t\|_{L^1\mathcal O}+\|w_t\|_{L^\infty\mathcal O})
 \lesssim(\|x\|+\eta_0)\eta_j^2,\\
 \|\Delta\langle\dot D_p,w_t\rangle\|_{L^1(I_j)}
 &\lesssim\eta_j^2\delta.
 \end{aligned}
\]
To estimate $\tau$, we subtract the moving and frozen scalar
columns in \eqref{non:eq:free-columns}. For $z\in\R^3$, let
\begin{equation}\label{end:eq:scalar-difference}
 \begin{aligned}
 \Gamma(t,s)z
 &=\langle V_*Z_*,\cos(sD)Z_*Az+S_{\rm free}(s)E_*z_y\rangle\\
 &\quad-\langle V_{p(t)}Z_{p(t)},
       \cos(sD)Z_{p(t-s)}Az+S_{\rm free}(s)E_{p(t-s)}z_y\rangle.
 \end{aligned}
\end{equation}
By \eqref{non:eq:scalar-columns} and \eqref{non:eq:path-differences},
we have
\[
 |\Gamma(t,s)|\lesssim dk_3(s),\qquad
 |\Delta\Gamma(t,s)|\lesssim\sqrt d\,k_3(s)\delta.
\]
We apply the product rule to $\Gamma(t,s)r(t-s)$ at fixed $s$.
The term $\partial_pE_p[p']r_y$ in \eqref{non:eq:core-pair}
cancels its counterpart from the product rule. For the remaining parameter derivatives, we use
\eqref{non:eq:scalar-columns} to obtain
\begin{equation}\label{end:eq:scalar-primitive}
 \tau=\partial_t\int_0^\infty\Gamma(t,s)r(t-s)\dd s
 +O\bigl(k_3*(|p'||r|)+|p'|(k_3*|r|)\bigr).
\end{equation}
We use Young's inequality on neighboring blocks and integrate
the tail of $k_3$ on the remaining blocks. By
\eqref{non:eq:path-bounds} and \eqref{non:eq:path-differences}, we have
\[
 \begin{aligned}
 \|k_3*|r|\|_{L^2(J_j)}
 &\lesssim\sqrt d\,\eta_j+T_j^{-\frac{3}{2}}\|r\|_\infty
 \lesssim\sqrt d\,\eta_j,\\
 \|k_3*|\Delta r|\|_{L^2(J_j)}&\lesssim\eta_j\delta.
 \end{aligned}
\]
By Lemma~\ref{time:convolution} and Cauchy--Schwarz, also after
varying one factor, we obtain
\[
 \begin{aligned}
 \left\|\tau-\partial_t\int_0^\infty
       \Gamma(t,s)r(t-s)\dd s\right\|_{L^1(I_j)}
 &\lesssim d\eta_j^2,\\
 \left\|\Delta\left(\tau-\partial_t\int_0^\infty
       \Gamma(t,s)r(t-s)\dd s\right)\right\|_{L^1(I_j)}
 &\lesssim\sqrt d\,\eta_j^2\delta.
 \end{aligned}
\]
For the instantaneous terms, we have
\[
 \|\Delta\langle D_p,w_t\rangle\|_{L^2(I_j)}
 +\left\|\Delta\int_0^\infty\Gamma(t,s)r(t-s)\dd s\right\|_{L^2(I_j)}
 \lesssim\sqrt d\,\eta_j\delta.
\]

We set
\[
 \begin{aligned}
 B(t)&=\int_0^\infty\Gamma(t,s)r(t-s)\dd s,\\
 \mathcal F&=\partial_t\bigl(\mathsf M b^++\langle D_p,w_t\rangle-B\bigr).
 \end{aligned}
\]
By \eqref{end:eq:momentum} and \eqref{end:eq:scalar-primitive},
$\mathcal F$ is the sum of the quadratic force difference, the moving
test term, and $\tau-B'$. The preceding estimates give
\[
 \begin{aligned}
 \|\Delta\mathcal F\|_{L^1(I_j)}&\lesssim\eta_j^2\delta,\\
 \left\|\int_t^\infty\Delta\mathcal F(s)\,\dd s\right\|_{L^2_t(I_j)}
 &\le |I_j|^{\frac{1}{2}}
      \sum_{i\ge j}\|\Delta\mathcal F\|_{L^1(I_i)}
 \lesssim T_j^{\frac{1}{2}}\eta_j^2\delta
 =\sqrt d\,\eta_j\delta .
 \end{aligned}
\]
We integrate to infinity using $b^+(t),B(t)\to0$ and a sequence on
which the integrable pairing $\langle D_p,w_t\rangle$ tends to zero.
Adding the instantaneous terms, we obtain
\begin{equation}\label{end:eq:velocity-bound}
 \begin{aligned}
 |b^+(t)|&\lesssim
 d\bigl(\|w_t(t)\|_{\mathcal O}+(k_3*|r|)(t)\bigr)
 +(d+\|x\|)\eta_j^2,\\
 \|\Delta b^+\|_{L^2(I_j)}
 &\lesssim\sqrt d\,\eta_j\delta
       +|I_j|^{\frac{1}{2}}\eta_j^2\delta
 \lesssim\sqrt d\,\eta_j\delta.
 \end{aligned}
\end{equation}
We take block $L^2$ norms and use
$|I_j|^{\frac{1}{2}}\eta_j\simeq\sqrt d$ for
\eqref{end:eq:update-size} and \eqref{end:eq:update-difference}.
For the remaining bound in \eqref{end:eq:strong-update}, we use
$\|w_t\|_{L^1\mathcal O}\lesssim\eta_0$ in
\eqref{end:eq:velocity-bound} to estimate
\[
 \sup_j\|b^+\|_{L^1(I_j)}
 \lesssim d\eta_0+d^2+d(d+\|x\|)
 \lesssim d(d+\|x\|).\qedhere
\]
\end{proof}

\subsection{Solution of the coupled equations}\label{end:fixed-section}
By Lemma~\ref{non:lem:comparison}, we may choose $d$ so that $\Phi$
preserves the ball \eqref{non:eq:trial-ball} and
$\operatorname{Lip}(\Phi)\le1/2$. The contraction mapping theorem
gives a unique fixed point. At this point, we have
\begin{equation}\label{end:eq:fixed-bounds}
 \|x\|\lesssim d,\qquad
 \|f\|_{\mathcal Y}+\|r'\|_{\mathcal T}\lesssim d\eta_0,
 \qquad \sup_j\|b\|_{L^1(I_j)}\lesssim d^2.
\end{equation}
We use \eqref{end:eq:update-size} and absorption for the first bound,
then \eqref{end:eq:strong-update} and $r'=b'+P(c)$ for the others.
The extension before zero obeys the same bounds by
\eqref{non:eq:velocity-tail}. By \eqref{end:eq:fixed-residual},
the corresponding field satisfies $\square u=u^5$ on $I$.
By \eqref{non:eq:vertices}, we also have
\begin{equation}\label{end:eq:path-reference}
 y(e_j)=\bar y(e_j),\qquad
 \|p-(\bar y,0)\|_{L^\infty(I)}\lesssim d^2.
\end{equation}

\section{Energy bounds and the asymptotic free wave}\label{end:global-section}
We now identify the outgoing free wave in
Proposition~\ref{construction:response}. At the fixed point,
\eqref{non:eq:core-sources} and \eqref{end:eq:fixed-bounds} give
\begin{equation}\label{end:eq:integrable-sources}
 \|Z_pAr'\|_{L^1\dot H^1}
 +\|(E_pr_y)'\|_{L^1L^2}\lesssim d\eta_0.
\end{equation}
The Duhamel term in \eqref{non:eq:core-pair} therefore scatters to
a free energy wave. For the remaining velocity term, we have
\[
 \|E_{p(t)}r_y(t)\|_2\lesssim|r(t)|\longrightarrow0.
\]
We add this free state to the outgoing state of $w$ in
Lemma~\ref{non:lem:retarded} and obtain one free wave $v_+$ satisfying
\begin{equation}\label{end:eq:physicalscattering}
 \|\vec u(t)-(Q_{p(t)},0)-\vec v_+(t)\|_{\mathcal H}
 \longrightarrow0.
\end{equation}
At zero, we have
\[
 y(0)=0,\qquad |r(0)|\lesssim d\eta_0+\eta_0^2.
\]
By the free energy and Strichartz estimates in \eqref{non:eq:core-pair}
and by \eqref{non:eq:wave-bounds}, we obtain
\begin{equation}\label{end:eq:initial}
 \|\vec u(0)-(Q,0)\|_{\mathcal H}
 +\|\vec v_+(0)\|_{\mathcal H}
 +\|u-Q_p\|_{L^5(I;L^{10})}\lesssim\eta_0=d^3.
\end{equation}
By the compactness of the stationary chart and
\eqref{end:eq:integrable-sources}, we have
\[
 \vec u\in C(I;\mathcal H),\qquad
 \sup_{t\in I}\|\vec u(t)\|_{\mathcal H}<\infty.
\]
Since $Q_p$ is uniformly in $L^{10}$, we also have
\[
 u\in L^5_{\rm loc}L^{10},\qquad u^5\in L^1_{\rm loc}L^2.
\]
Thus $u$ satisfies the Duhamel formula. By the local theory
\cite[Theorem~2.7 and Remark~2.16]{KM}, $u$ is the solution with data
$\vec u(0)$, and its energy \eqref{main:energy} is conserved.

For $C^1$ regularity of $p$, we use the integrable kernels in
\eqref{non:eq:free-columns}, translation continuity of
$\mathcal K_*$ in $L^1(Y\to X)$ and Lemma~\ref{non:lem:tame} to conclude that
\[
 H,H_*,\mathcal K_* *F\in C_tX,\qquad S\in C_tY.
\]
At the fixed point $F=S$ on $I$, so
\begin{equation}\label{end:eq:path-regularity}
 r'=\mathsf M^{-1}\langle Z_*,S\rangle\in C(I),
 \qquad p'=(r_y,0)+Ar'\in C(I).
\end{equation}
By \eqref{non:eq:path} and \eqref{end:eq:path-reference}, we have
\[
 p(e_j)-(\bar y(e_j),0)=Ar(e_j)\longrightarrow0.
\]
Finally, we sum \eqref{non:eq:path-bounds} in $j$ to obtain
\begin{equation}\label{end:eq:path-l2}
 \|p'\|_{L^2(I)}^2\lesssim d\sum_j\eta_j^2
 \lesssim d\eta_0^2,\qquad \|p'\|_{L^2(I)}\lesssim d^3.
\end{equation}
This completes the proof of Proposition~\ref{construction:response}.

\section{From a one-core sequence to full-time resolution in a conformal class}\label{cap:start}
Krieger, Nakanishi and Schlag proved a one-pass theorem in the
radiative distance to the ground state \cite[Theorem~6.3]{KNS}.
For Theorem~\ref{cap:main}, we identify the radiation by
\eqref{fp:free-separation}. We then bound the growing modes by exterior
energy in Lemma~\ref{cap:linear} and integrate \eqref{cap:modulation}
using \eqref{cap:control}. This yields the late-time estimate
\eqref{cap:late}.

\subsection{The free radiation and the limiting action}\label{cap:sec-radiation}
We fix $P,\ell,u,v$ as in Theorem~\ref{cap:main} and write $U=u-v$.
Constants in this section may depend on $P$ and $\ell$. In stereographic
coordinates, the conformal action is
\[
 (gP)(x)=|\det Dg^{-1}(x)|^{\frac{1}{6}}P(g^{-1}x),\qquad g\in G .
\]
Translations and dilations act simply transitively on hyperbolic
four-space, and the stabilizer of a base point in $G$ is $O(4)$.
Hence every element of $GP$ is a translate and dilate of some
$q\in\mathcal K_P=O(4)P$.
Duyckaerts, Kenig and Merle \cite[Theorem~2 and (1.9)]{DKMradiation}
proved that there is a unique free wave $v_{\rm can}$ for which
\begin{equation}\label{cap:radiation}
 U_0(-t)\vec u(t)\rightharpoonup\vec v_{\rm can}(0),\qquad
 \lim_{t\to\infty}\|\nabla_{t,x}(u-v_{\rm can})(t)\|_{L^2(|x|>t+A)}=0
 \quad(A\in\R).
\end{equation}
By compactness, \eqref{fp:free-separation} holds uniformly for
$\mathbf f\in\{\vec q_\ell:q\in\mathcal K_P\}$. Testing the decomposition
at $t_n$ against $U_0(t_n)\mathbf h$, we obtain, for every $\mathbf h\in\mathcal H$,
\[
 \big\langle\vec v_{\rm can}(0)-\vec v(0),\mathbf h\big\rangle_{\mathcal H}
 =\lim_{n\to\infty}\big\langle\vec U(t_n),U_0(t_n)\mathbf h\big\rangle_{\mathcal H}
 =0.
\]
In particular, $v=v_{\rm can}$. We set
\[
 D_\ell=\ell\cdot\nabla,\qquad A_\ell=(I-\ell\otimes\ell)^{-\frac{1}{2}},\qquad
 \gamma=\frac{1}{\sqrt{1-|\ell|^2}},\qquad \phi_q(y)=q(A_\ell y).
\]
We write the action and its quadratic part as
\[
 \begin{aligned}
 \mathcal J_\ell(f,k)&=E(f,k)+\int kD_\ell f\dd x
       =\mathcal J_{\ell,0}(f,k)-\frac{1}{6}\|f\|_6^6,\\
 2\mathcal J_{\ell,0}(f,k)&=\|k+D_\ell f\|_2^2+\|\nabla f\|_2^2-\|D_\ell f\|_2^2.
 \end{aligned}
\]
Since $-\Delta q=q^5$, we change variables $z=A_\ell y$ to obtain
\[
 \mathcal J_\ell'(\vec q_\ell)=0,\qquad
 \mathcal J_\ell(S_{a,\lambda}\vec q_\ell)
 =\frac{E(q,0)}{\gamma}=\frac{E(P,0)}{\gamma}=:J_* .
\]
We denote by $B_\ell$ the polarization of $\mathcal J_{\ell,0}$ and put
$S(I)=L_t^5L_x^{10}(I\times\R^3)$. By \eqref{cap:radiation}, free
dispersion, density and \eqref{fp:eq:free-strichartz}, we have
\[
 \begin{aligned}
 B_\ell(\vec U(t),\vec v(t))
 &=B_\ell(U_0(-t)\vec U(t),\vec v(0))\longrightarrow0,\\
 \|v(t)\|_6+\|v\|_{S([t,\infty))}&\longrightarrow0,\qquad
 \left|\|u(t)\|_6^6-\|U(t)\|_6^6\right|
 \lesssim \bigl(\|U(t)\|_6^5+\|v(t)\|_6^5\bigr)\|v(t)\|_6 .
 \end{aligned}
\]
By conservation of energy and momentum and the decomposition at $t_n$, we have
\[
 \lim_{t\to\infty}\mathcal J_\ell(\vec U(t))
 =\mathcal J_\ell(\vec u(0))-\mathcal J_{\ell,0}(\vec v(0))
 =\lim_{n\to\infty}\mathcal J_\ell(\vec U(t_n))=J_* .
\]
We shall use the decreasing bound
\begin{equation}\label{cap:eta}
 \eta(t):=\|v\|_{S([t,\infty))}
 +\sup_{s\ge t}\left(\|v(s)\|_6+
       \left|\mathcal J_\ell(\vec U(s))-J_*\right|^{\frac{1}{2}}\right)
 \longrightarrow0 .
\end{equation}

\subsection{The linearized estimate on the exterior cone}\label{cap:sec-exterior}
We use the exterior energy of the linearized wave to control the
growing modes in the kernel of the dual tests.

For $H=(h,k)$ we put $\widetilde h=h\circ A_\ell^{-1}$ and $L_q=-\Delta-5q^4$.
For the Hessian $\mathfrak b_q=\mathcal J_\ell''(\vec q_\ell)$ and
the generator $\mathcal A_q$ in the translating frame, we compute
\begin{equation}\label{cap:hessian}
 \begin{aligned}
 \mathfrak b_q(H,H)
 &=\|k+D_\ell h\|_2^2+\|\nabla h\|_2^2-\|D_\ell h\|_2^2-5\int\phi_q^4 h^2\dd y\\
 &=\|k+D_\ell h\|_2^2+\frac{1}{\gamma}\langle L_q\widetilde h,\widetilde h\rangle,\\
 \mathcal A_qH&=(k+D_\ell h,\Delta h+5\phi_q^4h+D_\ell k).
 \end{aligned}
\end{equation}
By Cauchy--Schwarz, the free form satisfies
\[
 \mathfrak b_{\rm free}(H,H)\ge(1-|\ell|)\|H\|_{\mathcal H}^2 .
\]
By local Rellich compactness and approximation in $L^{\frac{3}{2}}$,
the map $h\mapsto\phi_q^4h$ is compact from $\dot H^1$ to $\dot H^{-1}$.
By full conformal nondegeneracy, the spatial tangent space at $\phi_q$
and the kernel of the Hessian are
\[
 \mathcal N_q=\ker_{\dot H^1}(-\Delta+D_\ell^2-5\phi_q^4),\qquad
 \mathcal Z_q=\ker\mathfrak b_q=\{(Z,-D_\ell Z):Z\in\mathcal N_q\}.
\]
Conformal covariance gives the same finite index for every
$q\in\mathcal K_P$,
\[
 n=\operatorname{ind}_-(L_q)=\operatorname{ind}_-(\mathfrak b_q).
\]
Since $\langle L_q q,q\rangle=-4\int q^6\dd x<0$, a cutoff
approximation gives $n\ge1$. The potentials $5q^4$ vary smoothly in
$L^\infty$ on $\mathcal K_P$. By the min--max principle, the largest of the $n$ negative eigenvalues
is continuous in $q$. Compactness therefore gives a uniform gap from
zero. A local spectral projection gives a smooth orthonormal frame $Y$
for the negative eigenspace, with
\[
 L_qY=-YB^2,\qquad B=B^T\ge\kappa_0I ,
\]
where $\kappa_0>0$ is uniform on $\mathcal K_P$. We boost the modes
$Y(z)e^{\sigma B}$ and set
\begin{equation}\label{cap:modes}
 \begin{aligned}
 \psi(y)&=Y(A_\ell y)e^{-\gamma(\ell\cdot y)B},&
 \pi&=\frac{\psi B}{\gamma}-D_\ell\psi,\\
 V&=(\psi,\pi),\qquad \mathcal E_q=V\R^n,&
 \beta_q(H)&=\int(\psi^Tk-\pi^Th)\dd y .
 \end{aligned}
\end{equation}
For a spectral cluster with $|\ell|\kappa_{\max}<\alpha<\kappa_{\min}$,
we test the eigenfunction equation with $w_M^2Y$, where
$w_M=e^{\alpha\min(|z|,M)}$. We obtain
\[
 \|\nabla(w_MY)\|_2^2+(\kappa_{\min}^2-\alpha^2)\|w_MY\|_2^2
 \le 5\int q^4w_M^2|Y|^2\dd z .
\]
We absorb the integral outside a fixed ball and let $M\to\infty$.
After differentiating in a smooth spectral frame we use a slightly
smaller weight, still larger than $|\ell|\kappa_{\max}$. Thus $V$ and
its shape derivatives decay exponentially, uniformly on a finite cover.

By integration by parts in \eqref{cap:hessian} and \eqref{cap:modes}, we have
\[
 \mathcal A_q V=\frac{VB}{\gamma},\qquad
 \mathfrak b_q(V,H)=\frac{B\beta_q(H)}{\gamma},\qquad
 \beta_q(\mathcal A_qH)=-\frac{B\beta_q(H)}{\gamma}.
\]
The linear flow preserves $\mathfrak b_q$, fixes $\mathcal Z_q$ and
decays backward on $\mathcal E_q$. Thus $\mathcal E_q$ is isotropic,
$\mathcal E_q\cap\mathcal Z_q=\{0\}$ and $\beta_q(\mathcal Z_q)=0$.
By Cauchy uniqueness under the boost, $\dim\mathcal E_q=n$. If an
element of $\mathcal E_q$ has zero negative projection, isotropy places
it in $\mathcal Z_q$. Thus the negative projection maps $\mathcal E_q$
isomorphically onto the negative space of $\mathfrak b_q$. For $H\in\ker\beta_q$, we match its negative projection
with that of $V_0\in\mathcal E_q$. We have
\begin{equation}\label{cap:index}
 \begin{aligned}
 \mathfrak b_q(H,H)&=\mathfrak b_q(H-V_0,H-V_0)\ge0,\\
 \mathfrak b_q(H,H)=0&\quad\Longleftrightarrow\quad
                 H\in\mathcal Z_q+\mathcal E_q .
 \end{aligned}
\end{equation}

For $s\ge0$ and $I=[0,L]$, we set
\[
 \begin{aligned}
 P_q(s,y)&=\phi_q(y-\ell s),& \chi(s,y)&=\mathbf1_{\{|y|>s\}},\\
 e_z(s)&=\|\chi\nabla_{s,y}z(s)\|_2,&
 \|z\|_{\mathcal X(I)}&=\sup_{s\in I}e_z(s)+\|\chi z\|_{S(I)} .
 \end{aligned}
\]
Regularity at the stereographic pole \cite[Section~3.1]{DKMcompact}
gives
\[
 |\nabla^j\phi_q(y)|\lesssim_j\langle y\rangle^{-1-j}.
\]
The normalized tangents and their parameter derivatives obey the same
bound. Since $|y-\ell s|\ge(1-|\ell|)|y|$ on the exterior cone,
we obtain
\[
 \|\chi P_q(s)\|_{10}\lesssim\langle s\rangle^{-\frac{7}{10}},
 \qquad \|\chi P_q\|_{S([0,\infty))}\lesssim1.
\]
For $(\square-5P_q^4)z=F$, we solve $\square g=\chi(5P_q^4z+F)$
with $\vec g(0)=\vec z(0)$. By finite propagation, $\chi g=\chi z$.
On $J=[a,b]\subset I$, we apply H\"older to the potential term in
\eqref{fp:eq:free-strichartz} and obtain
\[
 \|\vec g\|_{L^\infty(J;\mathcal H)}+\|g\|_{S(J)}
 \lesssim\|\vec g(a)\|_{\mathcal H}+\|\chi F\|_{L^1_s(J)L^2_y}
                 +\|\chi P_q\|_{S(J)}^4\|g\|_{S(J)} .
\]
We split $I$ into a uniformly bounded number of intervals on which
$\|\chi P_q\|_{S(J)}$ is small, absorb the last term and obtain
\begin{equation}\label{cap:cone}
 \|z\|_{\mathcal X(I)}\lesssim\|\vec z(0)\|_{\mathcal H}
                            +\|\chi F\|_{L^1_s(I)L^2_y}.
\end{equation}

\begin{lemma}\label{cap:linear}
Let $z$ solve $(\square-5P_q^4)z=F$ with data
$\vec z(0)=H=(h,k)$ satisfying $h\perp_{\dot H^1}\mathcal N_q$.
Writing $[a]_+=\max(a,0)$, we have
\[
 \|H\|_{\mathcal H}^2\lesssim |\beta_q(H)|^2+[\mathfrak b_q(H,H)]_+
 +\|\chi F\|_{L^1_s L^2_y}^2+\limsup_{s\to\infty}e_z(s)^2 .
\]
\end{lemma}
\begin{proof}
If the estimate fails, we choose data $H_j$ of unit norm for which the
right-hand side tends to zero. After extraction, we have
\[
 q_j\to q,\qquad H_j\rightharpoonup H=(h,k),\qquad
 h\perp_{\dot H^1}\mathcal N_q,\qquad \beta_q(H)=0 .
\]
By compactness of the potential and \eqref{cap:index}, we have
\[
 0\ge\limsup_j\mathfrak b_{q_j}(H_j,H_j)
 \ge\mathfrak b_q(H,H)+(1-|\ell|)\limsup_j\|H_j-H\|_{\mathcal H}^2\ge0 .
\]
Thus $H_j\to H$ and $H\in\mathcal Z_q+\mathcal E_q$. For the homogeneous
solution $z$ with data $H$, we apply \eqref{cap:cone} to $z_j-z$. Writing
$\mathcal X=\mathcal X([0,\infty))$, we obtain
\[
 \begin{aligned}
 \|z_j-z\|_{\mathcal X}
 &\lesssim\|H_j-H\|_{\mathcal H}+\|\chi F_j\|_{L^1_sL^2_y}
           +\|\chi(P_{q_j}-P_q)\|_S\|\chi z\|_S\longrightarrow0,\\
 \limsup_{s\to\infty}e_z(s)&\le\|z-z_j\|_{\mathcal X}+\limsup_{s\to\infty}e_{z_j}(s)
 \longrightarrow0 .
 \end{aligned}
\]
Here we used dominated convergence in $S$ for the traveling profiles.

In rest coordinates, we write the solution as
\[
 \begin{aligned}
 (\sigma,x)&=\bigl(\gamma(s-\ell\cdot y),\,A_\ell y-\gamma\ell s\bigr),\\
 z^{\rm rest}(\sigma,x)&=Z(x)+\sum_\kappa e^{\kappa\sigma}Y_\kappa(x),
 \qquad L_qY_\kappa=-\kappa^2Y_\kappa .
 \end{aligned}
\]
Here $Z$ is a stationary tangent and the $\kappa$ are distinct. On
$\Omega_L=\{L<r=|x|<2L,\ 1<p=r-\sigma<2\}$, we have
\[
 \partial_\sigma z^{\rm rest}
   =\sum_\kappa\kappa e^{-\kappa p}e^{\kappa r}Y_\kappa(r\omega)
   =\gamma(\partial_s+D_\ell)z,\qquad |y|>s,\quad s\asymp_\ell L .
\]
Suppose some $Y_\kappa\ne0$. By \cite[Corollary~3.10]{DKMcompact},
the positive Gram matrix of $\{e^{-\kappa p}\}$ on $(1,2)$ and Lorentz
invariance of volume, we obtain
\[
 \begin{aligned}
 1&\lesssim\frac{1}{L}\sum_\kappa\int_L^{2L}
       r^2 e^{2\kappa r}\|Y_\kappa(r\,\cdot)\|_{L^2(\Sph^2)}^2\dd r\\
  &\lesssim\frac{1}{L}\iint_{\Omega_L}|\partial_\sigma z^{\rm rest}|^2\dd x\dd\sigma
   \lesssim_\ell\sup_{s\asymp_\ell L} e_z(s)^2\longrightarrow0 .
 \end{aligned}
\]
Hence $H\in\mathcal Z_q$. Since $h\perp_{\dot H^1}\mathcal N_q$, we
have $H=0$, contrary to $\|H\|_{\mathcal H}=1$.
\end{proof}

\subsection{Modulation and control of the error}\label{cap:sec-modulation}
When $\vec U(t)$ is close to $\mathcal M_\ell(P)$, we choose the
parameters by projecting $U(t)$ onto the spatial orbit in $\dot H^1$.
At the nearest point $\Phi=S_{a,\lambda}\phi_q$, we have
\begin{equation}\label{cap:gauge}
 \begin{gathered}
 d(t)=\operatorname{dist}_{\mathcal H}(\vec U(t),\mathcal M_\ell(P)),\qquad
 H=(h,k)=S_{a,\lambda}^{-1}\vec U-\vec q_\ell,\\
 h\perp_{\dot H^1}\mathcal N_q,\qquad
 d(t)\le\varepsilon(t):=\|H\|_{\mathcal H}\lesssim d(t).
 \end{gathered}
\end{equation}

We first normalize the center and scale. The normalized phases range
over the compact family $\{\phi_q:q\in\mathcal K_P\}$. By
\eqref{fp:dislocation}, escaping similarities applied to this family
have weak limit zero after extraction. Since a normalized phase is
nonzero, competitors at vanishing distance have convergent
representatives modulo the stabilizer. On a slice transverse to the
stabilizer, the derivative of the orthogonality equations is the
tangent Gram matrix. Constant rank and compactness give a uniform
bound for its inverse. Applying the implicit function theorem on a
finite cover, we obtain a unique nearest orbit point in a uniform
tube. Translations and dilations transport this tube isometrically.
For every other orbit position $\Psi$, we have
\[
 \|(\Phi-\Psi,-D_\ell(\Phi-\Psi))\|_{\mathcal H}
 \lesssim\|\Phi-\Psi\|_{\dot H^1}\le2\|U-\Psi\|_{\dot H^1}.
\]
This proves the last bound in \eqref{cap:gauge}. For a tangent $Z$,
we use
\[
 \langle h,Z\rangle_{\dot H^1}=\int h(-\Delta Z)\dd y .
\]
Since $U_t\in C_tL^2$ and $-\Delta Z$ and its parameter derivatives
belong to $L^{6/5}\cap L^2$, the projected position is $C^1$.

We write $q=RP$, $R\in O(4)$, and set
\[
 \dd\tau=\frac{\dd t}{\lambda},\qquad \mu=\lambda',\qquad
 \nu=a'-\ell,\qquad \Lambda=\frac{1}{2}+y\cdot\nabla,\qquad b=|\beta_q(H)|.
\]

In the frame translating with $\ell t$, we use the Euclidean metric
on the normalized translation and scale parameters and an invariant
metric on $O(4)$. The kernel of the orbit differential is the
infinitesimal stabilizer. On its orthogonal complement the differential
has a bounded inverse onto its range, uniformly in $R\in O(4)$ by
constant rank. We use this inverse to lift the tangent of the $C^1$
nearest point. Starting from the same representative, we continue the
local lifts on each compact time interval, since $O(4)$ is compact
and the similarity parameters cannot escape. For this $C^1$ lift,
we obtain
\[
 Z=\mu\Lambda\phi_q+\nu\cdot\nabla\phi_q-\partial_\tau\phi_q,\qquad
 \xi:=|\mu|+|\nu|+|\partial_\tau R|\lesssim\|Z\|_{\dot H^1}.
\]
With $\widetilde v(y)=\sqrt\lambda\,v(t,a+\lambda y)$, we compute
\[
 \begin{aligned}
 H_\tau&=\mathcal A_qH+(Z,-D_\ell Z)
       +\mu(\Lambda h,(\Lambda+1)k)+\nu\cdot\nabla H+(0,\mathcal N),\\
 \mathcal N&=(\phi_q+h+\widetilde v)^5-\phi_q^5-5\phi_q^4h,\qquad
 \|\mathcal N\|_{L^{\frac{6}{5}}}\lesssim\varepsilon^2+\eta
 \quad(\varepsilon,\eta\le1).
 \end{aligned}
\]
We differentiate \eqref{cap:gauge} against a tangent basis $Z_i$, move
the transport derivatives onto the tests, and use $\beta_q(Z,-D_\ell Z)=0$.
By \eqref{cap:modes} we obtain
\[
 \begin{aligned}
 \left|\langle Z,Z_i\rangle_{\dot H^1}+\int(k+D_\ell h)(-\Delta Z_i)\dd y\right|
      &\lesssim\varepsilon\xi,\\
 \left|\partial_\tau\beta_q(H)+\frac{B\beta_q(H)}{\gamma}\right|
      &\lesssim\varepsilon\xi+\varepsilon^2+\eta .
 \end{aligned}
\]
The first bound implies $\xi\lesssim\varepsilon+\varepsilon\xi$.
We absorb the last term and use $B\ge\kappa_0I$ in the second bound.
Since $b$ is independent of the orthonormal frame, its upper right
derivative satisfies
\begin{equation}\label{cap:modulation}
 \xi\lesssim\varepsilon,\qquad
 D_\tau^+b+c_0 b\lesssim\varepsilon^2+\eta,\qquad
 c_0=\frac{\kappa_0}{\gamma}.
\end{equation}

We fix a late $t$ and rescale about $(t,a(t))$. Writing
\[
 f^{(t)}(s,y)=\sqrt{\lambda(t)}\,f(t+\lambda(t)s,a(t)+\lambda(t)y),
 \qquad z=U^{(t)}-P_q ,
\]
we have
\[
 \vec z(0)=H,\qquad
 (\square-5P_q^4)z=F=(P_q+z+v^{(t)})^5-P_q^5-5P_q^4z .
\]
With $A=\|\chi v^{(t)}\|_S\le\eta(t)$ and $N_I=\|z\|_{\mathcal X(I)}$, we use
H\"older and \eqref{cap:cone} to estimate
\[
 \begin{aligned}
 \|\chi F\|_{L^1_s(I)L^2_y}&\lesssim A+(N_I+A)^2+(N_I+A)^5,\\
 N_I&\lesssim\varepsilon+\|\chi F\|_{L^1_s(I)L^2_y}.
 \end{aligned}
\]
For $\varepsilon+A$ small, we close the estimate by continuity in $L$ for $I=[0,L]$,
\[
 N_I\lesssim\varepsilon+A,\qquad
 \|\chi F\|_{L^1_sL^2_y}\lesssim A+\varepsilon^2.
\]
At physical time $t'=t+\lambda(t)s$, the cone lies in $|x|>t'-t-|a(t)|$.
By \eqref{cap:radiation} and $|y-\ell s|\ge(1-|\ell|)s$ there, we have
$e_z(s)\to0$. We expand the action at $\vec q_\ell$ and apply
Lemma~\ref{cap:linear} to obtain
\[
 \begin{aligned}
 \left|\mathcal J_\ell(\vec U(t))-J_*-\frac{1}{2}\mathfrak b_q(H,H)\right|
   &\lesssim\varepsilon^3,\\
 \varepsilon^2&\lesssim b^2+\eta(t)^2+\varepsilon^3 .
 \end{aligned}
\]
For $\varepsilon$ sufficiently small, we conclude that
\begin{equation}\label{cap:control}
 \varepsilon(t)\lesssim b(t)+\eta(t).
\end{equation}

\subsection{Convergence after one late time}\label{cap:sec-convergence}
We integrate the mode inequality in rescaled time and use it to
continue the modulation on a half-line.

\begin{proof}[Proof of the late-time estimate \eqref{cap:late} and decay]
We fix $\delta>0$ below the uniform modulation radius and small enough
to absorb the quadratic term. We then choose $T=t_n$ so late that
$d(T)+\eta(T)$ is small compared with $\delta$. On the component of
$\{t\ge T:d(t)<\delta\}$ containing $T$, \eqref{cap:gauge} and
\eqref{cap:control} give
\[
 \varepsilon^2\lesssim\delta(b+\eta).
\]
We absorb the term in $\delta b$ in \eqref{cap:modulation}. For
$T\le t_0\le t$ in this component, we integrate in $\tau$ and use the
monotonicity of $\eta$ to obtain, uniformly in $T$,
\begin{equation}\label{cap:late}
 b(t)\lesssim e^{-c(\tau(t)-\tau(t_0))}b(t_0)+\eta(t_0),
 \qquad c=\frac{c_0}{2}.
\end{equation}

Taking $t_0=T$ in \eqref{cap:late} and using \eqref{cap:gauge} and
\eqref{cap:control}, we find a uniform $C_*\ge1$ such that
\[
 d(t)\le C_*\,[d(T)+\eta(T)]<\frac{\delta}{2}
\]
throughout the component for all sufficiently large $n$. At a finite
right endpoint, continuity would give $d=\delta$, so the component
is $[T,\infty)$. By \eqref{cap:modulation},
$\lambda(t)\lesssim\lambda(T)+t-T$, which implies $\tau(t)\to\infty$.
Applying \eqref{cap:control} and \eqref{cap:late}, we obtain
\[
 \begin{aligned}
 \limsup_{t\to\infty}\varepsilon(t)&\lesssim\eta(t_0)\longrightarrow0
                  \qquad(t_0\to\infty),\\
 |\lambda'(t)|+|a'(t)-\ell|+\lambda(t)\|q'(t)\|_{\dot H^1}
    &\lesssim\xi(t)\longrightarrow0 .
 \end{aligned}
\]
Returning to the variables in \eqref{cap:gauge}, we obtain, for the
same free wave $v$ as at $t_n$,
\[
 \|\vec u(t)-\vec v(t)-S_{a(t),\lambda(t)}\vec q_\ell(t)\|_{\mathcal H}
       =\varepsilon(t)\longrightarrow0 .\qedhere
\]
\end{proof}

\appendix
\section{The Hopf stationary state and its conformal kernel}\label{cert:start}

We construct a balanced single-node profile in the isoparametric family
of Fern\'andez and Petean \cite{FP} by iterating a positive Green
operator between the barriers in \eqref{cert:eq:comparisons}. For the
kernel calculation, we use the angular conformal modes and radial
Wronskian identities of Neves and Pistoia \cite[Sections~4--8]{NP}.
The trial function for angular frequency eight in \eqref{cert:eq:negative-form},
the supersolution for frequency nine in \eqref{cert:eq:positive-form}, and the trace
comparison in Section~\ref{cert:count-section} exclude the remaining
one-sided modes. The surviving multiplicities give $2+2+4=8$,
as required in Theorem~\ref{cap:main}. We then use stereographic
projection to identify the $L^2$ tangents in the double-even sector.
We write
\begin{equation}\label{cert:eq:profile}
 \begin{gathered}
 c=\frac{3}{16},\qquad E(x)=1-x^2,\qquad
 V=c(5W^4-1),\qquad \rho(x)=\frac{1}{2(1-x)},\\
 (EW')'+c(W^5-W)=0,\qquad W(-x)=-W(x).
 \end{gathered}
\end{equation}
For the one-sided nonzero angular modes, we use the Friedrichs
form domain obtained by closing
$C_c^\infty(-1,1)$ in the norm
\[
 \left(\int_{-1}^1(E|h'|^2+\rho|h|^2)\dd x\right)^{\frac{1}{2}}.
\]
For the one-sided problem, we use
\[
 -(Eh')'-Vh+\Lambda\rho h=0
\]
and its quadratic form
\begin{equation}\label{cert:eq:form}
 \mathfrak q_\Lambda[h]=
 \int_{-1}^1(E|h'|^2-V|h|^2+\Lambda\rho|h|^2)\dd x.
\end{equation}
At $x=-1$, a logarithmic cutoff $\chi_\delta$, vanishing for
$x+1\le\delta^2$ and equal to one for $x+1\ge\delta$, satisfies
\begin{equation}\label{cert:eq:endpoint-cutoff}
 \int_{-1}^1E|\chi_\delta'|^2\dd x
 \lesssim\frac{1}{|\log\delta|}\longrightarrow0.
\end{equation}
This permits approximation at $-1$ of smooth functions with nonzero
endpoint value. The logarithmic solution branch has infinite derivative
energy. At $x=1$, the regular branch for
$\Lambda=\ell^2$, $\ell\ge1$, has leading power $(1-x)^{\frac{\ell}{2}}$
and finite form norm there. The negative-power branch has infinite form norm.
\begin{theorem}\label{cert:thm:gap}
There exists a smooth odd solution $W$ of \eqref{cert:eq:profile}
such that $W(x)>0$ on $(0,1]$. Let $w(r)=W(\cos r)$, let $q$ be its
lift in \eqref{geo:eq:spherical-profile}, and let $Q$ be the normalized Euclidean
representative from Section~\ref{geo:chart-section}. Then
\[
 \ker_{H^1(\Sph^3)}\!\left(-\Delta_{\Sph^3}+\tfrac34-5q^4\right)
   =T_q(Gq),\qquad \dim T_q(Gq)=8 .
\]
For $L=-\Delta-5Q^4$, the sector even under both reflections satisfies
\[
 \dim\ker_{\dot H^1,\mathrm{even}} L=3,\qquad
 \dim\ker_{L^2,\mathrm{even}} L=2 .
\]
\end{theorem}

\subsection{Construction of the single-node profile}\label{cert:profile-section}
We construct upper and lower barriers for the Green iteration, then
use the equation at $x=1$ to bound the amplitude. For $b,d>0$ with
$b+d<1$, we define
\begin{equation}\label{cert:eq:comparisons}
 \begin{gathered}
 D_{b,d}(x)=1-bx^2-dx^4,\qquad
 R_{b,d}(x)=\frac{x}{\sqrt{D_{b,d}(x)}},\\
 p=R_{\frac{5}{7},\frac{1}{16}},\qquad s=R_{\frac{35}{48},\frac{2}{35}},\qquad
 P=\frac{p}{p(1)},\quad S=\frac{s}{s(1)}.
 \end{gathered}
\end{equation}
We differentiate
\[
 R_{b,d}'=(1+dx^4)D_{b,d}^{-\frac{3}{2}}
\]
and write the normalized residual as
\begin{equation}\label{cert:eq:residual-derivative}
 \begin{gathered}
 F_{b,d}=\frac{-(ER_{b,d}')'+cR_{b,d}}{R_{b,d}^5},\qquad
 F_{b,d}'(x)=\frac{N_{b,d}(x^2)}{4x^5},\\
 N_{b,d}(z)=48b-35+(80d-5b)z-(21bd+16d^2)z^3+3d^2z^4.
 \end{gathered}
\end{equation}
Subscripts $p,s$ denote the parameter pairs in \eqref{cert:eq:comparisons}.
For $0\le z\le1$, we have
\[
 N_p=- \frac{5}{7}+\frac{10}{7}z-z^3+\frac{3}{256}z^4,
 \qquad z^4\le z^3,
\]
\[
 z^3-\frac{4}{3}z+\frac{16}{27}
 =\left(z-\frac{2}{3}\right)^2\left(z+\frac{4}{3}\right)\ge0.
\]
Consequently
\begin{equation}\label{cert:eq:comparison-signs}
 \begin{aligned}
 N_p&\le- \frac{5}{7}+\frac{253}{432}
    +\left(\frac{10}{7}-\frac{253}{192}\right)z
    \le- \frac{215}{12096}<0,\\
 N_s&=z\left(\frac{311}{336}-\frac{9087}{9800}z^2+\frac{12}{1225}z^3\right),\\
 \frac{311}{336}-\frac{9087}{9800}+\frac{12}{1225}&=\frac{479}{58800}>0.
 \end{aligned}
\end{equation}
The polynomial $\frac{N_s}{z}$ is decreasing. By
\eqref{cert:eq:comparison-signs}, we have
\[
 F_p'<0<F_s'\quad(0<x<1).
\]
We also have
\[
 \left(\log\frac{s}{p}\right)'
 =\frac{x[(b_s-b_p)+2(d_s-d_p)x^2+(b_sd_p-b_pd_s)x^4]}{D_{b_s,d_s}D_{b_p,d_p}}>0\quad(0<x<1).
\]
Here we use
\[
 b_sd_p-b_pd_s>0,\qquad
 (b_s-b_p)+2(d_s-d_p)=\frac{1}{240}>0.
\]

Let $G$ be the inverse of $Lu=-(Eu')'+cu$, with Dirichlet data at
zero and the finite-energy condition at one. We use the coercive form
\[
 \int_0^1(E|u'|^2+c|u|^2)\dd x,\qquad u(0)=0.
\]
For continuous $f\ge0$, we test with the negative part and compare
with $\frac{\|f\|_\infty}{c}$.
Finite energy forces the flux at one to vanish, so
\begin{equation}\label{cert:eq:green-estimates}
 \begin{gathered}
 0\le Gf\le\frac{\|f\|_\infty}{c},\\
 (Gf)'(x)=\frac{1}{1+x}\int_0^1
 (f-cGf)(x+(1-x)t)\dd t.
 \end{gathered}
\end{equation}
By \eqref{cert:eq:green-estimates}, $Gf$ is Lipschitz up to one.
For $f>0$ on $(0,1]$, we apply the maximum principle in the interior.
At $x=1$, we have
\[
 Gf(1)=0\quad\Longrightarrow\quad(Gf)'(1)=\frac{f(1)}{2}>0,
\]
contrary to $Gf\ge0$. Hence $Gf>0$ on $(0,1]$.

For positive functions normalized at one, we write $v\preceq w$
when $\frac{v}{w}$ is nondecreasing, and set
\begin{equation}\label{cert:eq:normalized-iteration}
 Tv=\frac{G(v^5)}{G(v^5)(1)}.
\end{equation}
For the comparison, let $U,H>0$ on $(0,1]$, with zero value at zero
and zero flux at one. We fix
\[
 \gamma=\frac{U(x_0)}{H(x_0)},\qquad 0<x_0<1.
\]
If $L(U-\gamma H)$ changes sign at most once, from negative to
positive, it is nonpositive on $(0,x_0)$ or nonnegative on $(x_0,1)$.
By the maximum principle, we have
\[
 U-\gamma H\le0\text{ on }(0,x_0)\quad\text{or}\quad
 U-\gamma H\ge0\text{ on }(x_0,1).
\]
In either case we obtain
\[
 \left(\frac{U}{H}\right)'(x_0)\ge0.
\]

We apply this comparison to
\[
 \begin{aligned}
 L\bigl(G(v^5)-\gamma G(w^5)\bigr)
 &=w^5\left[\left(\frac{v}{w}\right)^5-\gamma\right],\\
 L\bigl(G(P^5)-\gamma P\bigr)
 &=P^5\bigl[1-\gamma p(1)^4F_p\bigr],\\
 L\bigl(S-\gamma G(S^5)\bigr)
 &=S^5\bigl[s(1)^4F_s-\gamma\bigr].
 \end{aligned}
\]
By the first identity, $T$ preserves $\preceq$. Using $S\preceq P$ and $F_p'<0<F_s'$ in the last two identities,
we obtain
\begin{equation}\label{cert:eq:iteration-order}
 \begin{gathered}
 S\preceq T(S)\preceq T(P)\preceq P,\\
 S\preceq T^{n+1}(P)\preceq T^n(P)\preceq P\quad(n\ge0).
 \end{gathered}
\end{equation}
The iterates decrease pointwise and lie between $S$ and $P$. By
\eqref{cert:eq:green-estimates} and the lower bound
\[
 G((T^nP)^5)(1)\ge G(S^5)(1)>0 ,
\]
they are uniformly Lipschitz. They therefore converge uniformly,
and continuity of $T$ gives a fixed point $v=Tv$. We set
\begin{equation}\label{cert:eq:profile-comparison}
 W=\bigl(cG(v^5)(1)\bigr)^{-\frac{1}{4}}v,\qquad
 LW=cW^5,\qquad
 \frac{p'}{p}\le\frac{W'}{W}\le\frac{s'}{s}.
\end{equation}
In particular, $W'>0$. We bootstrap \eqref{cert:eq:green-estimates}
to obtain smoothness at $x=1$. At $x=0$, regularity of the equation
and uniqueness give the odd extension. At $x=1$, we compute
\[
 \frac{W'(1)}{W(1)}=\frac{c}{2}\bigl(W(1)^4-1\bigr),\qquad
 \frac{R_{b,d}'(1)}{R_{b,d}(1)}=\frac{1+d}{1-b-d}.
\]
By \eqref{cert:eq:profile-comparison}, we obtain
\begin{equation}\label{cert:eq:amplitude}
 \begin{gathered}
 \frac{155}{3}<\frac{3883}{75}\le W(1)^4
       \le\frac{19303}{359}<\frac{269}{5},\\
 W(1)S(x)\le W(x)\le W(1)P(x)\quad(0\le x\le1).
 \end{gathered}
\end{equation}
\subsection{Bounds for the angular quadratic forms}\label{cert:bounds-section}
We prove positivity at angular frequency nine by a ground-state transform
and construct a polynomial negative direction at angular frequency eight.

For $0<B<1$ and $0\le z\le1$, we use
\[
 \begin{aligned}
 &(1-b-d)(1-Bz)-(1-B)(1-bz-dz^2)\\
 &\qquad=(1-z)[B-b-d-(1-B)dz],\\
 &\frac{4}{5}(1+d_s)-b_s-2d_s=\frac{19}{8400}>0,\qquad
 b_p+d_p-\frac{3}{4}=\frac{3}{112}>0
 \end{aligned}
\]
with $z=x^2$ to obtain
\[
 S(x)\ge\frac{x}{\sqrt{5-4x^2}},\qquad
 P(x)\le\frac{x}{\sqrt{4-3x^2}}\quad(0\le x\le1).
\]
By \eqref{cert:eq:amplitude} and oddness, we have
\begin{equation}\label{cert:eq:potentials}
 \begin{gathered}
 V_-\le V\le V_+,\\
 V_-=\frac{775x^4}{16(5-4x^2)^2}-c,\qquad
 V_+=\frac{807x^4}{16(4-3x^2)^2}-c.
 \end{gathered}
\end{equation}
For the upper bound, we choose
\[
 \Psi(x)=\frac{(1-x)^{\frac{9}{2}}}{(4-3x^2)^{\frac{13}{6}}}>0
 \quad(-1<x<1).
\]
By logarithmic differentiation, we find
\[
 E\frac{\Psi'}{\Psi}
 =-\frac{9}{2}(1+x)+\frac{13x(1-x^2)}{4-3x^2}.
\]
After one more differentiation, we obtain
\begin{equation}\label{cert:eq:upper-residual}
 \frac{-(E\Psi')'-V_+\Psi+81\rho\Psi}{\Psi}
 =\frac{J(-x)}{2(4-3x^2)^2},
\end{equation}
where, writing $y=1-t$, we have
\begin{equation}\label{cert:eq:quartic-positive}
 \begin{aligned}
 J(t)&=694-936t-365t^2+702t^3-94t^4\\
 &=t^2(1-31y)^2+y^2(91+603y+1055ty)>0
       \quad(0\le t\le1).
 \end{aligned}
\end{equation}
We also have
\[
 J(-t)-J(t)=468t(4-3t^2)\ge0\quad(0\le t\le1).
\]
The residual is continuous and positive on $[-1,1]$.
By the ground-state transform and $V\le V_+$, we obtain
\begin{equation}\label{cert:eq:positive-form}
 \mathfrak q_{81}[h]
 \ge\int_{-1}^1\frac{J(-x)}{2(4-3x^2)^2}|h|^2\dd x>0
 \qquad(h\ne0).
\end{equation}
We first prove this for compactly supported $h$, then pass to the
Friedrichs form domain by closure.
For the lower bound, we take
\begin{equation}\label{cert:eq:trial}
 t=\frac{1-x}{2},\qquad h(x)=(5-4x^2)(t^{16}+3t^{48}).
\end{equation}
The function $h$ vanishes at $x=1$ and belongs to the form domain
by \eqref{cert:eq:endpoint-cutoff}. Its factor $5-4x^2$ cancels the
denominator of $V_-$. We expand in powers of $t$ and use
\[
 \int_0^1t^j\dd t=\frac{1}{j+1}\qquad(j=0,1,2,\ldots).
\]
By $V\ge V_-$, we have
\begin{equation}\label{cert:eq:negative-form}
 \begin{aligned}
 \mathfrak q_{64}[h]
 &\le\int_{-1}^1\bigl[E|h'|^2+(64\rho-V_-)|h|^2\bigr]\dd x\\
 &=\frac{98576}{145299}-6\frac{1186819}{4886310}
       +9\frac{3089990}{47525247}\\
 &<\frac{7}{10}-6\frac{6}{25}+9\frac{1}{15}
       =-\frac{7}{50}<0.
 \end{aligned}
\end{equation}
\subsection{The radial kernel}\label{cert:radial-section}
We use two Wronskian identities to exclude both radial parities.
Throughout this subsection, $0<x<1$ and $W'>0$. The comparison
construction gives
\[
 \frac{xW'}{W}\ge\frac{xP'}{P}
 =\frac{1+d_px^4}{1-b_px^2-d_px^4}>1,
\]
where we use the monotonicity of $\frac{W}{P}$. In particular, $W<xW'$.
We also have
\begin{equation}\label{cert:eq:energy-sign}
 \begin{gathered}
 \left[\frac{1}{2}E(W')^2+c\left(\frac{1}{6}W^6-\frac{1}{2}W^2\right)\right]'
 =x(W')^2,\\
 E(W')^2+c(W^6-W^2)>\frac{2c}{3}W^6>0,
 \end{gathered}
\end{equation}
because the expression in brackets equals $\frac{1}{2}W'(0)^2>0$ at zero.

We take the solution regular at $x=1$ of
\[
 L_0\phi=-(E\phi')'-V\phi=0,\qquad \phi(1)=1,
\]
and set
\[
 \begin{gathered}
 g=c(W^5-W),\qquad \Phi=\frac{x}{2}(1-W^{-4}),\\
 A=E(W\phi'-W'\phi),\qquad B=E(EW'\phi'+g\phi).
 \end{gathered}
\]
By the profile equation and $W<xW'$, we have
\begin{equation}\label{cert:eq:wronskians}
 \begin{gathered}
 A'=-4cW^5\phi,\qquad B'=\Phi A',\qquad A(1)=B(1)=0,\\
 \Phi'=\frac{1}{2}+W^{-4}\left(\frac{2xW'}{W}-\frac{1}{2}\right)>0.
 \end{gathered}
\end{equation}
We first prove $A>0$ on $(0,1)$. This holds near one. At a first
zero $s$ moving left, we would have
\[
 A(s)=0,\qquad A'(s)\ge0,\qquad \phi(s)<0,\qquad
 \frac{\phi'(s)}{\phi(s)}=\frac{W'(s)}{W(s)}.
\]
If $\phi(s)=0$, then $A(s)=0$ forces $\phi'(s)=0$, contrary to
$\phi(1)=1$. Since $(B-\Phi A)'=-\Phi'A$, we can evaluate $B(s)$
in two ways. By \eqref{cert:eq:energy-sign} and \eqref{cert:eq:wronskians}, we obtain
\[
 0>\frac{E(s)\phi(s)}{W(s)}
       \bigl[E(W')^2+Wg\bigr](s)=B(s)=\int_s^1\Phi'A\dd x>0.
\]
Thus $A>0$, and we have
\[
 \left(\frac{\phi}{W}\right)'=\frac{A}{EW^2}>0.
\]
Let $W(x_1)=1$. If $\phi(x_1)\ge0$, then $\phi>0$ on $(x_1,1)$.
We use
\[
 g(x_1)=0,\qquad g'(x_1)>0,\qquad L_0g=-20cEW^3(W')^2
\]
in Green's identity and reach the contradiction
\[
 0\le E(x_1)g'(x_1)\phi(x_1)
   =-20c\int_{x_1}^1EW^3(W')^2\phi\dd x<0.
\]
Thus $\phi<0$ and $g\phi>0$ on $(0,x_1)$.
Using $\Phi(x_1)=0$ and $B'=-2xg\phi$, we integrate separately
on $(0,x_1)$ and $(x_1,1)$ to obtain
\begin{equation}\label{cert:eq:endpoint-slope}
 \begin{aligned}
 B(0)&=\int_{x_1}^1\Phi'A\dd x+2\int_0^{x_1}xg\phi\dd x>0,\\
 \phi'(0)&=\frac{B(0)}{W'(0)}>0.
 \end{aligned}
\end{equation}
Since $\phi<0$ near zero and $\phi'(0)>0$, we have $\phi(0)<0$.
By endpoint uniqueness and the evenness of $V$, a radial kernel vector
must be even or odd. The even condition $\phi'(0)=0$ and the odd
condition $\phi(0)=0$ are both excluded.

\subsection{The full and double-even kernel counts}\label{cert:count-section}
In the coordinates of \eqref{geo:eq:spherical-profile}, we set
$x=\cos(2\chi)$ and divide the spherical linearization by four. In
angular mode $(j,\ell)$ we obtain
\begin{equation}\label{cert:eq:separated}
 L_{j\ell}=L_0+\frac{j^2}{2(1+x)}+\frac{\ell^2}{2(1-x)},
 \qquad j,\ell=0,1,2,\ldots.
\end{equation}
At each endpoint, the regular Frobenius branch has the corresponding
nonnegative angular power. The singular or logarithmic branch is excluded.
For the mixed mode, we have
\[
 L_{11}(\sqrt E\,W')=0,\qquad \sqrt E\,W'>0\quad(-1<x<1).
\]
By the ground-state transform, $L_{11}$ is nonnegative, with a
one-dimensional radial zero space. For $j,\ell\ge1$ and
$(j,\ell)\ne(1,1)$, we have
\[
 L_{j\ell}-L_{11}
 =\frac{j^2-1}{2(1+x)}+\frac{\ell^2-1}{2(1-x)}>0.
\]
Only $(1,1)$ contributes mixed modes, with four real directions.

For the pure degree-one mode, we have
\begin{equation}\label{cert:eq:pure-mode}
 \zeta(x)=\sqrt{1-x}\left[2(1+x)W'(x)+\frac{W(x)}{2}\right],
 \qquad (L_0+\rho)\zeta=0.
\end{equation}
This function is positive on $[0,1)$. At $x=-t$, its zeros are those
of $4(1-t)W'(t)-W(t)$, which is positive at zero and negative at one.
At each interior zero, we compute from \eqref{cert:eq:profile}
\[
 \bigl(4(1-t)W'-W\bigr)'
 =-\frac{W}{2(1-t^2)}-\frac{3W^5}{4(1+t)}<0.
\]
Every crossing is downward, so $\zeta$ has exactly one interior zero.
Endpoint regularity determines the radial solution up to a scalar.
Thus each of $(0,1)$ and $(1,0)$ contributes two real directions.

Let $x_0\in(-1,0)$ be the zero of $\zeta$. For $\Lambda>1$ and nonzero
$f$ in the form domain with $f(x_0)=0$, we apply the ground-state transform
on the two nodal intervals and pass to the form closure to obtain
\begin{equation}\label{cert:eq:nodal-positive}
 \mathfrak q_\Lambda[f]
 =\int_{-1}^1\left[
 E\zeta^2\left|\left(\frac{f}{\zeta}\right)'\right|^2
 +(\Lambda-1)\rho|f|^2\right]\dd x>0.
\end{equation}
To justify the closure at $x_0$, we use the zero trace and local $H^1$
regularity to estimate
\begin{equation}\label{cert:eq:nodal-cutoff}
 \frac{1}{r^2}\int_{|x-x_0|<2r}|f(x)|^2\dd x
 \lesssim\int_{|x-x_0|<2r}|f'(x)|^2\dd x\longrightarrow0.
\end{equation}
Cutting off $f$ near $x_0$ therefore converges in the form norm,
which justifies \eqref{cert:eq:nodal-positive} on the two nodal intervals.
For $1<\Lambda\le64$, the trial function in \eqref{cert:eq:trial}
satisfies
\[
 \mathfrak q_\Lambda[h]\le\mathfrak q_{64}[h]<0.
\]
If a nonzero kernel vector $u$ existed, then $h,u$ would be independent.
By form-orthogonality, we would have
\[
 \mathfrak q_\Lambda[\alpha h+\beta u]
 =|\alpha|^2\mathfrak q_\Lambda[h]\le0.
\]
A nonzero element of this two-dimensional span vanishes at $x_0$,
contrary to \eqref{cert:eq:nodal-positive}.
For $\Lambda\ge81$, \eqref{cert:eq:positive-form} excludes a zero mode.
Every integer $\ell\ge2$ satisfies $\ell^2\le64$ or $\ell^2\ge81$.
Consequently
\begin{equation}\label{cert:eq:angular-exclusion}
 \ker(L_0+\ell^2\rho)=\{0\},\qquad \ell=2,3,\ldots.
\end{equation}
Together with the radial exclusion, we have
\[
 \dim\ker\left(-\Delta_{\Sph^3}+\frac{3}{4}-5q^4\right)
 =\underbrace{2}_{(1,0)}+\underbrace{2}_{(0,1)}
   +\underbrace{4}_{(1,1)}=8.
\]
The pure degree-one vectors are conformal boosts, and the four mixed
vectors are rotations between the two coordinate planes.

By conformal covariance, the Euclidean energy kernel corresponds to the
spherical $H^1$ kernel. The pole has zero $H^1$ capacity, so it is
removable, and elliptic regularity makes each kernel vector smooth.
The sector even under both imaginary-coordinate reflections contains
one real direction from each of $(1,0)$, $(0,1)$ and $(1,1)$.
For a spherical kernel vector $v$, stereographic projection gives
\[
 \left(\frac{2}{1+|x|^2}\right)^{\frac{1}{2}}v(\sigma(x))
 =\frac{\sqrt2\,v(p)}{|x|}+O(|x|^{-2}).
\]
Its Euclidean representative belongs to $L^2$ exactly when $v(p)=0$.
Evaluation at $p$ has rank one in the double-even sector, since the
dilation tangent has nonzero monopole coefficient. Consequently
\[
 \dim\ker_{\dot H^1,\mathrm{even}}L=3,\qquad
 \dim\ker_{L^2,\mathrm{even}}L=3-1=2.
\]
This completes the proof of Theorem~\ref{cert:thm:gap}.

\Needspace{30\baselineskip}
\section{Notation for the construction and the single-core theorem}\label{notation:start}
The table lists the notation for the construction and for
Theorem~\ref{cap:main}, together with the sections or equations where
we define it.
\begingroup
\small
\begin{longtable}{@{}L{2.7cm}L{8.5cm}L{3.5cm}@{}}
\toprule
Notation & Definition or role & Location\\
\midrule
\endfirsthead
\toprule Notation & Definition or role & Location\\\midrule
\endhead
\bottomrule\endlastfoot
$\mathcal H,\mathrm{Sim}$ & Energy space and spatial similarities. & Section~\ref{main:start}\\
$Q_p,Z_p,E_p,R_p$ & Stationary chart, its tangent frame, finite tangents and monopole resonance. & Section~\ref{geo:validation-section}\\
$L_*,\mathcal A_*,\mathcal K_*$ & Frozen linearization, threshold residue and decaying Green remainder. & \eqref{fr:eq:main}\\
$X_2,Y_2,X,Y,\mathcal O$ & Weighted field, source and observation spaces. & \eqref{fr:eq:weighted-spaces}, \eqref{lin:eq:spaces}\\
$\Pi,\mu$ & Compatible source projection and compact dual tests. & \eqref{lin:eq:projection}\\
$\mathsf M,A$ & Frozen mass matrix and the chart representation of the simple-pole residue. & \eqref{construction:eq:frozen-frame}, \eqref{non:eq:residue-coordinates}\\
$\eta_j,T_j,I_j,e_j$ & Packet amplitudes, geometric blocks and their right endpoints. & Section~\ref{construction:controls-section}\\
$\mathcal T$ & Time norm for the field, source and velocity derivative. & \eqref{time:eq:norms}\\
$v_c,P(c),e_c$ & Packet train, its quadratic force and the integrated impulse error. & \eqref{cat:eq:train}, \eqref{cat:eq:integrated-error}\\
$\bar y,b,r,p$ & Reference polygon, velocity correction, corrected velocity and actual chart path. & \eqref{non:eq:path}\\
$f,F$ & Compatible source and full frozen source $F=f+\mu\mathsf M r'$. & \eqref{non:eq:wave-definition}\\
$H,H_*,w$ & Retarded core correction, its frozen counterpart and the scattering field. & Section~\ref{non:retarded-section}\\
$\mathcal B,\mathcal Y,\|x\|$ & Field, source and coupled trial norms, with $x=(f,b)$. & \eqref{non:eq:contraction-spaces}, \eqref{non:eq:trial-ball}\\
$S$ & Nonlinear source, with $\square u-u^5=F-S$. & \eqref{end:eq:residual}\\
$v_+$ & Asymptotic free wave of the counterexample. & Section~\ref{end:global-section}\\
$\mathcal K_P,\mathcal M_\ell(P)$ & Compact orbit of stationary representatives and the associated traveling-wave family. & Section~\ref{main:classes}\\
$\mathcal J_\ell,\eta(t)$ & Action at velocity $\ell$ and bound for the radiation and action errors. & Section~\ref{cap:sec-radiation}\\
$\varepsilon,b,\tau$ & Norm of the modulated error, magnitude of its mode coefficients and rescaled time. & Section~\ref{cap:sec-modulation}\\
\end{longtable}
\endgroup
\Needspace{9\baselineskip}
\section*{AI disclosure}
In 2014, J.H. and X.C. organized a seminar on soliton resolution,
initially focused on the nonlinear Schr\"odinger equation. Inspired by the
literature discussed in \cite[Section~1.1]{CHS5D}, the authors conceived
the idea to consider whether radiation could drive a stationary core along
the conformal orbit of a sign changing state, producing recurrent changes
of shape modulo spatial similarities. We constructed a counterexample to
the fixed-profile formulation for the nonlinear Schr\"odinger equation,
which is available upon request. We then turned to the nonlinear wave
equation, where we expected radiation to have a weaker effect on the core.

We noticed that increasing the separation between paired impulses so that
progressively weaker radiation packets could produce fixed core
displacements, and developed the five-dimensional construction in
\cite{CHS5D} with AI assistance. The present paper implements this
mechanism in three dimensions. It combines the radiation-identification
and profile-collapse arguments of \cite{CHS5D} into a single energy-space
norm identity. Two inequivalent limiting profiles suffice for the
fixed-profile obstruction considered here.

For the revised three-dimensional construction,\footnote{In v1 version of
this paper, the authors came up with a construction which was using odd
symmetry to cancel the leading monopole contribution and explicitly
separating the remaining dipole response. That construction is now in
\cite{CHS3D}.} we decided to replace the odd-symmetric core used in v1
with a different stationary state that has even symmetry, which gave
another route to the counterexample and enabled a more analytic proof
identifying the full kernel of the linearized operator. A further
insight is that the two limiting shapes belong to the same conformal
class. This observation led the authors to propose a reformulation of
soliton resolution in which each core may change shape within a fixed
conformal class. We also made progress toward establishing this
formulation in the nonradial one-core setting.

The main ideas of this paper are due to the authors. We directed ChatGPT
Astra 6.2 to assist in deriving and writing out specific technical
estimates. We also used ChatGPT to help polish the language and check
the proof.

\end{document}